\documentclass[10pt]{amsart}
\usepackage{amsmath,amssymb, amsbsy}
\usepackage[dvips]{graphicx}
\usepackage{color}
\usepackage[active]{srcltx}
\usepackage{wrapfig}
\usepackage{graphicx}
\usepackage[usenames,dvipsnames]{pstricks}
\usepackage{epsfig}
\usepackage{pst-grad} 
\usepackage{pst-plot} 
\usepackage{epic}
\usepackage{hyperref}
\usepackage[normalem]{ulem}
\usepackage{color}
\newcommand{\dyle}{\displaystyle}

\newcommand{\R}{{\mathbb R}}
\newcommand{\ren}{\R^{N}}

\renewcommand{\l }{\lambda }

\newcommand{\iy}{\infty}

\newcommand{\inn}{\text{  in   }}

\newcommand{\dint}{\dyle\int}

\newcommand{\e }{\varepsilon}

\renewcommand{\ge }{\geqslant}
\renewcommand{\geq }{\geqslant}
\renewcommand{\le }{\leqslant}
\renewcommand{\leq }{\leqslant}

\newenvironment{pfn}[1]{\noindent{\it Proof of
        {#1}.\enspace}}{\hfill\qed\medskip}\newtheorem{Theorem}{Theorem}[section]

\newtheorem{Lemma}[Theorem]{Lemma}
\newtheorem{Proposition}[Theorem]{Proposition}
\theoremstyle{definition}
\newtheorem{Definition}[Theorem]{Definition}

\newtheorem{remark}[Theorem]{Remark}

\begin{document}
    \title[]
    {Existence and nonexistence results for a mixed\\ local--nonlocal operator with a Hardy potential}
    \thanks{The first author is partially supported by Project PID2023-150166NB-100, MINECO, Spain.}
    \author[B. Abdellaoui, K.Biroud, E.-H. Laamri, L. Montoro ]{Boumediene Abdellaoui, Kheireddine Biroud, El-Haj Laamri, Luigi Montoro}

    \address{\hbox{\parbox{5.7in}{\medskip\noindent{{ B. Abdellaoui: Laboratoire d'Analyse Nonlin\'eaire et Math\'ematiques
                        Appliqu\'ees, \hfill \break\indent Universit\'e AbouBekr Belkaid,
                        Tlemcen, \hfill\break\indent Tlemcen 13000, Algeria. }}}}}\email{{\tt boumediene.abdellaoui@inv.uam.es}}

    \address{\hbox{\parbox{5.7in}{\medskip\noindent{{K. Biroud: Laboratoire d'Analyse Nonlin\'eaire et Math\'ematiques
                        Appliqu\'ees, \hfill \break\indent Ecole Sup\'rieure de Management
                        de Tlemcen, \hfill \break\indent No. 01, Rue Barka Ahmed
                        Bouhannak Imama, \hfill\break\indent Tlemcen 13000, Algeria.    }}}}}\email{{\tt kh$_{-}$biroud@yahoo.fr}}

    \address{\hbox{\parbox{5.7in}{\medskip\noindent {{E.-H. Laamri: Institut Elie Cartan,\\ Universit\'{e} Lorraine,\\
                        B. P. 239, 54506 Vand{\oe}uvre l\'es
                        Nancy, France. \\[3pt]
                        \em{E-mail addresses: }{\tt el-haj.laamri@univ lorraine.fr}.}}}}}

    \address{\hbox{\parbox{5.7in}{\medskip\noindent{L. Montoro: Dipartimento di Matematica e Informatica, Universit\`a  della Calabria,
                    \hfill \break\indent Ponte Pietro Bucci 31 B, 87036 Arcavacata di
                    Rende, Cosenza, Italy.}}}} \email{{\tt luigi.montoro@unical.it}}


    \begin{abstract} In this paper, we investigate a class of problems involving a mixed local–nonlocal operator with a Hardy-type potential. More precisely, we consider a problem driven by the operator \(-\Delta+(-\Delta)^s\), where \(s\in(0,1)\), in a bounded smooth domain \(\Omega\subset\mathbb{R}^N\) containing the origin, with \(N\geq 3\). The problem involves a positive parameter \(\lambda\) and a nonnegative measurable function \(f\) satisfying suitable assumptions.
            \\
            We first investigate the associated linear problem with a prescribed datum. By means of a suitable change of variables and within an appropriate framework of weighted Sobolev spaces, we establish a weak Harnack inequality for the resulting weighted mixed operator. As a consequence, we derive the precise asymptotic behavior of nonnegative solutions near the singularity at the origin. Our analysis also yields a necessary and sufficient condition for the existence of weak solutions.
            \\
            We next address the nonlinear case, focusing on two classes of nonlinearities:
            \[
            f(x,u)=u^q+g(x)
            \qquad\text{and}\qquad
            f(x,u)=\frac{h(x)}{u^\sigma}.
            \]
            For the first class, we identify a critical exponent governing the existence of solutions. For the second class, we establish existence results for a broad class of data $h$.
    \end{abstract}

 \maketitle

    \section{Introduction and statement of the problem}

    This work deals with the following problem:
    \begin{equation}\label{tlem}
        \left\{
        \begin{array}{rclll}
            -\Delta u+(-\Delta)^s u
            &=&
            \lambda \dfrac{u}{|x|^{2}}+f(x,u)
            &\hbox{ in } &\Omega,\\[4pt]
            u&>&0
            &\hbox{ in }&\Omega,
            \\[4pt]
            u&=&0
            &\hbox{ in }& \mathbb{R}^N\setminus\Omega,
        \end{array}
        \right.
    \end{equation}
    where $s\in(0,1)$, $N\geq 3$, $\lambda>0$, $\Omega\subset \mathbb{R}^N$ is a bounded smooth domain containing the origin, and $f$ is a measurable function satisfying suitable assumptions to be specified below.

    Throughout this paper, $(-\Delta)^s$ denotes the integral (or Riesz) fractional Laplacian, defined for any $u\in\mathcal{S}(\mathbb{R}^N)$ by
    \begin{equation}\label{fraclap}
        (-\Delta)^s u(x):=a_{N,s}\,\mathrm{P.V.}\int_{\mathbb{R}^N}\frac{u(x)-u(y)}{|x-y|^{N+2s}}\,dy,
    \end{equation}
    where $a_{N,s}$ is the normalization constant
    \begin{equation}\label{aNs}
        a_{N,s}:=\frac{s\,2^{2s}\,\Gamma\!\left(\frac{N+2s}{2}\right)}{\pi^{\frac{N}{2}}\,\Gamma(1-s)}.
    \end{equation}
    Equivalently, $(-\Delta)^s$ is the Fourier multiplier with symbol $|\xi|^{2s}$. The choice of $a_{N,s}$ makes the definition consistent with the limiting cases $s=0$ and $s=1$, in the sense that
    \[
    \lim_{s\to 0^+}(-\Delta)^s u=u
    \qquad\text{and}\qquad
    \lim_{s\to 1^-}(-\Delta)^s u=-\Delta u
    \]
    (see, for instance, \cite[Proposition~4.4]{dine} and \cite[Proposition~2.1]{DaoLaam}).\\
    The operator $(-\Delta)^s$ is the infinitesimal generator of a symmetric $2s$-stable L\'evy process, whose trajectories may jump outside $\Omega$ without touching its boundary. For this reason, the homogeneous Dirichlet condition in \eqref{tlem} is prescribed on the whole complement $\mathbb{R}^N\setminus\Omega$, and not only on $\partial\Omega$; this corresponds to killing the process as soon as it leaves $\Omega$. This is in contrast with the spectral fractional Laplacian, which is defined through the spectral decomposition of the Dirichlet Laplacian in $\Omega$ and only involves the values of $u$ inside $\Omega$. For further background on fractional Laplacians in bounded domains, and in particular on the distinction between the spectral and integral fractional Laplacians, we refer the reader to \cite{DaoLaam} and the references therein.\

    \medskip
    For notational convenience, we write $\mathcal{L}_s:=-\Delta+(-\Delta)^s$ throughout the paper. In this mixed operator, the local part accounts for a Brownian diffusion and the nonlocal part for long-range jumps.

    \medskip
    Due to the presence of the Hardy potential, Problem~\eqref{tlem} is closely related to the following mixed Hardy constant:
    \[
    \overline{\Lambda}_{N,s}
    =
    \inf_{\phi\in C_0^\infty(\Omega)}
    \frac{
        \displaystyle
        \int_{\Omega}|\nabla \phi|^2\,dx
        +\frac{a_{N,s}}{2}
        \iint_{\mathbb{R}^N\times \mathbb{R}^N}
        \frac{|\phi(x)-\phi(y)|^2}{|x-y|^{N+2s}}\,dx\,dy
    }
    {
        \displaystyle
        \int_{\Omega}\frac{\phi^2}{|x|^2}\,dx
    }.
    \]
    As proved in \cite{Biagi1}, $\overline{\Lambda}_{N,s}$ is independent of $\Omega$ and $s$, is never attained, and coincides with the classical Hardy constant $\Lambda_N$. More generally, Hardy inequalities play an important role in many areas of mathematical analysis and mathematical physics, including the stability of relativistic matter; see \cite{FLS}.

    In the critical case \(\lambda=\Lambda_N\), the improved Hardy inequality established in \cite[Theorem~2.2]{VZ}
    (see also \cite{VZ} and \cite{ACP}) yields, for every \(\phi\in C_0^\infty(\Omega)\) and every \(p<2\),

    \begin{equation}\label{imp00}
        \begin{aligned}
            &\int_{\Omega} |\nabla \phi|^2\,dx
            +\frac{a_{N,s}}{2}
            \iint_{\mathbb{R}^N\times \mathbb{R}^N}
            \frac{|\phi(x)-\phi(y)|^2}{|x-y|^{N+2s}}\,dx\,dy
            -\Lambda_N\int_{\Omega}\frac{\phi^2}{|x|^2}\,dx\\
            &\qquad\geq C(\Omega,N)
            \left(\int_{\Omega}|\nabla \phi|^p\,dx\right)^{\frac{2}{p}}.
        \end{aligned}
    \end{equation}
    This estimate will be used repeatedly throughout the paper, especially when dealing with the critical case $\lambda=\Lambda_N$.

    \medskip
    As we shall see in Remark~\ref{rem:NNN}, Problem~\eqref{tlem} admits no positive supersolution when $\lambda>\Lambda_N$. Hence, throughout this paper, we assume that $0<\lambda\leq\Lambda_N$.

    \medskip
    When $f$ is independent of $u$ (i.e., $f(x,u)=f(x)$), Problem~\eqref{tlem} reduces to the linear problem
    \begin{equation}\label{prob00}
        \left\{
        \begin{array}{rclll}
            -\Delta u + (-\Delta)^s u&=&
            \lambda \dfrac{u}{|x|^{2}}+f(x)
            &\hbox{ in } & \Omega,\\[4pt]
            u&=&0
            &\hbox{ in } &\mathbb{R}^N\setminus\Omega,
        \end{array}
        \right.
    \end{equation}
    and our first objective is to determine the optimal summability properties of the solution $u$ in terms of the summability of the datum $f$ and the value of the parameter $\lambda$. Since Problem~\eqref{prob00} is linear, we may assume without loss of generality that $f$ is nonnegative, and we focus on positive solutions.

    The case $\lambda=0$ was studied in \cite{Anto,Arora,Biagi,Biagi01,Bir}, where Calder\'on--Zygmund-type regularity estimates were established for the mixed operator.

    The influence of the Hardy potential in the purely local case, i.e., for the problem
    \begin{equation}\label{prob00PP}
        -\Delta u=\lambda \dfrac{u}{|x|^{2}}+f(x)
        \quad\hbox{ in } \Omega,
    \end{equation}
    was studied in detail in \cite{BOP}. Looking for radial solutions of the form $w(x)=|x|^{-\gamma}$ to
    \[
    -\Delta w=\lambda\frac{w}{|x|^2}
    \quad \text {in } \mathbb{R}^N\setminus\{0\},
    \]
    one finds $w_i(x)=|x|^{-\gamma_i}$, $i=1,2$, where
    \begin{equation}\label{gamma1}
        \gamma_{1}=\frac{N-2}{2}-\sqrt{\Lambda_N-\lambda},
        \qquad
        \gamma_{2}=\frac{N-2}{2}+\sqrt{\Lambda_N-\lambda},
    \end{equation}
    are the two roots of the algebraic equation
    \begin{equation}\label{polyn}
        \beta^2-(N-2)\beta-\lambda=0.
    \end{equation}
    Consequently, for Problem~\eqref{prob00PP} with nonnegative datum $f$, a necessary and sufficient condition for the existence of a positive solution is the weighted integrability condition
    \[
    \int_{\Omega} f(x)\,|x|^{-\gamma_1}\,dx<\infty.
    \]
    Moreover, assume that $f\in L^m(\Omega)$ with $m>1$, and define
    \[
    \lambda(m):=\frac{N(m-1)(N-2m)}{m^2}.
    \]
    If $0<\lambda<\lambda(m)$, the solution of Problem~\eqref{prob00PP} satisfies the same Calder\'on--Zygmund estimates as in the case $\lambda=0$. In contrast, when $\lambda\geq\lambda(m)$, counterexamples can be constructed showing that these regularity properties fail (see \cite{BOP} for details). In particular, when $m>\frac{N}{2}$, solutions need not be bounded, whereas existence is not guaranteed in general when $m=1$.

    \medskip
    The situation for the mixed operator $\mathcal{L}_s$ is considerably more delicate. In order to extend the above results to this setting, we use a suitable ground-state transformation and establish a \emph{weak Harnack inequality} for a singular weighted mixed operator (introduced below). The proof, given in Section~\ref{sec-harnack}, combines several technical estimates, including potential-type estimates in weighted Sobolev spaces and sharp bounds for certain radial singular integrals. As a consequence, we derive the precise asymptotic behavior of positive supersolutions near the origin. This result is a cornerstone of the subsequent existence and regularity analysis.

    \medskip
    As a first application of the weak Harnack inequality, we return to the linear Problem~\eqref{prob00}. We derive a sharp upper bound for its solution near the origin, and we give a necessary and sufficient condition, in the form of a weighted integrability condition on $f$, for the existence of a weak solution, together with regularity properties of the minimal solution.

    \medskip
    As a second application, we consider the semilinear problem
    \begin{equation}\label{tlem001}
        \left\{
        \begin{array}{rclll}
            \mathcal{L}_{s}u
            &=&\lambda \dfrac{u}{|x|^{2}}
            +u^q+g(x)
            &\hbox{ in }  &\Omega,\\[4pt]
            u&>&0
            &\hbox{ in }  & \Omega,\\[4pt]
            u&=&0
            &\hbox{ in }  & \mathbb{R}^N\setminus\Omega.
        \end{array}
        \right.
    \end{equation}
    For the local operator $-\Delta$, it was shown in \cite{BDT} that the Hardy potential gives rise to a critical exponent $q^+(\lambda)$ that governs the existence of positive supersolutions. This result relies on the precise behavior of solutions near the origin, together with suitable comparison principles. Our aim is to identify the corresponding critical exponent for $\mathcal{L}_s$. More precisely, we prove that if $q\geq q^+(\lambda)$, Problem~\eqref{tlem001} admits no positive weak supersolution. In the supercritical case $q>q^+(\lambda)$, the proof relies on a suitable variant of Picone's inequality combined with the local behavior of any supersolution near the origin. The critical case $q=q^+(\lambda)$ requires a finer description of the behavior of supersolutions near the singularity, which is provided by the weak Harnack inequality. We also establish a complete blow-up result for solutions to approximating problems. Conversely, for $1<q<q^+(\lambda)$, we prove an existence result showing the optimality of $q^+(\lambda)$.

    \smallbreak
    \noindent We also extend a result of Brezis and Cabr\'e by proving a nonexistence theorem for the problem
    \[
    \mathcal{L}_{s}u = \dfrac{u^q}{|x|^{2}}
    \quad\hbox{ in } \Omega.
    \]
    Specifically, we show that this problem admits no positive weak supersolution for any $q>1$. We further discuss connections between Problem~\eqref{tlem001} and results of Baras--Pierre, Brezis--Cabr\'e and Kalton--Verbitsky type.

    \medskip
    Finally, as a third application, we study the singular semilinear problem
    \begin{equation}\label{tlem002}
        \left\{
        \begin{array}{rclll}
            \mathcal{L}_{s}u
            &=&\lambda \dfrac{u}{|x|^{2}}+\dfrac{h(x)}{u^\sigma}
            &\hbox{ in } & \Omega,\\[4pt]
            u&>&0
            &\hbox{ in } &\Omega,\\[4pt]
            u&=&0
            &\hbox{ in } & \mathbb{R}^N\setminus\Omega.
        \end{array}
        \right.
    \end{equation}
    The purely local case without Hardy potential, i.e., for the operator $-\Delta$, was studied in \cite{BO}. In that paper, the authors proved the existence of a distributional solution for every nonnegative $h\in L^1(\Omega)$, with additional regularity properties depending on the summability of $h$ and on the value of the singular exponent $\sigma$. The purely fractional case without Hardy potential, i.e., for the operator $(-\Delta)^s$, was investigated in \cite{BBMI,CMoSS}. The quasilinear case with a Hardy potential was treated in \cite{AA0}. For the mixed operator $\mathcal{L}_s$, we prove existence and regularity results depending on the value of $\sigma$ and the summability of $h$, and we show that, when $h$ is a measure concentrated on a set of zero capacity, the solutions of the approximating problems converge to zero.

    \medskip
    In summary, in all the problems considered in this paper, the simultaneous presence of the classical Laplacian, the fractional Laplacian and the Hardy potential gives rise to new analytical difficulties, particularly in the study of the behavior of solutions near the singularity and in the derivation of optimal regularity estimates. The weak Harnack inequality established in Section~\ref{sec-harnack} is the common tool that allows us to overcome these difficulties, since it provides the precise behavior of supersolutions near the origin.

    \medskip
    The paper is organized as follows. Section~\ref{prim} introduces the functional framework and the notions of solution used throughout the paper. It also presents several preliminary results, including Picone-type inequalities, compactness properties, and useful algebraic estimates.

    Section~\ref{sec-harnack} is devoted to the proof of the weighted weak Harnack inequality. The strategy is to reduce the problem, via a ground-state change of variables, to a weak Harnack inequality for a weighted mixed operator without a Hardy singularity, and then to establish this inequality by means of a Moser-type iteration. As a consequence, we derive the precise asymptotic behavior of any positive supersolution of Problem~\eqref{prob00} near the origin.

    Section~\ref{s4} is devoted to the three applications described above: the linear Problem~\eqref{prob00} in Subsection~\ref{subsec:linear}, Problem~\eqref{tlem001} in Subsection~\ref{subsec:application}, and the singular Problem~\eqref{tlem002} in Subsection~\ref{sec:main}.

    The paper concludes with a section discussing open problems and directions for future research.


    \section{Functional setting and useful tools}\label{prim}

    We begin this section by introducing the functional framework and the auxiliary tools that will be used throughout the paper. We first recall the definition of fractional Sobolev spaces.

    Let $s\in (0,1)$. For any $p\in [1,\infty)$ and $\Omega\subseteq\mathbb{R}^N$, we define
    \[
    W^{s,p}(\Omega) := \left\{ u\in L^{p}(\Omega) : \,
    \frac{|u(x)-u(y)|}{|x-y|^{\frac{N}{p}+s}} \in
    L^{p}(\Omega\times\Omega)\right\}.
    \]
    For $p=2$, the fractional Sobolev space $H^{s}(\Omega):=W^{s,2}(\Omega)$ is a Hilbert space. When $\Omega=\mathbb{R}^N$, the Fourier transform provides an equivalent characterization of this space.

    \begin{Definition}\label{Hs}
        For $0<s<1$, we define the fractional Sobolev space of order $s$ by
        \[
        H^s(\mathbb{R}^N) :=
        \bigl\{ u\in L^2(\mathbb{R}^N) :
        |\xi|^{s}\mathcal{F}(u)(\xi)\in L^2(\mathbb{R}^N) \bigr\}.
        \]
    \end{Definition}

    For $u\in H^s(\mathbb{R}^N)$, we denote by
    \[
    [u]^2_s:=\iint_{\mathbb{R}^N\times\mathbb{R}^N}
    \frac{|u(x)-u(y)|^2}{|x-y|^{N+2s}}\,dx\,dy
    \]
    the Gagliardo seminorm of $u$. By Plancherel's identity, we obtain the following equivalent formulation; see, e.g., \cite[Lemma~3.1]{FLS} or \cite[Proposition~3.4]{dine}.

    \begin{Proposition}\label{norma}
        Let $N\geq 1$ and $0<s<1$. Then, for every $u\in H^s(\mathbb{R}^N)$,
        \[
        \int_{\mathbb{R}^N} |\xi|^{2s}|\mathcal{F}(u)(\xi)|^2\,d\xi
        =\frac{a_{N,s}}{2}\,[u]^2_s,
        \]
        where $a_{N,s}$ is the constant defined in \eqref{aNs}.
    \end{Proposition}


    Now, for $\Omega\subset\mathbb{R}^N$, we define the space $H^s_0(\Omega)$ as the set of functions $u\in H^s(\mathbb{R}^N)$ such that $u=0$ a.e.\ in $\mathbb{R}^N\setminus\Omega$. For such functions,
    \begin{equation}\label{norm2}
        [u]^2_s
        =\iint_{D_\Omega}
        \frac{|u(x)-u(y)|^2}{|x-y|^{N+2s}}\,dx\,dy,
    \end{equation}
    where $D_\Omega := \mathbb{R}^N\times\mathbb{R}^N\setminus\bigl(\mathcal{C}\Omega\times\mathcal{C}\Omega\bigr)$ and $\mathcal{C}\Omega:= \mathbb{R}^N\setminus\Omega$.

    We define the mixed energy space $\mathcal{H}_0^1(\Omega)$ by
    \[
    \mathcal{H}_0^1(\Omega)
    :=\Bigl\{ u\in H^1(\mathbb{R}^N) :
    u\big|_\Omega \in H^1_0(\Omega),\quad
    u=0 \text{ a.e.\ in } \mathbb{R}^N\setminus\Omega \Bigr\}.
    \]
    The space $\mathcal{H}_0^1(\Omega)$ is a Hilbert space with respect to the norm
    \[
    \|u\|_{\mathcal{H}_0^1(\Omega)}
    :=\Bigl( \|\nabla u\|^2_{L^2(\Omega)}
    + \frac{a_{N,s}}{2}\,[u]^2_s \Bigr)^{1/2},
    \]
    where $[u]_s$ is given by \eqref{norm2}. We also introduce the local variant
    \[
    \mathcal{H}_{loc}^1(\Omega)
    :=\Bigl\{ u:\mathbb{R}^N\to\mathbb{R} :
    |\nabla u|\in L^2(K),\;
    [u]^2_{s,K}<\infty \;\text{ for every } K\subset\subset\Omega \Bigr\},
    \]
    where
    \[
    [u]^2_{s,K}
    := \iint_{K\times K}
    \frac{|u(x)-u(y)|^2}{|x-y|^{N+2s}}\,dx\,dy.
    \]
    Here $K\subset\subset\Omega$ means that $K$ is compactly contained in $\Omega$, so that $K\times K\subset D_\Omega$.


    \medskip
    Before turning to the problem with the Hardy potential, we make precise the notion of solution for the mixed operator $\mathcal{L}_{s}$ with general data. Consider the model problem
    \begin{equation}\label{ENT00}
        \left\{
        \begin{array}{rcll}
            \mathcal{L}_{s}u &=& f &\text{in }\Omega,\\
            u &=& 0 &\text{in }\mathbb{R}^N\setminus\Omega.
        \end{array}
        \right.
    \end{equation}
    For $f\in\mathcal{H}^{-1}(\Omega)$, where $\mathcal{H}^{-1}(\Omega)$ denotes the dual space of $\mathcal{H}_0^1(\Omega)$, the Lax--Milgram theorem yields the existence of a unique finite energy solution to \eqref{ENT00}.

    This finite energy framework is, however, not sufficient to cover the general data considered in this paper. We therefore introduce a weaker notion of solution, suited to data of much lower regularity. For $k>0$, we recall the truncation operator
    \begin{equation}\label{f-trun00}
        T_k(a)=\left\{
        \begin{array}{cl}
            a, &\text{if }|a|\leq k,\\[2mm]
            k\,\dfrac{a}{|a|}, &\text{if }|a|>k,
        \end{array}
        \right.
    \end{equation}
    and we set $G_k(a):=a-T_k(a)$.

    \begin{Definition}\label{def:entropy}
        Let $f\in L^{1}(\Omega)$. We say that $u\in L^1(\Omega)$, with $u=0$ a.e.\ in $\mathbb{R}^N\setminus\Omega$, is an \emph{entropy solution} to Problem~\eqref{ENT00} if $T_k(u)\in\mathcal{H}_0^1(\Omega)$ for all $k>0$ and, for every $\varphi\in\mathcal{H}_0^1(\Omega)\cap L^\infty(\Omega)$ and every $k>0$,
        \begin{equation}\label{eq:alcala}
            \begin{split}
                &\int_\Omega \nabla u\cdot\nabla T_k(u-\varphi)\,dx
                +\frac{a_{N,s}}{2}\iint_{D_\Omega}
                \frac{(u(x)-u(y))
                    [T_k(u(x)-\varphi(x))-T_k(u(y)-\varphi(y))]}{|x-y|^{N+2s}}\,dx\,dy \\
                &\qquad\leq
                \int_\Omega f(x)\,T_k(u(x)-\varphi(x))\,dx.
            \end{split}
        \end{equation}
    \end{Definition}

    Throughout the paper, we work with two notions of solution: \emph{finite energy solutions}, associated with data in the dual of the energy space, and \emph{weak solutions}, namely entropy solutions or solutions obtained as limits of approximations, following the standard approach of \cite{BZ,LPPS,AAB}. In this framework, one has the following existence and uniqueness result.

    \begin{Theorem}\label{SOLA1}
        Assume that $f\in L^1(\Omega)$. Then, Problem~\eqref{ENT00} has a unique entropy solution $u$, which is also a solution obtained as limit of approximations, and $u\in W^{1,q}_0(\Omega)$ for every $q<\frac{N}{N-1}$. Moreover, $T_k(u)\in\mathcal{H}_0^1(\Omega)$ for all $k>0$, and, for every $\varphi\in C_0^\infty(\Omega)$,
        \begin{equation}\label{test3}
            \begin{split}
                &\int_\Omega \nabla u\cdot\nabla\varphi\,dx
                +\frac{a_{N,s}}{2}\iint_{D_\Omega}
                \frac{(u(x)-u(y))(\varphi(x)-\varphi(y))}{|x-y|^{N+2s}}\,dx\,dy \\
                &\qquad =\int_\Omega f(x)\,\varphi(x)\,dx.
            \end{split}
        \end{equation}
    \end{Theorem}

    \medskip
    Denoting by $\mathcal{G}_{N,s}$ the Green function of the mixed operator $\mathcal{L}_{s}$ in $\Omega$, the solution to \eqref{ENT00} may be represented as
    \[
    u(x)=\int_{\Omega}\mathcal{G}_{N,s}(x,y)\,f(y)\,dy.
    \]
    It is known that $\mathcal{G}_{N,s}\simeq\mathcal{G}_N$, where $\mathcal{G}_N$ is the Green function of the Laplacian in $\Omega$; see, e.g., \cite{Arora}. This equivalence yields the following useful proposition.

   \begin{Proposition}\label{BC}
        Let $f\in L^1(\Omega)$ be a nonnegative function and let $u$ be the unique entropy solution to \eqref{ENT00}. Then, there exists a positive constant $C=C(N,s,\Omega)$ such that, for a.e.\ $x\in\Omega$,
        \[
        u(x)\geq C\,\delta(x)\int_{\Omega}f(y)\,\delta(y)\,dy,
        \]
        where $\delta(x):=\operatorname{dist}(x,\partial\Omega)$.
    \end{Proposition}

    \begin{proof}
        Let $\tilde u$ be the solution to $-\Delta\tilde u=f$ in $\Omega$, $\tilde u=0$ on $\partial\Omega$. Since $\mathcal{G}_{N,s}\geq c\,\mathcal{G}_N$ for some $c>0$, we have
        \[
        u(x)=\int_{\Omega}\mathcal{G}_{N,s}(x,y)f(y)\,dy
        \geq c\int_{\Omega}\mathcal{G}_{N}(x,y)f(y)\,dy
        =c\,\tilde u(x).
        \]
        The conclusion then follows from the corresponding estimate for the Laplacian, due to Brezis and Cabr\'e \cite[Lemma~3.2]{BC},  applied to $\tilde u$.
    \end{proof}

    \begin{remark}\label{BBB}
        If, in addition, $f\in L^\infty(\Omega)$,  the properties of $\mathcal{G}_{N,s}$ yield $u(x)\simeq\delta(x)$ for a.e.\ $x\in\Omega$.
    \end{remark}


    \medskip
    We now turn to the operator with the Hardy potential, namely Problem~\eqref{prob00}. We introduce the notion of \emph{finite energy solution}, which is appropriate for data in the dual of the energy space and allows one to work within the natural variational framework.

    \begin{Definition}\label{energia}
        Assume $0<\lambda<\Lambda_N$ and let $f\in\mathcal{H}^{-1}(\Omega)$. We say that $u\in\mathcal{H}_0^1(\Omega)$ is a \emph{finite energy solution} to \eqref{prob00} if
        \[
        \int_{\Omega}\nabla u\cdot\nabla\phi\,dx
        +\frac{a_{N,s}}{2}
        \iint_{D_\Omega}
        \frac{(u(x)-u(y))(\phi(x)-\phi(y))}{|x-y|^{N+2s}}\,dx\,dy
        -\lambda\int_{\Omega}\frac{u\,\phi}{|x|^2}\,dx
        =\langle f,\phi\rangle
        \]
        for all $\phi\in\mathcal{H}_0^1(\Omega)$.
    \end{Definition}

    If $\lambda<\Lambda_{N}$, the coercivity of the associated bilinear form, together with the Lax--Milgram theorem, yields the existence and uniqueness of a finite energy solution $u\in\mathcal{H}_0^1(\Omega)$ for every $f\in\mathcal{H}^{-1}(\Omega)$.

    \begin{remark}\label{HHH}
        For $\lambda=\Lambda_N$, existence and uniqueness still hold, but the natural solution space is enlarged to
        \[
        \mathcal{H}(\Omega)
        :=\overline{C^\infty_0(\Omega)}^{\|\cdot\|_{\mathcal{H}(\Omega)}},
        \]
        where
        \[
        \|\phi\|^2_{\mathcal{H}(\Omega)}
        :=\int_{\Omega}|\nabla\phi|^2\,dx
        +\frac{a_{N,s}}{2}
        \iint_{D_\Omega}
        \frac{|\phi(x)-\phi(y)|^2}{|x-y|^{N+2s}}\,dx\,dy
        -\Lambda_N\int_{\Omega}\frac{\phi^2}{|x|^2}\,dx.
        \]
        By the improved Hardy inequality \eqref{imp00}, we have $\mathcal{H}(\Omega)\subset W^{1,p}_0(\Omega)$ for all $p<2$.
    \end{remark}

    Consider the space of test functions
    \begin{equation}\label{test}
        \mathcal{X}:=\Bigl\{\phi:\mathbb{R}^N\to\mathbb{R} :
        \mathcal{L}_{s}\phi=\vartheta \text{ in }\Omega \text{ for some } \vartheta\in L^\infty(\Omega)\cap C^\alpha(\Omega),\;
        0<\alpha<1,\;
        \phi=0 \text{ in }\mathbb{R}^N\setminus\Omega\Bigr\}.
    \end{equation}
    Following the same approach as in \cite{LPPS}, one deduces that $\mathcal{X}\subset L^\infty(\Omega)$ and that each $\phi\in\mathcal{X}$ is a classical solution of the corresponding problem.
    \begin{Definition}\label{veryweak}
        Let $f\in L^{1}(\Omega)$ be a nonnegative function. Suppose that $u\in L^{1}(\Omega)$ is nonnegative, $\frac{u}{|x|^2}\in L^1(\Omega)$, and $u=0$ in $\mathbb{R}^N\setminus\Omega$. We say that $u$ is a \emph{weak supersolution} (resp.\ \emph{subsolution}) of Problem~\eqref{prob00} if, for all nonnegative $\phi\in\mathcal{X}$,
        \[
        \int_{\Omega}u\,\mathcal{L}_{s}\phi\,dx
        -\lambda\int_{\Omega}\frac{u\,\phi}{|x|^{2}}\,dx
        \geq\;(\leq)\;
        \int_{\Omega}f\,\phi\,dx.
        \]
        If $u$ is both a supersolution and a subsolution, we call it a \emph{weak solution}.
    \end{Definition}

    Under the above conditions, if $u$ is a nonnegative supersolution of \eqref{prob00} in the sense of Definition~\ref{veryweak}, then
    \[
    u(x)\geq\int_{\Omega}\mathcal{G}_{N,s}(x,y)
    \Bigl(\lambda\,\frac{u(y)}{|y|^2}+f(y)\Bigr)\,dy.
    \]
    The following proposition justifies the direct use of the weak formulation with test functions, and will be applied systematically throughout the paper.

    \begin{Proposition}\label{aprggg}
        Assume that $u$ is a nonnegative function with $\dfrac{u}{|x|^2}\in L^1(\Omega)$ and $u=0$ in $\mathbb{R}^N\setminus\Omega$. Suppose that
        \begin{equation}\label{super10}
            \mathcal{L}_{s}u-\lambda\,\frac{u}{|x|^{2}}\geq 0
            \quad\text{in }\Omega
        \end{equation}
        in the sense of Definition~\ref{veryweak}. Then there exists $w\in W^{1,q}_0(\Omega)$ for every $q<\frac{N}{N-1}$ such that
        \begin{equation}\label{sola00}
            \left\{
            \begin{array}{rclll}
                \mathcal{L}_{s}w &=& \eta &\text{\rm in }&\Omega,\\
                w &>& 0 &\text{\rm in } &\Omega,\\
                w &=& 0 &\text{\rm in } &\mathbb{R}^N\setminus\Omega,
            \end{array}
            \right.
        \end{equation}
        where $\eta\in L^1(\Omega)$ satisfies $\eta\geq\lambda\,\dfrac{w}{|x|^2}$ a.e.\ in $\Omega$. In particular, $w$ is a supersolution of \eqref{super10} with $w\leq u$.
    \end{Proposition}

    \begin{proof}
        For $n\in\mathbb{N}^*$, let $w_n$ be the unique solution to
        \begin{equation}\label{solannn}
            \left\{
            \begin{array}{rclll}
                \mathcal{L}_{s}w_n &=&
                T_n\!\left(\lambda\,\dfrac{u}{|x|^2}\right)
                &\text{\rm in }&\Omega,\\[4pt]
                w_n &=& 0 &\text{\rm in } &\mathbb{R}^N\setminus\Omega.
            \end{array}
            \right.
        \end{equation}
        The existence of $w_n$ follows from the Lax--Milgram theorem; moreover $w_n\in C^{1,\alpha}(\overline{\Omega})$ (see for instance~\cite[Theorem~1.1]{Anto}). The sequence $\{w_n\}_n$ is increasing in $n$. Since $T_n(\lambda\,u/|x|^2)\leq\lambda\,u/|x|^2$ for all $n$, the comparison principle gives $w_n\leq u$, and hence $w_n/|x|^2\leq u/|x|^2$. By Theorem~\ref{SOLA1}, $w_n$ converges to the entropy solution $w$ of \eqref{sola00} with $\eta=\lambda\,u/|x|^2$. Since $\lambda\,u/|x|^2\geq\lambda\,w/|x|^2$, we conclude that $w$ is the desired supersolution of \eqref{super10}, obtained as a limit of approximations.
    \end{proof}


    \medskip
    Consider a solution $u$ of Problem~\eqref{tlem} with $\lambda\leq\Lambda_N$, and define $w := |x|^{\gamma_1}u$, where $\gamma_1$ is given by \eqref{gamma1}. Then $w$ solves
    \begin{equation}\label{equi}
        \left\{
        \begin{array}{rcll}
            L_{\gamma_1}w &=&
            |x|^{-\gamma_1}f(x,|x|^{-\gamma_1}w)=:g(x,w)
            &\text{in }\Omega,\\[2pt]
            w &=& 0 &\text{in }\mathbb{R}^N\setminus\Omega,
        \end{array}
        \right.
    \end{equation}
    where
    \begin{equation}\label{LgO}
        L_{\gamma_1}w(x)
        :=-\operatorname{div}\bigl(|x|^{-2\gamma_1}\nabla w\bigr)
        +\frac{a_{N,s}}{|x|^{\gamma_1}}\,\mathrm{P.V.}
        \int_{\mathbb{R}^N}
        \frac{|x|^{-\gamma_1}w(x)-|y|^{-\gamma_1}w(y)}{|x-y|^{N+2s}}\,dy.
    \end{equation}

    To analyze the regularity of $u$ through the auxiliary function $w$, we work in weighted Sobolev spaces with admissible weights. Throughout this work, we set
    \[
    d\mu := \frac{dx}{|x|^{2\gamma_1}},\qquad
    d\nu := \frac{dx\,dy}{|x-y|^{N+2s}|x|^{\gamma_1}|y|^{\gamma_1}},\qquad
    d\tau := \frac{dx\,dy}{|x-y|^{N+2s}}.
    \]
    For $\Omega\subseteq\mathbb{R}^N$, we define the weighted Sobolev space
    \[
    Y^{1,\gamma_1}(\Omega)
    := \Bigl\{ \phi\in L^2(\Omega,d\mu) :
    \|\phi\|_{Y^{1,\gamma_1}(\Omega)}<\infty \Bigr\},
    \]
    where
    \begin{equation}\label{norm1}
        \|\phi\|_{Y^{1,\gamma_1}(\Omega)}
        :=\Bigl(
        \int_{\Omega}(|\phi|^2+|\nabla\phi|^2)\,d\mu
        +\iint_{\Omega\times \Omega}
        (|x|^{-\gamma_1}\phi(x)-|y|^{-\gamma_1}\phi(y))^2\,d\tau
        \Bigr)^{1/2}.
    \end{equation}
    Equipped with the inner product induced by \eqref{norm1}, $Y^{1,\gamma_1}(\Omega)$ is a Hilbert space. Similarly, we define
    \[
    Y^{s,\gamma_1}(\Omega)
    := \Bigl\{ \phi\in L^2(\Omega,d\mu) :
    \|\phi\|_{Y^{s,\gamma_1}(\Omega)}<\infty \Bigr\},
    \]
    with
    \begin{equation}\label{norm1ss}
        \|\phi\|_{Y^{s,\gamma_1}(\Omega)}
        :=\Bigl(
        \int_{\Omega}|\phi|^2\,d\mu
        +\iint_{\Omega\times \Omega}
        (|x|^{-\gamma_1}\phi(x)-|y|^{-\gamma_1}\phi(y))^2\,d\tau
        \Bigr)^{1/2},
    \end{equation}
    which is also a Hilbert space.

    \begin{remark}\label{remmm}
        By \cite[Theorem~2.10]{AABP}, applied with $\beta=\gamma_1$ and $\alpha=-\gamma_1$, one has, for every $\phi\in C^\infty_0(\mathbb{R}^N)$,
        \[
        \iint_{\mathbb{R}^N\times \mathbb{R}^N}
        (|x|^{-\gamma_1}\phi(x)-|y|^{-\gamma_1}\phi(y))^2\,d\tau
        \simeq
        \iint_{\mathbb{R}^N\times \mathbb{R}^N}
        (\phi(x)-\phi(y))^2\,d\nu.
        \]
    \end{remark}

    The following extension lemma is proved by combining Remark~\ref{remmm} with the arguments of \cite{Adams} (see also \cite{dine}).

    \begin{Lemma}\label{ext}
        Let $\Omega\subset\mathbb{R}^N$ be a smooth domain. For every $w\in Y^{1,\gamma_1}(\Omega)$, there exists an extension $\tilde{w}\in Y^{1,\gamma_1}(\mathbb{R}^N)$ with $\tilde{w}|_\Omega=w$ and
        \[
        \|\tilde{w}\|_{Y^{1,\gamma_1}(\mathbb{R}^N)}
        \leq C\,\|w\|_{Y^{1,\gamma_1}(\Omega)},
        \]
        where $C=C(N,s,\Omega,\gamma_1)>0$.
    \end{Lemma}

    \begin{Lemma}\label{new010}
        Let $s\in(0,1)$ and let $\Omega\subset\mathbb{R}^N$ be a smooth bounded domain. {Assume that $0<\lambda<\Lambda_N$.} Then there exists a constant $C=C(N,s,\Omega,\gamma_1)$ such that
        \begin{equation}\label{ee2e}
            \iint_{\mathbb{R}^N\times \mathbb{R}^N}
            (|x|^{-\gamma_1}\phi(x)-|y|^{-\gamma_1}\phi(y))^2\,d\tau
            \leq C\int_{\Omega}|\nabla\phi|^2\,d\mu
        \end{equation}
        for every $\phi\in C^\infty_0(\Omega)$.
    \end{Lemma}

    \begin{proof}
        Let $\phi\in C^\infty_0(\Omega)$. By \cite[Lemma~2.1]{Butch},
        \begin{equation}\label{eee}
            \iint_{\mathbb{R}^N\times \mathbb{R}^N}
            \frac{(|x|^{-\gamma_1}\phi(x)-|y|^{-\gamma_1}\phi(y))^2}
            {|x-y|^{N+2s}}\,dx\,dy
            \leq C\int_{\Omega}|\nabla(|x|^{-\gamma_1}\phi)|^2\,dx.
        \end{equation}
        A direct computation yields
        \[
        |\nabla(|x|^{-\gamma_1}\phi)|^2
        \leq C\bigl(|x|^{-2\gamma_1}|\nabla\phi|^2
        +|x|^{-2\gamma_1-2}|\phi|^2\bigr).
        \]
        Applying the weighted Hardy inequality of \cite{ACS1}, we obtain
        \begin{equation}\label{eee1}
            \int_{\Omega}|\nabla(|x|^{-\gamma_1}\phi)|^2\,dx
            \leq C\int_{\Omega}|x|^{-2\gamma_1}|\nabla\phi|^2\,dx.
        \end{equation}
        Combining \eqref{eee} and \eqref{eee1} yields \eqref{ee2e}.
    \end{proof}
    \medbreak

    We define $Y^{1,\gamma_1}_0(\Omega)$ as the completion of $C^\infty_0(\Omega)$ with respect to the norm of $Y^{1,\gamma_1}(\Omega)$. On a bounded domain $\Omega$, by Lemma~\ref{new010} and the weighted Hardy inequality, the norm
    \begin{equation}\label{equiv-norm}
        \|\phi\|_{Y^{1,\gamma_1}_0(\Omega)}
        :=\Bigl(\int_{\Omega}|x|^{-2\gamma_1}|\nabla\phi|^2\,dx\Bigr)^{1/2}
    \end{equation}
    is equivalent to the norm of $Y^{1,\gamma_1}(\Omega)$ on $Y^{1,\gamma_1}_0(\Omega)$. It is also equivalent to the energy norm
    \begin{equation}\label{energy-norm}
        [\![\phi]\!]^2
        :=\int_{\Omega}|x|^{-2\gamma_1}|\nabla\phi|^2\,dx
        +\frac{a_{N,s}}{2}\iint_{D_\Omega}
        (|x|^{-\gamma_1}\phi(x)-|y|^{-\gamma_1}\phi(y))^2\,d\tau,
    \end{equation}
    which is associated with the operator $L_{\gamma_1}$.

    \medskip
    The following Caffarelli--Kohn--Nirenberg inequality is recalled from \cite{CKN}.

    \begin{Lemma}\label{caff}
        Let $p\geq 1$ and $\alpha,\gamma\in\mathbb{R}$ with $\alpha<N$ and $\gamma<(N-p)/p$. Set
        \[
        p^*_{\alpha,\gamma}
        :=\frac{(N-\alpha)p}{N-(\gamma+1)p}.
        \]
        There exists a positive constant $C$, depending only on the parameters and independent of $\Omega$, such that
        \begin{equation}\label{Caff}
            \Bigl(\int_{\Omega}|\phi|^{p^*_{\alpha,\gamma}}
            |x|^{-\alpha}\,dx\Bigr)^{1/p^*_{\alpha,\gamma}}
            \leq C\Bigl(\int_{\Omega}|\nabla\phi|^p
            |x|^{-\gamma p}\,dx\Bigr)^{1/p}
        \end{equation}
        for all $\phi\in C_0^\infty(\Omega)$ if and only if
        \[
        \alpha-p\leq\gamma p\leq\frac{(N-p)\alpha}{N}.
        \]
    \end{Lemma}

    We also record the following weighted Poincar\'e--Wirtinger inequality; see \cite[Appendix~B]{AMPP} for a proof.

    \begin{Theorem}\label{PW}
        Let $r>0$ and $w\in Y^{s,\gamma_1}(B_r)$. Assume that $\psi$ is a radially symmetric decreasing function with $\operatorname{supp}\psi\subset B_r$ and $0\lneqq\psi\leq 1$. Define
        \[
        W_\psi
        :=\frac{\displaystyle\int_{B_r}w(x)\,\psi(x)\,d\mu}
        {\displaystyle\int_{B_r}\psi(x)\,d\mu}.
        \]
        Then
        \[
        \int_{B_r}(w(x)-W_\psi)^2\psi(x)\,d\mu
        \leq C\,r^{2s}\iint_{B_r\times B_r}
        (w(x)-w(y))^2\min\{\psi(x),\psi(y)\}\,d\nu.
        \]
        Moreover, for every $w\in Y^{1,\gamma_1}(B_r)$,
        \[
        \int_{B_r}(w(x)-W_\psi)^2\psi(x)\,d\mu
        \leq C\,r^{2}\int_{B_r}|\nabla w|^2\psi(x)\,d\mu.
        \]
    \end{Theorem}

    Finally, we define the local weighted space
    \begin{align*}
        Y^{1,\gamma_1}_{loc}(\Omega) &:= \Bigl\{
        u\in L^2_{loc}(\Omega,d\mu) :\;
        \iint_{\Omega_1\times \Omega_1}
        (|x|^{-\gamma_1}u(x)-|y|^{-\gamma_1}u(y))^2\,d\tau\\
        &\qquad\qquad+\int_{\Omega_1}(u^2+|\nabla u|^2)\,d\mu
        <+\infty\ \text{ for every } \Omega_1\subset\subset\Omega\Bigr\}.
    \end{align*}


    \medskip
    We next establish a weighted Picone-type inequality for the operator $L_{\gamma_1}$.

    \begin{Theorem}[Picone's inequality]\label{Picone}
        Let $v\in Y^{1,\gamma_1}(\Omega)$ be such that $v>0$ in $\Omega$ and $L_{\gamma_1}(v)=\zeta$, with $\zeta\in L^1_{loc}(\mathbb{R}^N)$ and $\zeta\gneqq 0$. Then, for every $w\in C^\infty_0(\Omega)$,
        \begin{equation}\label{picone1}
            \int_{\Omega}|x|^{-2\gamma_1}|\nabla w|^2\,dx
            +\frac{a_{N,s}}{2}
            \iint_{D_\Omega}
            \frac{\bigl||x|^{-\gamma_1}w(x)-|y|^{-\gamma_1}w(y)\bigr|^2}
            {|x-y|^{N+2s}}\,dx\,dy
            \geq
            \Bigl\langle L_{\gamma_1}(v),\frac{w^2}{v}\Bigr\rangle,
        \end{equation}
        where $L_{\gamma_1}$ is defined by \eqref{LgO}.
    \end{Theorem}

    \begin{proof}
        The case $\gamma_1=0$ is treated in \cite{Amun}. We provide the necessary adaptations for general $\gamma_1$.
         Notice that
        \begin{equation*}
            \Bigl\langle L_{\gamma_1}(v),\bar{u}\Bigr\rangle
            = \underbrace{
                \int_{\Omega}\nabla v\cdot\nabla\bar{u}\,d\mu
            }_{I_1}
            +\underbrace{
                \frac{a_{N,s}}{2}\iint_{D_\Omega}
                (|x|^{-\gamma_1}v(x)-|y|^{-\gamma_1}v(y))
                (|x|^{-\gamma_1}\bar{u}(x)-|y|^{-\gamma_1}\bar{u}(y))
                \,d\tau
            }_{I_2}.
        \end{equation*}

        \noindent$\bullet$\; \emph{Estimate of $I_1$.} Adapting the proof of \cite[Theorem~2.3]{AP}, we obtain
        \begin{equation}\label{I_1}
            I_1
            =\int_{\Omega}\nabla v\cdot\nabla\!\left(\frac{w^2}{v}\right)d\mu
            \leq\int_{\Omega}|\nabla w|^2\,d\mu.
        \end{equation}

        \noindent$\bullet$\; \emph{Estimate of $I_2$.} We have
        \begin{eqnarray*}
        &(|x|^{-\gamma_1}v(x)-|y|^{-\gamma_1}v(y))
        \left(\dfrac{w^2(x)|x|^{-\gamma_1}}{v(x)}
        -\dfrac{w^2(y)|y|^{-\gamma_1}}{v(y)}\right)\\
&=|x|^{-2\gamma_1}w^2(x)+|y|^{-2\gamma_1}w^2(y)-
|x|^{-\gamma_1}|y|^{-\gamma_1}\bigg(\dfrac{v(x)}{v(y)}w^2(x)-
\dfrac{v(y)}{v(x)}w^2(y)\bigg).
\end{eqnarray*}
Since $\frac{v(x)}{v(y)}w^2(x)+
\frac{v(y)}{v(x)}w^2(y)\ge 2w(x)w(y)$, it follows that
\begin{eqnarray*}
    &(|x|^{-\gamma_1}v(x)-|y|^{-\gamma_1}v(y))
    \left(\dfrac{w^2(x)|x|^{-\gamma_1}}{v(x)}
    -\dfrac{w^2(y)|y|^{-\gamma_1}}{v(y)}\right)\\
&\le |x|^{-2\gamma_1}w^2(x)+|y|^{-2\gamma_1}w^2(y)-2|x|^{-\gamma_1}|y|^{-\gamma_1}w(x)w(y)=(|x|^{-\gamma_1}w(x)-|y|^{-\gamma_1}w(y))^2.
\end{eqnarray*}
Thus
        \begin{equation}\label{I_2}
            \frac{a_{N,s}}{2}
            \iint_{D_\Omega}
            \bigl|w(x)|x|^{-\gamma_1}-w(y)|y|^{-\gamma_1}\bigr|^2\,d\tau
            \geq
            \Bigl\langle L_{\gamma_1,s}(v),\frac{w^2}{v}\Bigr\rangle,
        \end{equation}
        where
        \[
        L_{\gamma_1,s}w(x)
        :=\frac{a_{N,s}}{|x|^{\gamma_1}}\,
        \mathrm{P.V.}
        \int_{\mathbb{R}^N}
        \frac{|x|^{-\gamma_1}w(x)-|y|^{-\gamma_1}w(y)}{|x-y|^{N+2s}}\,dy
        \]
        is the nonlocal part of $L_{\gamma_1}$. Combining \eqref{I_1} and \eqref{I_2} yields \eqref{picone1}.
    \end{proof}

    \begin{remark}\label{rem:NNN}
        A direct consequence of the Picone inequality (Theorem~\ref{Picone} with $\gamma_1=0$, see \cite{Amun}) is that, if $\lambda>\Lambda_N$, the inequality
        \begin{equation}\label{NNN}
            \mathcal{L}_{s}u-\lambda\,\frac{u}{|x|^{2}}\geq 0
            \quad\text{in }\Omega
        \end{equation}
        admits no positive weak supersolution $u$ in the sense of Definition~\ref{veryweak}, with $\dfrac{u}{|x|^2}\in L^1(\Omega)$ and $u=0$ in $\mathbb{R}^N\setminus\Omega$. Indeed, argue by contradiction. If \eqref{NNN} had a positive weak supersolution, then, by Proposition~\ref{aprggg}, there would exist $w\in W^{1,q}_0(\Omega)$, for every $q<\frac{N}{N-1}$, satisfying
        \[
        \left\{
        \begin{array}{rcll}
            \mathcal{L}_{s}w &=& \eta &\text{in }\Omega,\\
            w &>& 0 &\text{in }\Omega,\\
            w &=& 0 &\text{in }\mathbb{R}^N\setminus\Omega,
        \end{array}
        \right.
        \]
        where $\eta\in L^1(\Omega)$ satisfies $\eta\geq\lambda\,\dfrac{w}{|x|^2}\geq 0$ in $\Omega$. By an approximation argument, we may apply the Picone inequality to $w$, which yields, for every $\psi\in C^\infty_0(\Omega)$,
        \[
        \|\psi\|^2_{\mathcal{H}^1_0(\Omega)}
        \geq
        \int_\Omega
        \frac{\mathcal{L}_s(w)}{w}\,\psi^2\,dx
        \geq
        \lambda\int_\Omega\frac{\psi^2}{|x|^2}\,dx.
        \]
        Since $\lambda>\Lambda_N$, this contradicts the optimality of the Hardy constant.
    \end{remark}

    For this reason, throughout the paper we assume that $\lambda\leq\Lambda_N$.


    \medskip
    The following comparison principle extends the classical result of Brezis and Kamin \cite{BK}.

    \begin{Lemma}\label{compar}
        Let $\Omega$ be a bounded domain and let $F:[0,+\infty)\to[0,+\infty)$ be a continuous function with $F(r)>0$ for $r>0$ and such that $r\mapsto F(r)/r$ is strictly decreasing on $(0,\infty)$. Let $u,v\in Y^{1,\gamma_1}_0(\Omega)$ with $u,v>0$ in $\Omega$ and
        \[
        \left\{
        \begin{array}{rclll}
            L_{\gamma_1}(u) &\geq& F(u)&\text{\rm in } &\Omega,\\
            L_{\gamma_1}(v) &\leq& F(v)&{\rm in }&\Omega.
        \end{array}
        \right.
        \]
        Then $u\geq v$ in $\Omega$.
    \end{Lemma}

    \begin{proof}
        By an approximation argument, since $u,v>0$, we have
        \begin{equation}\label{rr1}
            \frac{L_{\gamma_1}(u)}{u}
            -\frac{L_{\gamma_1}(v)}{v}
            \geq\frac{F(u)}{u}-\frac{F(v)}{v}.
        \end{equation}
        Define $\psi:=(v^2-u^2)^+$. Then
        \begin{equation}\label{rr2}
            \int_{\Omega}\!\left(\frac{F(u)}{u}-\frac{F(v)}{v}\right)\psi\,dx
            \leq
            \int_{\Omega}\!\left(\frac{L_{\gamma_1}(u)}{u}
            -\frac{L_{\gamma_1}(v)}{v}\right)\psi\,dx=:J.
        \end{equation}
        Since $F(r)/r$ is strictly decreasing, we have $\frac{F(u)}{u}-\frac{F(v)}{v}>0$ on the set $\{v>u\}$, where $\psi>0$; hence the integrand on the left-hand side of \eqref{rr2} is nonnegative.

        We now decompose the right-hand side:
        \begin{equation*}
            \begin{split}
                J &=
                \underbrace{
                    \Bigl\langle{-\operatorname{div}(|x|^{-2\gamma_1}\nabla u)},
                    \tfrac{\psi}{u}\Bigr\rangle
                    +\frac{a_{N,s}}{2}
                    \iint_{D_\Omega}
                    (|x|^{-\gamma_1}u(x)-|y|^{-\gamma_1}u(y))
                    \!\left(\tfrac{\psi(x)|x|^{-\gamma_1}}{u(x)}
                    -\tfrac{\psi(y)|y|^{-\gamma_1}}{u(y)}\right)d\tau
                }_{J_1}\\
                &\quad-
                \underbrace{
                    \Bigl\langle{-\operatorname{div}(|x|^{-2\gamma_1}\nabla v)},
                    \tfrac{\psi}{v}\Bigr\rangle
                    +\frac{a_{N,s}}{2}
                    \iint_{D_\Omega}
                    (|x|^{-\gamma_1}v(x)-|y|^{-\gamma_1}v(y))
                    \!\left(\tfrac{\psi(x)|x|^{-\gamma_1}}{v(x)}
                    -\tfrac{\psi(y)|y|^{-\gamma_1}}{v(y)}\right)d\tau
                }_{J_2}.
            \end{split}
        \end{equation*}
For the local part, we have
        \begin{equation}\label{j11}
            \begin{array}{lll}
        & & \Bigl\langle{-\operatorname{div}(|x|^{-2\gamma_1}\nabla u)},
            \tfrac{\psi}{u}\Bigr\rangle-
                \Bigl\langle{-\operatorname{div}(|x|^{-2\gamma_1}\nabla v)},
            \tfrac{\psi}{v}\Bigr\rangle\\
        &=& \dyle\int_{\Omega\cap\{v>u\}} |x|^{-2\gamma_1}\Bigg\{\bigg(2\frac{v}{u}\nabla u\nabla v-\frac{v^2}{u^2}|\nabla u|^2-|\nabla u|^2\bigg)+\bigg(2\frac{u}{v}\nabla u\nabla v-\frac{u^2}{v^2}|\nabla v|^2-|\nabla v|^2\bigg)\Bigg\}\,dx\\
    &=& -\dyle\int_{\Omega\cap\{v>u\}} |x|^{-2\gamma_1}\Bigg\{|\nabla u+\frac{u}{v}\nabla v|^2+|\nabla v-\frac{v}{u}|\nabla u|^2\Bigg\}\,dx\le 0.\\
\end{array}
        \end{equation}
Respect to the nonlocal part, we have
$$\begin{array}{lll}
                & & (|x|^{-\gamma_1}u(x)-|y|^{-\gamma_1}u(y))
                (\tfrac{\psi(x)|x|^{-\gamma_1}}{u(x)}
                -\tfrac{\psi(y)|y|^{-\gamma_1}}{u(y)})\\
&-& (|x|^{-\gamma_1}v(x)-|y|^{-\gamma_1}v(y))
                (\tfrac{\psi(x)|x|^{-\gamma_1}}{v(x)}           -\tfrac{\psi(y)|y|^{-\gamma_1}}{v(y)})\\
            &= & |x|^{-\gamma_1}|y|^{-\gamma_1}\bigg(
            2v(x)v(y)+u(x)u(y)-\frac{u(x)}{u(y)}v^2(x)-
            \frac{u(y)}{u(x)}v^2(y)
        -\frac{v(x)}{v(y)}u^2(x)-
            \frac{v(y)}{v(x)}u^2(y)\Bigg)
        \end{array}
$$
Using the fact that
        $$
        \frac{u(x)}{u(y)}v^2(x)+
        \frac{u(y)}{u(x)}v^2(y)\ge 2 v(x)v(y)\mbox{  and }\frac{v(x)}{v(y)}u^2(x)+
        \frac{v(y)}{v(x)}u^2(y)\ge 2 u(x)u(y),$$
        it fellows that
\begin{equation}\label{j22}
(|x|^{-\gamma_1}u(x)-|y|^{-\gamma_1}u(y))
    (\tfrac{\psi(x)|x|^{-\gamma_1}}{u(x)}
    -\tfrac{\psi(y)|y|^{-\gamma_1}}{u(y)})-
    (|x|^{-\gamma_1}v(x)-|y|^{-\gamma_1}v(y))
    (\tfrac{\psi(x)|x|^{-\gamma_1}}{v(x)}
    -\tfrac{\psi(y)|y|^{-\gamma_1}}{v(y)})\le 0.
\end{equation}
    Combining \eqref{j11}--\eqref{j22}, it holds that
        \[
        \int_{\Omega}\!\left(\frac{F(u)}{u}-\frac{F(v)}{v}\right)\psi\,dx
        \leq 0.
        \]
        Since the integrand is nonnegative and strictly positive wherever $\psi>0$, this forces $\psi\equiv 0$, that is, $v\leq u$ in $\Omega$.
    \end{proof}

    \begin{remark}
        The comparison principle remains valid when $F(r)$ is replaced by a nonnegative function $H(x,r)$ that is continuous in $r$ for a.e.\ $x\in\Omega$, provided that $r\mapsto H(x,r)/r$ is strictly decreasing on $(0,\infty)$.
    \end{remark}


    \medskip
    We conclude this section with a compactness result that will be needed in the sequel.

    \begin{Lemma}\label{compa}
        Let $\{u_n\}_{n\in\mathbb{N}}\subset Y^{1,\gamma_1}_0(\Omega)$ be an increasing sequence of nonnegative functions with $L_{\gamma_1}(u_n)\geq 0$. If $u_n\rightharpoonup u$ weakly in $Y^{1,\gamma_1}_0(\Omega)$, then $u_n\to u$ strongly in $Y^{1,\gamma_1}_0(\Omega)$.
    \end{Lemma}

    \begin{proof}
        Since the sequence is increasing, we have $u_n\leq u$ a.e.\ in $\Omega$. Set $w_n:=u_n-u\leq 0$. Since $L_{\gamma_1}(u_n)\geq 0$, we obtain $\langle L_{\gamma_1}(u_n),w_n\rangle\leq 0$, that is,
        \[
        \int_{\Omega}|x|^{-2\gamma_1}\nabla u_n\cdot\nabla w_n\,dx
        +\frac{a_{N,s}}{2}\iint_{D_\Omega}
        (u_n(x)|x|^{-\gamma_1}-u_n(y)|y|^{-\gamma_1})
        (w_n(x)|x|^{-\gamma_1}-w_n(y)|y|^{-\gamma_1})\,d\tau\leq 0.
        \]
        Denoting by $(\cdot,\cdot)$ the scalar product associated with the energy norm \eqref{energy-norm}, this reads $[\![u_n]\!]^2\leq(u_n,u)$. By the Cauchy--Schwarz inequality, $[\![u_n]\!]^2\leq[\![u_n]\!]\,[\![u]\!]$, hence $[\![u_n]\!]\leq[\![u]\!]$. On the other hand, by weak lower semicontinuity of the norm, $[\![u]\!]\leq\liminf_n[\![u_n]\!]$. Therefore $[\![u_n]\!]\to[\![u]\!]$, and the strong convergence follows from the weak convergence in the Hilbert space $Y^{1,\gamma_1}_0(\Omega)$.
    \end{proof}

    \medskip
    We conclude this section by stating the following elementary algebraic inequality, which will be used in several subsequent arguments. The proof is left to the reader.

    \begin{Lemma}\label{alg}
        Let $a,b\geq 0$, $p\geq 1$ and $\alpha>0$. There exists a constant $C=C(p,\alpha)>0$ such that
        \begin{equation}\label{alg0}
            |a-b|^{p-2}(a-b)(a^\alpha-b^\alpha)
            \geq C\,\Bigl|a^{\frac{p+\alpha-1}{p}}-b^{\frac{p+\alpha-1}{p}}\Bigr|^p.
        \end{equation}
    \end{Lemma}
    %
    %
    %

    \section{Weak Harnack inequality and local behavior of nonnegative supersolutions}\label{sec-harnack}

    To handle the difficulties caused by the singular term, it is essential to obtain precise information on the behavior of solutions near the origin. When $\gamma_1=0$ in \eqref{LgO}, that is, when $\lambda=0$, such information follows directly from the strong maximum principle, which is itself a consequence of the weak Harnack inequality.

    The derivation of a weak Harnack inequality for the \emph{mixed} operator $L_{\gamma_1}$ involves difficulties that cannot be addressed by simply combining the corresponding local and nonlocal arguments. The main challenge arises from the interaction between two distinct weighted structures. The local operator is coupled with the Hardy singularity through the measure $d\mu=|x|^{-2\gamma_1}dx$, which leads naturally to weighted gradient estimates in $Y^{1,\gamma_1}_0(\Omega)$. In contrast, the fractional part involves the genuinely nonlocal quantity $v(x)|x|^{-\gamma_1}-v(y)|y|^{-\gamma_1}$, which couples the values of $v$ at two different points with different singular weights and does not admit a straightforward factorization. This difficulty becomes apparent when using the classical test function $\psi^2v^{-1}$: in the nonlocal term, it produces cross-terms involving $(|x|^{-\gamma_1}-|y|^{-\gamma_1})^2$, which yield an additional singular contribution of order $r^{N-2s-2\gamma_1}$ (see \eqref{new1}). Since $s<1$, this contribution is dominated, for small radii, by the local remainder, of order $r^{N-2-2\gamma_1}$ (see \eqref{e01}). The resulting estimate \eqref{elli30} is therefore genuinely mixed, and this interaction is captured neither by the purely local nor by the purely nonlocal theory. Although our approach is inspired by \cite{CKP} and \cite{AMPP}, the coupling of the two weighted structures requires new estimates, for which we provide complete proofs.

    \medskip
    Throughout this section, we consider the problem
    \begin{equation}\label{equi1}
        \left\{\begin{array}{rclll}
            L_{\gamma_1} v&=& h&\text{in }&\Omega,\\
            v&=&0&\text{in } &\mathbb{R}^{N}\setminus\Omega,
        \end{array}\right.
    \end{equation}
    with $h\geq 0$. We first make precise the notion of supersolution to Problem~\eqref{equi1}.

    \begin{Definition}\label{super00}
        Let $v\in Y^{1,\gamma_1}_{0}(\Omega)$. We say that $v$ is a supersolution to Problem~\eqref{equi1} if
        \begin{equation}\label{super}
            \int_{\Omega}\nabla v\cdot\nabla \varphi\, d\mu
            +\frac{a_{N,s}}{2}\iint_{\mathbb{R}^N\times\mathbb{R}^N}(v(x)|x|^{-\gamma_1} - v(y)|y|^{-\gamma_1})(\varphi(x)|x|^{-\gamma_1}- \varphi(y)|y|^{-\gamma_1})\,d\tau
            \geq \int_{\Omega}h\varphi\,dx
        \end{equation}
        for every nonnegative $\varphi \in Y_{0}^{1,\gamma_1} (\Omega)$.
    \end{Definition}

    The main result of this section is the following.

    \begin{Theorem}[Weak Harnack inequality]\label{harnack}
        Let $0<r\leq 1$ be such that $B_{2r}\subset\Omega$. Assume that $h\geq 0$ and let $v\in Y^{1,\gamma_1}_0(\Omega)$ be a supersolution to \eqref{equi1} with $v\gneqq 0$ in $\mathbb{R}^N$. Then, for every $q<\frac{N}{N-2}$, there exists a positive constant $C=C(N,s,\gamma_1,q)$ such that
        \begin{equation}\label{main}
            \Big(\frac{1}{|B_r|_{d\mu}}\int_{B_r}v^q\,d\mu\Big)^{\frac 1q}\leq C
            \inf_{B_{\frac{3}{2}r}}v.
        \end{equation}
    \end{Theorem}

\begin{remark}
    \

    \begin{enumerate}
    \item If $0\notin \overline{B_{3r/2}}$, then the desired inequality \eqref{main} follows directly from the weak Harnack inequality established in \cite{FM}.
        \item Consequently, the only case requiring further consideration is that of balls centered at the origin. Throughout the remainder of the proof, we therefore assume that
    \[
    B_r=B_r(0)\subset B_{2r}(0)\subset\Omega,
    \]
    where $r>0$ is chosen sufficiently small, with its precise value being specified as needed at each stage of the argument.
\end{enumerate}
\end{remark}

    The proof is carried out in several steps, through a series of technical lemmas.

    \begin{Lemma}\label{lemma1}
        Let $0<R\leq 1$ be such that $B_R\subset\subset \Omega$. Suppose that $v\in Y^{1,\gamma_1}_0(\Omega)$ is a nonnegative supersolution to \eqref{equi1} such that, for some $k>0$ and some $\kappa\in (0,1]$,
        \begin{equation}\label{elli1}
            |B_r\cap \{v\geq k\}|_{d\mu}\geq \kappa\,|B_r|_{d\mu},
        \end{equation}
        with $0<r<\frac{R}{16}$. Then there exists a positive constant $C=C(N,s,\gamma_1)$ such that
        \begin{equation}\label{elli2}
            |B_{6r}\cap \{v\leq 2\delta k\}|_{d\mu}\leq
            \frac{C}{\kappa\log(\frac{1}{2\delta})}|B_{6r}|_{d\mu}
        \end{equation}
        for all $\delta\in (0,\frac 14)$.
    \end{Lemma}
    \begin{proof}
        Let $\psi\in C^\infty_0(B_R)$ be such that $0\leq \psi\leq 1$, $\operatorname{supp}\psi\subset B_{7r}$, $\psi=1$ in $B_{6r}$ and $|\nabla \psi|\leq \frac{C}{r}$. Without loss of generality, we may assume that $v>0$ in $B_R$; otherwise, we work with $v+\varepsilon$ and let $\varepsilon\to 0$ at the end.
        \\
        Using $\psi^{2}v^{-1}$ as a test function in \eqref{super}, and since $h\geq 0$, we obtain
        \[
        \int_{B_{8r}} \nabla v\cdot \nabla \left(\frac{\psi^2}{v}\right)d\mu+\iint_{\mathbb{R}^{N}\times\mathbb{R}^{N}}
        (v(x)|x|^{-\gamma_1}-v(y)|y|^{-\gamma_1})\left(\frac{\psi^{2}(x)|x|^{-\gamma_1}}{v(x)}-\frac{\psi^{2}(y)|y|^{-\gamma_1}}{v(y)}\right) d\tau\geq 0.
        \]
        Thus
        \begin{equation}\label{eq01}
            \begin{split}
                0 &\leq \int_{B_{8r}} \nabla v\cdot \nabla \left(\frac{\psi^2}{v}\right)d\mu
                + 2\int_{\mathbb{R}^N \setminus B_{8r}}\int_{B_{8r}}(v(x)|x|^{-\gamma_1}-v(y)|y|^{-\gamma_1})\frac{\psi^{2}(x)|x|^{-\gamma_1}}{v(x)}\,d\tau\\
                &\quad+\int_{B_{8r}}\int_{B_{8r}}
                (v(x)|x|^{-\gamma_1}-v(y)|y|^{-\gamma_1})\left(\frac{\psi^{2}(x)|x|^{-\gamma_1}}{v(x)}-\frac{\psi^{2}(y)|y|^{-\gamma_1}}{v(y)}\right)
                d\tau \\
                & =: J_1+J_2+J_3.
            \end{split}
        \end{equation}

        To estimate $J_1$, we follow closely the proof of \cite[Lemma~3.4]{Kuusi}. Taking into account the properties of $\psi$, we obtain a positive constant $c$, independent of the data, such that
        \begin{equation}\label{e01}
            J_1\leq -c\int_{B_{6r}}|\nabla \log v|^2\, d\mu +c\,r^{N-2\gamma_1-2}.
        \end{equation}

        We now deal with $J_2$. Write $x=|x|x'$ and $y=\rho y'$, where $|x'|=|y'|=1$ and $\rho:=|y|$. Since $v\geq 0$, we have
        \begin{equation*}
            \begin{split}
                J_2&\leq 2\int_{\mathbb{R}^N \setminus B_{8r}}\int_{B_{8r}}\psi^{2}(x)|x|^{-2\gamma_1}\,d\tau\\
                &\leq 2\int_{B_{7r}}
                \frac{\psi^2(x)}{|x|^{2\gamma_1}}\int_{8r}^{\infty}
                \frac{\rho^{N-1}}{|x|^{N+2s}}
                \Big(\int_{\mathbb{S}^{N-1}}\frac{dy'}{|\frac{\rho}{|x|}y'-x'|^{N+2s}}\Big)d\rho\,dx.
            \end{split}
        \end{equation*}
        Setting $t:=\frac{\rho}{|x|}$, we obtain
        \[
        J_2\leq
        C\int_{B_{7r}}\frac{\psi^{2}(x)}{|x|^{2\gamma_1+2s}}\int_{\frac 87}^{\infty} t^{N-1}
        \Big(\int_{\mathbb{S}^{N-1}}\frac{dy'}{|t y'-x'|^{N+2s}}\Big)dt\,dx
        = C\int_{B_{7r}}\frac{\psi^{2}(x)}{|x|^{2\gamma_1+2s}}\int_{\frac 87}^{\infty} t^{N-1}D(t)\,dt\,dx,
        \]
        where
        \begin{equation}\label{DTT}
        D(t):=2\frac{\pi^{\frac{N-1}{2}}}{\Gamma(\frac{N-1}{2})}\int_0^\pi
        \frac{\sin^{N-2}(\theta)}{(1-2t \cos(\theta)+t^2)^{\frac{N+2s}{2}}}\,d\theta.
    \end{equation}
        Taking into account the behavior of $D$ near $0$, near $1$ and at infinity (see \cite{FV}), we obtain
        \[
        \int_{\frac 87}^{\infty} t^{N-1}D(t)\,dt\leq C.
        \]
        Since $2\gamma_1+2s<N$, we conclude that
        \begin{equation}\label{eq002}
            J_2\leq  C\,r^{N-2s-2\gamma_1}.
        \end{equation}

        To estimate $J_3$, we adopt the approach of the proof of \cite[Lemma~1.3]{AMPP}. More precisely,
        \begin{equation}\label{eq003}
            \begin{split}
                J_3&=\int_{B_{6r}}\int_{B_{6r}}
                (v(x)|x|^{-\gamma_1}-v(y)|y|^{-\gamma_1})\left(\frac{\psi^{2}(x)|x|^{-\gamma_1}}{v(x)}-\frac{\psi^{2}(y)|y|^{-\gamma_1}}{v(y)}\right)
                d\tau \\
                &\quad+
                \iint_{(B_{8r}\times B_{8r}) \setminus (B_{6r}\times
                    B_{6r})}(v(x)|x|^{-\gamma_1}-v(y)|y|^{-\gamma_1})\left(\frac{\psi^{2}(x)|x|^{-\gamma_1}}{v(x)}-\frac{\psi^{2}(y)|y|^{-\gamma_1}}{v(y)}\right)
                d\tau\\
                &=J_{31}+\tilde{J}_{31}.
            \end{split}
        \end{equation}
We claim that
$$
\tilde{J}_{31}\le C\,r^{N-2s-2\gamma_1}.
$$
Notice that
$$
\begin{array}{lll}
& &(v(x)|x|^{-\gamma_1}-v(y)|y|^{-\gamma_1})\left(\frac{\psi^{2}(x)|x|^{-\gamma_1}}{v(x)}-\frac{\psi^{2}(y)|y|^{-\gamma_1}}{v(y)}\right)\le (\psi(x)|x|^{-\gamma_1}-\psi(y)|y|^{-\gamma_1})^2\\ \\
&\le & 2|x|^{-2\gamma_1}(\psi(x)-\psi(y))^2+2 \psi^2(y)(|x|^{-\gamma_1}-|y|^{-\gamma_1})^2.
\end{array}
$$
Hence
$$
\begin{array}{lll}
\tilde{J}_{31} &=& \dyle\iint_{(B_{8r}\times B_{8r}) \setminus (B_{6r}\times
    B_{6r})} |x|^{-2\gamma_1}(\psi(x)-\psi(y))^2
d\tau+\dyle\iint_{(B_{8r}\times B_{8r}) \setminus (B_{6r}\times
    B_{6r})} \psi^2(y)(|x|^{-\gamma_1}-|y|^{-\gamma1})^2
d\tau\\
&=&\tilde{J}_{311}+\tilde{J}_{312}.
\end{array}
$$
Respect to $\tilde{J}_{311}$, we have
$$
\begin{array}{lll}
\tilde{J}_{311}&\le & \dyle\int_{B_{8r}}|x|^{-2\gamma_1}\int_{|y-x|<16r}\frac{(\psi(x)-\psi(y))^2}{|x-y|^{N+2s}}dydx\\
    &\le & \dyle\frac{C||\nabla \psi||_{\infty}}{r^2}\int_{B_{8r}}|x|^{-2\gamma_1}\int_{|y-x|<16r}\frac{dy dx}{|x-y|^{N+2s-2}}\le  Cr^{N-2s-\gamma_1}.
\end{array}
$$
We consider now the term $\tilde{J}_{311}$. Using the radial computation as in estimation the term $J_2$, setting $\sigma=|x|, \rho=|y|$, it follows that
$$
\tilde{J}_{312}\le \dyle\int_0^{8r}\int_0^{8r}
(\sigma^{-\gamma_1}-\rho^{-\gamma_1})^2\sigma^{N-1}\rho^{N-1}\Big(\int_{\mathbb{S}^{N-1}}\frac{dy'}{|\rho y'-\sigma x'|^{N+2s}}\Big)d\rho\,d\sigma.
$$
Setting $t=\frac{\rho}{\sigma}$, we reach that
$$
\begin{array}{lll}
\tilde{J}_{312} &\le &  \dyle\int_0^{8r}\sigma^{N-1-2s-2\gamma_1}\int_{0}^{\frac{8r}{\sigma}}(t^{\gamma_1}-1)^2t^{N-1-2\gamma_1}D(t) dt d\sigma \\
    &\le & \dyle\int_0^{8r}\sigma^{N-1-2s-2\gamma_1}\int_{0}^{\infty}(t^{\gamma_1}-1)^2t^{N-1-2\gamma_1}D(t) dt d\sigma.
\end{array}
$$
where $D$ is defined in \eqref{DTT}. According to the behavior of $D$ near $1$ and $\infty$, we reach that
$\int_{0}^{\infty}(t^{\gamma_1}-1)^2t^{N-1-2\gamma_1}D(t) dt <\infty$. Hence
$$
\tilde{J}_{312} \le C\int_0^{8r}\sigma^{N-1-2s-2\gamma_1}d\sigma \le Cr^{N-2s-2\gamma_1}.
$$
Hence the claim follows.
Thus
  $$
            J_3\leq J_{31}+C\,r^{N-2s-2\gamma_1},
            $$
        where we have used that $\psi\equiv 1$ in $B_{6r}$; the integral over $(B_{8r}\times B_{8r}) \setminus (B_{6r}\times B_{6r})$ is estimated by means of the radial symmetry of the kernel, reproducing the computation carried out for \eqref{eq002}.

 We now treat $J_{31}$. A direct computation gives
        \begin{equation}\label{new1}
            \begin{split}
                J_{31}&=\int_{B_{6r}}\int_{B_{6r}}
                (v(x)-v(y))\left(\frac{1}{v(x)}-\frac{1}{v(y)}\right) \frac{d\tau}{|x|^{\gamma_1}|y|^{\gamma_1}}+\int_{B_{6r}}\int_{B_{6r}}\left(\frac{1}{|x|^{\gamma_1}}-\frac{1}{|y|^{\gamma_1}}\right)^{2}d\tau\\
                &\leq \int_{B_{6r}}\int_{B_{6r}}
                (v(x)-v(y))\left(\frac{1}{v(x)}-\frac{1}{v(y)}\right)d\nu +C\,r^{N-2s-2\gamma_1},
            \end{split}
        \end{equation}
        where the last integral is finite, and scales as indicated, since $2\gamma_1+2s<N$. By the proof of \cite[Lemma~1.3]{CKP1}, there exist constants $C_1, C_2>0$ such that
        \begin{equation}\label{eqlog}
            (v(x)-v(y))\left(\frac{\psi^{2}(x)}{v(x)}-\frac{\psi^{2}(y)}{v(y)}\right)\leq
            -C_1\bigl(\log v(x)-\log v(y)\bigr)^{2}\psi^{2}(y)+C_2\bigl(\psi(x)-\psi(y)\bigr)^{2}.
        \end{equation}
        Applying \eqref{eqlog} with $\psi\equiv 1$ in $B_{6r}$, we get
        \begin{equation}\label{eq004}
            \int_{B_{6r}}\int_{B_{6r}}
            (v(x)-v(y))\left(\frac{1}{v(x)}-\frac{1}{v(y)}\right)d\nu
            \leq -C_1\int_{B_{6r}}\int_{B_{6r}}\bigl(\log v(x)-\log v(y)\bigr)^{2}\,d\nu.
        \end{equation}
        Combining \eqref{eq003}--\eqref{eq004}, we obtain
        \begin{equation}\label{J3}
            J_{3} \leq -c\int_{B_{6r}}\int_{B_{6r}}(\log v(x)-\log v(y))^2 \,d\nu+c\, r^{N-2s-2\gamma_1}.
        \end{equation}
        Gathering \eqref{eq01}, \eqref{e01}, \eqref{eq002} and \eqref{J3}, and using that $r^{N-2s-2\gamma_1}\leq r^{N-2-2\gamma_1}$ since $s<1$ and $r\leq 1$, we arrive at
        \begin{equation}\label{elli30}
            \int_{B_{6r}} |\nabla \log v|^{2}\, d\mu
            + \int_{B_{6r}}\int_{B_{6r}} \bigl(\log v(x)-\log v(y)\bigr)^{2}\, d\nu
            \leq C\, r^{N-2-2\gamma_1},
        \end{equation}
        where $C=C(N,s,\gamma_1)>0$.

        Let $\delta\in (0,\frac14)$, and set $w:=\min\big\{\log\big(\frac{1}{2\delta}\big),\log\big(\frac{k}{v}\big)\big\}_+$. By \eqref{elli30},
        \begin{equation}\label{elli4}
            \int_{B_{6r}}|\nabla w|^2\, d\mu\leq \int_{B_{6r}}|\nabla \log v|^2\, d\mu\leq C\,r^{N-2-2\gamma_1}.
        \end{equation}
        In view of \eqref{elli4}, the H\"older inequality and the weighted Poincar\'e inequality of Theorem~\ref{PW} yield
        \begin{equation}\label{eq6}
            \int_{B_{6r}}|w-\langle w\rangle_{B_{6r}}|\,d\mu\leq C\,r^{1+\frac{N}{2}-\gamma_1} \left(\int_{B_{6r}}|\nabla w|^2\, d\mu \right)^{\frac 12}\leq C\,|B_{6r}|_{d\mu},
        \end{equation}
        where
        $\langle w\rangle_{B_{6r}}:=\frac{1}{|B_{6r}|_{d\mu}}\dyle\int_{B_{6r}}w\,d\mu$.

        On the other hand, $\{x\in\Omega:\;w(x)=0\}=\{x\in\Omega:\;v(x)\geq k\}$. Hence, by \eqref{elli1},
        \[
        |B_{6r}\cap \{w=0\}|_{d\mu}\geq
        \frac{\kappa}{6^{N-2\gamma_1}}|B_{6r}|_{d\mu}.
        \]
        Proceeding as in the proof of \cite[Lemma~3.1]{CKP}, we obtain
        \[
        \Big|B_{6r}\cap \Big\{w=\log\Big(\frac{1}{2\delta}\Big)\Big\}\Big|_{d\mu}\leq
        \frac{6^{N-2\gamma_1}}{\kappa \log\big(\frac{1}{2\delta}\big)}\int_{B_{6r}}|w-\langle w\rangle_{B_{6r}}|\,d\mu.
        \]
        The result then follows from \eqref{eq6}, since
        \[
        B_{6r}\cap \{v\leq 2\delta k\} = B_{6r}\cap \Big\{w=\log\Big(\frac{1}{2\delta}\Big)\Big\}.
        \]
    \end{proof}

The following lemma plays a central role in our proof.

\begin{Lemma}\label{lema2}
    Suppose that the assumptions of Lemma~\ref{lemma1} hold. Then there exists $\delta \in (0,\frac{1}{2})$, depending only on $N$, $s$, $\gamma_1$ and $\kappa$, such that
    \begin{equation}\label{estim2}
        \inf_{B_{4r}} v \geq \delta k.
    \end{equation}
\end{Lemma}

\begin{proof}
    Let $w := (\ell-v)_{+}$, where $\ell\in (\delta k,2\delta k)$, and let $\psi \in C^\infty_0(B_\rho)$ with $r \leq \rho < 6r$. Using $w\psi^{2}$ as a test function in \eqref{super}, and following the computations of \cite{AMPP,CKP,Garain}, we obtain
    \begin{equation}\label{last}
        \begin{split}
            & \int_{B_{\rho}}  |\nabla w|^2 \psi^2\, d\mu +\int_{B_{\rho}}\int_{B_\rho}
            (w(x)\psi(x)-w(y)\psi(y))^2\,d\nu\\
            &\leq C\int_{B_{\rho}}   w^2 |\nabla\psi |^2\,d\mu +C\int_{B_{\rho}}\int_{B_\rho}\max\{w(x),w(y)\}^2(\psi(x)-\psi(y))^2\,d\nu \\
            &+\ell^2|B_\rho\cap \{v<\ell\}|_{d\mu}\, \sup_{x\in
                \operatorname{supp}\psi}\int_{\mathbb{R}^N\setminus
                B_{\rho}}\frac{dy}{|x-y|^{N+2s}}\\
            & =: I_1+I_2+I_3.
        \end{split}
    \end{equation}
    Consider the sequences $\{\ell_j\}_{j\in\mathbb{N}}$, $\{\rho_j\}_{j\in\mathbb{N}}$ and $\{\bar{\rho}_j\}_{j\in\mathbb{N}}$ defined by
    \[
    \ell_j:=\delta k+2^{-j-1}\delta k,\qquad \rho_j:=4r+2^{1-j}r,\qquad
    \bar{\rho}_j:=\frac{\rho_j+\rho_{j+1}}{2}.
    \]
    Then $\ell_j\in (\delta k,2\delta k)$, $\rho_j, \bar{\rho}_j \in (4r,6r]$, and
    \[
    \ell_j-\ell_{j+1}=2^{-j-2}\delta k\geq 2^{-j-3}\ell_j.
    \]
    Set $w_j:=(\ell_j-v)_+$ and $B_j:=B_{\rho_j}$, so that
    \[
    w_j\geq 2^{-j-3}\ell_j\, \chi_{\{v<\ell_{j+1}\}}.
    \]
    Let $\psi_j\in C_0^\infty(B_{\bar{\rho}_j})$ be such that $0\leq \psi_j\leq 1$, $\psi_j\equiv 1$ in $B_{j+1}$ and $|\nabla\psi_j|\leq 2^{j+3}/r$. Then
    \begin{equation}\label{eqy1}
        \int_{B_j}   w_j^2 |\nabla\psi_j |^2\,d\mu\leq C\, 2^{2j}\ell_j^2 r^{-2}|B_j\cap \{v<\ell_{j}\}|_{d\mu}.
    \end{equation}
    Following the proof of \cite[Lemma~3.6]{AMPP}, and since $r\leq 1$, we have
    \begin{equation}\label{eqy2}
        \int_{B_j}\int_{B_j}\max\{w_j(x),w_j(y)\}^2(\psi_j(x)-\psi_j(y))^2\,d\nu\leq C\, 2^{2j}\ell_j^2 r^{-2}|B_j\cap \{v<\ell_{j}\}|_{d\mu}.
    \end{equation}
    For the term $I_3$, following closely \cite{Garain} and \cite{AMPP}, we obtain
    \begin{equation}\label{eqy3}
        I_3=\ell_j^2|B_j\cap \{v<\ell_j\}|_{d\mu}\, \sup_{x\in \operatorname{supp}\psi_j}\int_{\mathbb{R}^N\setminus
            B_j}\frac{dy}{|x-y|^{N+2s}}\leq C\,  2^{j(N+2s)}\ell_j^2 r^{-2}|B_{j}\cap \{v<\ell_j\}|_{d\mu},
    \end{equation}
    with $C=C(N,s)$. Inserting \eqref{eqy1}, \eqref{eqy2} and \eqref{eqy3} into \eqref{last}, we obtain
    \begin{equation}\label{est01}
        \int_{B_{j}}  |\nabla w_j|^2 \psi_j^2\, d\mu\leq C\, \ell_j^2\, 2^{j(N+2s+2)} r^{-2}|B_{j}\cap \{v<\ell_j\}|_{d\mu}.
    \end{equation}
    Let $2^*:=\frac{2N}{N-2}$. The Caffarelli--Kohn--Nirenberg inequality of Lemma~\ref{caff} gives
    \[
    S\Big(\int_{B_j}
    \frac{|w_j\psi_j|^{2^*}}{|x|^{2^*\gamma_1}}\,dx\Big)^{\frac{2}{2^*}}\leq
    \int_{B_{j}}  |\nabla (w_j \psi_j)|^2 |x|^{-2\gamma_1}\,dx.
    \]
    Since $w_j\psi_j\geq \ell_j-\ell_{j+1}$ in $B_{j+1}\cap \{v<\ell_{j+1}\}$, and $|x|^{-2^*\gamma_1}\geq \bar{C}r^{-(2^*-2)\gamma_1}|x|^{-2\gamma_1}$ in $B_j$, with $\bar{C}$ independent of $j$, we deduce that
    \[
    \Big(\int_{B_j} \frac{|w_j\psi_j|^{2^*}}{|x|^{2^*\gamma_1}}\,dx\Big)^{\frac{2}{2^*}}\geq \frac{C}{r^{\frac{2}{2^*}(2^*-2)\gamma_1}}(\ell_j-\ell_{j+1})^2|B_{j+1} \cap \{v<\ell_{j+1}\}|^{\frac{2}{2^*}}_{d\mu}.
    \]
    Since $|B_{j+1}|_{d\mu}\simeq r^{N-2\gamma_1}$, we have
    \[
    \frac{r^{\frac{2}{2^*}(2^*-2)\gamma_1}}
    {\left(|B_{j+1}|_{d\mu}\right)^{\frac{2}{2^*}}}
    \leq C(N,\gamma_1)\,
    r^{-(N-2-2\gamma_1)}.
    \]
    Hence
    \begin{equation}\label{plus}
        (\ell_j-\ell_{j+1})^2\Big(\frac{|B_{j+1} \cap
            \{v<\ell_{j+1}\}|_{d\mu}}{|B_{j+1}|_{d\mu}}\Big)^{\frac{2}{2^*}}\leq
        C(N,\gamma_1)\,r^{-(N-2\gamma_1-2)}\int_{B_{j}}  |\nabla (w_j \psi_j)|^2\, d\mu.
    \end{equation}
    Combining \eqref{eqy1}, \eqref{est01} and \eqref{plus}, we conclude that
    \begin{equation}\label{last1}
        (\ell_j-\ell_{j+1})^2\Big(\frac{|B_{j+1} \cap
            \{v<\ell_{j+1}\}|_{d\mu}}{|B_{j+1}|_{d\mu}}\Big)^{\frac{2}{2^*}}\leq C\, 2^{j(N+2+2s)}\ell_j^2\, \frac{|B_{j} \cap
            \{v<\ell_{j}\}|_{d\mu}}{|B_{j}|_{d\mu}},
    \end{equation}
    and thus, setting $A_j:=\frac{|B_{j} \cap \{v<\ell_j\}|_{d\mu}}{|B_{j}|_{d\mu}}$,
    \begin{equation}\label{lastt1}
        A_{j+1}^{\frac{2}{2^*}}\leq C\, 2^{j(N+4+2s)}A_j.
    \end{equation}
    By Lemma~\ref{lemma1}, $A_0\leq \frac{C}{\kappa\log(\frac{1}{2\delta})}$, which can be made arbitrarily small by choosing $\delta$ small enough, depending only on $N$, $s$, $\gamma_1$ and $\kappa$. Using the same iteration argument as in the proof of \cite[Lemma~7.1]{Garain}, we deduce from \eqref{lastt1} that $A_j\to 0$ as $j\to\infty$, that is, $v\geq \delta k$ a.e.\ in $B_{4r}$.
\end{proof}

We now establish the following \emph{reverse H\"older inequality} for $v$.

\begin{Lemma}\label{dos}
    Let $0<r\leq 1$ be such that $B_{3r/2}\subset\Omega$, and assume that $v$ is a supersolution to \eqref{equi1}. Then, for every $0<\beta_1<\beta_2<\frac{N}{N-2}$,
    \begin{equation}\label{est3}
        \left(\frac{1}{|B_{r}|_{d\mu}}\int_{B_{r}}v^{\beta_2}\,d\mu\right)^{\frac{1}{\beta_2}}
        \leq C \left(\frac{1}{|B_{3r/2}|_{d\mu}}\int_{B_{3r/2}}v^{\beta_1}\,d\mu\right)^{\frac{1}{\beta_1}},
    \end{equation}
    with $C=C(N,s,\gamma_1,\beta_1,\beta_2)>0$.
\end{Lemma}

\begin{proof}
    Let $\theta\in (1,2)$ and $d>0$. Set $\bar{v}:=v+d$, and let $\psi\in C^\infty_0(\Omega)$ be such that $\operatorname{supp}\psi\subset B_{\tau r}$, $\psi=1$ in $B_{\tau' r}$ and $|\nabla \psi|\leq \frac{C}{(\tau-\tau')r}$, where $\frac 12 \leq \tau'<\tau\leq\frac 32$. Using $\bar{v}^{1-\theta}\psi^2$ as a test function in \eqref{super}, we obtain
    \begin{equation}\label{0aaa}
        \begin{split}
            0 &\leq \int_{B_{\tau r}} \nabla v\cdot \nabla \left(\frac{\psi^2}{\bar{v}^{\theta-1}}\right)d\mu +\int_{B_{\tau r}}\int_{B_{\tau r}}
            (\bar{v}(x)|x|^{-\gamma_1}-\bar{v}(y)|y|^{-\gamma_1})\left(\frac{\psi^{2}(x)|x|^{-\gamma_1}}{\bar{v}^{\theta-1}(x)}-\frac{\psi^{2}(y)|y|^{-\gamma_1}}{\bar{v}^{\theta-1}(y)}\right)
            d\tau \\
            &\quad+2\int_{\mathbb{R}^N \setminus B_{\tau r}}\int_{B_{\tau r}}(\bar{v}(x)|x|^{-\gamma_1}-\bar{v}(y)|y|^{-\gamma_1})\frac{\psi^{2}(x)|x|^{-\gamma_1}}{\bar{v}^{\theta-1}(x)}\,d\tau\\
            &=: I_1+I_2+I_3.
        \end{split}
    \end{equation}

    To estimate $I_1$, we use the same arguments as in \cite{Garain} to obtain
    \begin{equation}\label{I11N}
        \begin{split}
            I_1&\leq -c_1 \int_{B_{\tau r}} |\nabla(\bar{v}^{\frac{2-\theta}{2}})|^2 \psi^2\, d\mu + c_2\int_{B_{\tau r}} |\nabla\psi|^2 \bar{v}^{2-\theta}\, d\mu\\
            &\leq -c_1 \int_{B_{\tau r}} |\nabla(\bar{v}^{\frac{2-\theta}{2}})|^2 \psi^2\, d\mu + \frac{c_3\, r^{-2}}{(\tau-\tau')^2}\int_{B_{\tau r}} \bar{v}^{2-\theta}\, d\mu,
        \end{split}
    \end{equation}
    where $c_1$, $c_2$ and $c_3$ are positive constants depending only on $\theta$.

We now consider $I_2$. Since
$$
\begin{array}{lll}
    & & (|x|^{-\gamma_1}\bar v(x)-|y|^{-\gamma_1}\bar v(y))\Big(\frac{|x|^{-\gamma_1}\psi^2(x)}{\bar v^{\theta-1}(x)}-\frac{|y|^{-\gamma_1}\psi^2(y)}{\bar v^{\theta-1}(y)}\Big)
    \\ &= & |x|^{-\gamma_1}|y|^{-\gamma_1}\,(\bar v(x)-\bar v(y))\Big(\frac{\psi^2(x)}{\bar v^{\theta-1}(x)}-\frac{\psi^2(y)}{\bar v^{\theta-1}(y)}\Big)\\
    &+ & (|x|^{-\gamma_1}-|y|^{-\gamma_1})\big(|x|^{-\gamma_1}\psi^2(x)\bar v^{2-\theta}(x)-|y|^{-\gamma_1}\psi^2(y)\bar v^{2-\theta}(y)\big),
\end{array}
$$
then
\begin{equation}\label{rur}
    I_2\leq \int_{B_{\tau r}}\int_{B_{\tau r}}
    (\bar{v}(x)-\bar{v}(y))\left(\frac{\psi^{2}(x)}{\bar{v}^{\theta-1}(x)}-\frac{\psi^{2}(y)}{\bar{v}^{\theta-1}(y)}\right)d\nu+\mathcal{R},
\end{equation}
where \begin{equation}\label{eq:R_definition}
    \begin{aligned}
        \mathcal R=
        \iint_{B_{\tau r}\times B_{\tau r}}
        \left(
        |x|^{-\gamma_1}-|y|^{-\gamma_1}
        \right)
        \left(
        \frac{\psi^2(x)\bar v^{2-\theta}(x)}
        {|x|^{\gamma_1}}
        -
        \frac{\psi^2(y)\bar v^{2-\theta}(y)}
        {|y|^{\gamma_1}}
        \right)
        \frac{dx\,dy}{|x-y|^{N+2s}}.
    \end{aligned}
\end{equation}
By the pointwise inequality established in \cite[Lemma~3.3]{FK},
there exist positive constants $C_1$ and $C_2$, depending on
$\theta$, such that
\begin{equation}\label{off}
    \begin{split}
        &\int_{B_{\tau r}}\int_{B_{\tau r}}
        (\bar{v}(x)-\bar{v}(y))
        \left(\frac{\psi^{2}(x)}{\bar{v}^{\theta-1}(x)}
        -\frac{\psi^{2}(y)}{\bar{v}^{\theta-1}(y)}\right)d\nu\\
        &\quad\leq -C_1\int_{B_{\tau r}}\int_{B_{\tau r}}
        (\bar{v}^{\frac{2-\theta}{2}}(x)\psi(x)
        -\bar{v}^{\frac{2-\theta}{2}}(y)\psi(y))^2\,d\nu +C_2\int_{B_{\tau r}}\int_{B_{\tau r}}
        (\bar{v}^{2-\theta}(x)+\bar{v}^{2-\theta}(y))
        (\psi(x)-\psi(y))^2\,d\nu\\
        &\quad= -C_1\int_{B_{\tau r}}\int_{B_{\tau r}}
        (\bar{v}^{\frac{2-\theta}{2}}(x)\psi(x)
        -\bar{v}^{\frac{2-\theta}{2}}(y)\psi(y))^2\,d\nu
        +2C_2\int_{B_{\tau r}}\int_{B_{\tau r}}
        \bar{v}^{2-\theta}(x)(\psi(x)-\psi(y))^2\,d\nu.
    \end{split}
\end{equation}
where, in the last equality, we have used the symmetry of the integrand.

Following the same radial computation as in the proof of the
previous lemma, we obtain
\begin{equation}\label{off88}
    \begin{array}{lll}
        &\dyle\int_{B_{\tau r}}\int_{B_{\tau r}}
        (\bar{v}(x)-\bar{v}(y))
        \left(\frac{\psi^{2}(x)}{\bar{v}^{\theta-1}(x)}
        -\frac{\psi^{2}(y)}{\bar{v}^{\theta-1}(y)}\right)d\nu\\
        &\le -C_1\dyle\int_{B_{\tau r}}\int_{B_{\tau r}}
        (\bar{v}^{\frac{2-\theta}{2}}(x)\psi(x)
        -\bar{v}^{\frac{2-\theta}{2}}(y)\psi(y))^2\,d\nu
        +\frac{Cr^{-2}}{(\tau-\tau')^2}
        \int_{B_{\tau r}}\bar{v}^{2-\theta}\,d\mu.
    \end{array}
\end{equation}

\medskip

We now turn to the estimate of $\mathcal{R}$. By H\"older's
inequality, we obtain
\[
|\mathcal R|\le\dyle\bigg(
\iint_{B_{\tau r}\times B_{\tau r}}
\frac{(|x|^{-\gamma_1}-|y|^{-\gamma_1})^2}{|x-y|^{N+2s}}dxdy
\bigg)^{\frac12}
\bigg(\iint_{B_{\tau r}\times B_{\tau r}}
\frac{(\psi^2(x)\bar v^{2-\theta}(x)|x|^{-\gamma_1}
    -\psi^2(y)v^{2-\theta}(y)|y|^{-\gamma_1})^2}
{|x-y|^{N+2s}}dxdy\bigg)^{\frac12}.
\]
Since $v\in Y^{1,\gamma_1}_0(\Omega)$ and
$\psi\in\mathcal{C}^{\infty}_0(B_{\tau r})$, the same radial
computation as in the proof of Lemma~\ref{lemma1} yields
$|\mathcal R|<\infty$. We now proceed to estimate
$\mathcal R$. To this end, observe that
\begin{align}
    &\left(|x|^{-\gamma_1}-|y|^{-\gamma_1}\right)
    \left(
    \frac{\psi^2(x)\bar v^{2-\theta}(x)}{|x|^{\gamma_1}}
    -\frac{\psi^2(y)\bar v^{2-\theta}(y)}{|y|^{\gamma_1}}
    \right)\notag\\
    &=
    \frac{\psi^2(x)\bar v^{2-\theta}(x)}{|x|^{\gamma_1}}
    \left(|x|^{-\gamma_1}-|y|^{-\gamma_1}\right)\notag +
    \frac{\psi^2(y)\bar v^{2-\theta}(y)}{|y|^{\gamma_1}}
    \left(|y|^{-\gamma_1}-|x|^{-\gamma_1}\right).
    \label{eq:R_decomposition}
\end{align}
Consequently,
\[
\begin{array}{lll}
    \mathcal R &=& 2\dyle\int_{B_{\tau r}}
    \frac{\psi^2(x)\bar v^{2-\theta}(x)}{|x|^{\gamma_1}}
    \times
    \left[
    \int_{\substack{y\in B_{\tau r}}}
    \frac{|x|^{-\gamma_1}-|y|^{-\gamma_1}}{|x-y|^{N+2s}}\,dy
    \right]dx.
\end{array}
\]
For $x\in B_{\tau r}\setminus\{0\}$ and
$y\in\mathbb{R}^N\backslash B_{\tau r}$, we have
\[
\left(|x|^{-\gamma_1}-|y|^{-\gamma_1}\right)
\left(
\frac{\psi^2(x)\bar v^{2-\theta}(x)}{|x|^{\gamma_1}}
-\frac{\psi^2(y)\bar v^{2-\theta}(y)}{|y|^{\gamma_1}}
\right)\ge0.
\]
It follows that
\[
\mathcal R\le 2\dyle\int_{B_{\tau r}}
\frac{\psi^2(x)\bar v^{2-\theta}(x)}{|x|^{\gamma_1}}
\times
\left[
\int_{\mathbb{R}^N}
\frac{|x|^{-\gamma_1}-|y|^{-\gamma_1}}{|x-y|^{N+2s}}\,dy
\right]dx.
\]
Taking into account that
\[
\begin{array}{lll}
    \operatorname{PV}\dyle\int_{B_{\tau r}}
    \frac{|x|^{-\gamma_1}-|y|^{-\gamma_1}}{|x-y|^{N+2s}}\,dy
    &\le&
    \operatorname{PV}\dyle\int_{\mathbb R^N}
    \frac{|x|^{-\gamma_1}-|y|^{-\gamma_1}}{|x-y|^{N+2s}}\,dy\\
    &=&(-\Delta)^s|x|^{-\gamma_1}=C|x|^{-\gamma_1-2s},
\end{array}
\]
we obtain
\begin{align*}
    \mathcal R
    &\le C\int_{B_{\tau r}}
    \frac{\psi^2(x)\bar v^{2-\theta}(x)}{|x|^{\gamma_1}}
    |x|^{-\gamma_1-2s}\,dx\\
    &\le C\int_{B_{\tau r}}
    \frac{\psi^2(x)\bar v^{2-\theta}(x)}{|x|^{2s}}\,d\mu
    \le Cr^{2-2s}
    \int_{B_{\tau r}}
    \frac{\psi^2(x)\bar v^{2-\theta}(x)}{|x|^2}\,d\mu.
\end{align*}

Now, applying the classical weighted Hardy inequality, we obtain
\[
\begin{array}{lll}
    \mathcal R
    &\le&
    \dyle Cr^{2-2s}
    \int_{B_{\tau r}}
    |\nabla(\psi\bar v^{\frac{2-\theta}{2}})|^2\,d\mu\\
    &\le&
    \dyle Cr^{2-2s}
    \int_{B_{\tau r}}\psi^2
    \left|\nabla\left(\bar v^{\frac{2-\theta}{2}}\right)\right|^2\,d\mu
    +
    Cr^{2-2s}
    \int_{B_{\tau r}}\bar v^{2-\theta}|\nabla\psi|^2\,d\mu.
    \label{eq:R_before_cutoff}
\end{array}
\]
Recalling that
\[
|\nabla\psi|\le\frac{C}{(\tau-\tau')r},
\]
we deduce that
\begin{equation}\label{R_final}
    \begin{aligned}
        \mathcal R\le{}&
        Cr^{2-2s}\int_{B_{\tau r}}\psi^2
        \left|\nabla\left(\bar v^{\frac{2-\theta}{2}}\right)\right|^2\,d\mu
        +\frac{Cr^{-2s}}{(\tau-\tau')^2}
        \int_{B_{\tau r}}\bar v^{2-\theta}\,d\mu.
    \end{aligned}
\end{equation}
Combining \eqref{off88} and \eqref{R_final}, we conclude that
\begin{equation}\label{I22N}
    \begin{array}{lll}
        I_2&\leq&-C_1\dyle\int_{B_{\tau r}}\int_{B_{\tau r}}
        (\bar{v}^{\frac{2-\theta}{2}}(x)\psi(x)
        -\bar{v}^{\frac{2-\theta}{2}}(y)\psi(y))^2\,d\nu\\
        &&+\dyle \frac{Cr^{-2}}{(\tau-\tau')^2}
        \int_{B_{\tau r}}\bar{v}^{2-\theta}\,d\mu
        +Cr^{2-2s}\int_{B_{\tau r}}\psi^2
        \left|\nabla\left(\bar v^{\frac{2-\theta}{2}}\right)\right|^2\,d\mu\\
        &&+\dyle \frac{Cr^{-2s}}{(\tau-\tau')^2}
        \int_{B_{\tau r}}\bar v^{2-\theta}\,d\mu.
    \end{array}
\end{equation}
We now estimate $I_3$. Since $\bar v>0$, we have
\[
I_3\leq2\int_{B_{\tau r}}\bar v^{2-\theta}(x)\psi^2(x)
\left(\int_{\mathbb R^N\setminus B_{\tau r}}
\frac{dy}{|x-y|^{N+2s}}\right)d\mu.
\]
Notice that, for $x\in\operatorname{supp}(\psi)$ and
$y\in\mathbb R^N\setminus B_{\tau r}$, we have
\[
\int_{\mathbb R^N\setminus B_{\tau r}}
\frac{dy}{|x-y|^{N+2s}}
\le\frac{C}{\operatorname{dist}^{2s}(x,\partial B_{\tau r})}.
\]
Consequently,
\[
I_3\leq2\int_{B_{\tau r}}
\frac{\bar v^{2-\theta}(x)\psi^2(x)}
{\operatorname{dist}^{2s}(x,\partial B_{\tau r})}\,d\mu
\le2(\tau r)^{2-2s}\int_{B_{\tau r}}
\frac{\bar v^{2-\theta}(x)\psi^2(x)}
{\operatorname{dist}^2(x,\partial B_{\tau r})}\,d\mu.
\]
Applying the Hardy inequality, we obtain
\[
\begin{array}{lll}
    I_3&\le&
    \dyle Cr^{2-2s}\int_{B_{\tau r}}
    |\nabla(\psi \bar{v}^{\frac{2-\theta}{2}})|^2\,d\mu\\
    &\le&
    \dyle Cr^{2-2s}\int_{B_{\tau r}}\psi^2
    \left|\nabla\left(\bar v^{\frac{2-\theta}{2}}\right)\right|^2\,d\mu
    +
    Cr^{2-2s}\int_{B_{\tau r}}\bar v^{2-\theta}|\nabla\psi|^2\,d\mu.
\end{array}
\]
Hence,
\begin{equation}\label{I33N}
    I_3\le\dyle Cr^{2-2s}
    \int_{B_{\tau r}}\psi^2
    \left|\nabla\left(\bar v^{\frac{2-\theta}{2}}\right)\right|^2\,d\mu
    +\frac{Cr^{-2s}}{(\tau-\tau')^2}
    \int_{B_{\tau r}}\bar v^{2-\theta}\,d\mu.
\end{equation}
Combining estimates \eqref{I11N}, \eqref{I22N} and \eqref{I33N}, choosing $r$ small enough, we conclude that
    \[
    \int_{B_{\tau r}} |\nabla(\bar{v}^{\frac{2-\theta}{2}})|^2 \psi^2 \, d\mu \leq \frac{C\, r^{-2}}{(\tau -\tau ')^2} \int_{B_{\tau r}}\bar{v}^{2 - \theta} \, d\mu.
    \]
    Therefore, by this inequality and the Caffarelli--Kohn--Nirenberg inequality of Lemma~\ref{caff}, arguing as in the proof of Lemma~\ref{lema2},
    \[
    \left(\frac{1}{|B_{\tau' r}|_{d\mu}}\int_{B_{\tau'r}}\bar{v}^{\frac{(2-\theta)N}{N-2}}\,d\mu\right)^{\frac{N-2}{N}}
    \leq
    \frac{C}{(\tau-\tau')^2}\,\frac{1}{|B_{\tau r}|_{d\mu}}\int_{B_{\tau r}}\bar{v}^{2-\theta}\,d\mu.
    \]
    This gives \eqref{est3} for $\bar v$, with $\beta_1=2-\theta$ and $\beta_2=\frac{(2-\theta)N}{N-2}$. The general case $0<\beta_1<\beta_2<\frac{N}{N-2}$ follows by a finite iteration of this estimate on a sequence of radii between $r$ and $\frac{3}{2}r$, together with H\"older's inequality. Finally, letting $d\to 0$ and invoking the monotone convergence theorem, we conclude the proof.
\end{proof}

In order to derive the weak Harnack inequality, we establish the following estimate, along the lines of \cite{GI}.

\begin{Lemma}\label{tres}
    Let $0<r\leq 1$ be such that $B_r\subset\Omega$, and assume that $v$ is a supersolution to \eqref{equi1}. Then there exist constants $\eta\in (0,1)$ and $C>0$, depending only on $N$, $s$ and $\gamma_1$, such that
    \begin{equation}\label{est31}
        \Big(\frac{1}{|B_r|_{d\mu}}\int_{B_r}v^{\eta}\,d\mu\Big)^{\frac{1}{\eta}}\leq
        C\inf_{B_r} v.
    \end{equation}
\end{Lemma}

To prove this lemma, we use the following covering result, established in \cite{KS} (see also \cite{CKP}), which is derived from the Krylov--Safonov theory.

\begin{Lemma}\label{cover}
    Let $E\subset B_r(x_0)$ be a measurable set. For $\bar{\delta}\in (0,1)$, define
    \[
    [E]_{\bar{\delta}}:= \bigcup_{\rho>0}\bigl\{B_{3\rho}(x)\cap
    B_r(x_0) : x\in B_r(x_0),\ |E\cap B_{3\rho}(x)|_{d\mu}>\bar{\delta}\,|
    B_{\rho}(x)|_{d\mu}\bigr\}.
    \]
    Then either
    \begin{enumerate}
        \item $|[E]_{\bar{\delta}}|_{d\mu}\geq\frac{\tilde{C}}{\bar{\delta}}|E|_{d\mu}$, or
        \item $[E]_{\bar{\delta}}=B_r(x_0)$,
    \end{enumerate}
    where $\tilde{C}$ depends only on $N$, $s$ and $\gamma_1$.
\end{Lemma}

\begin{pfn}{Lemma~\ref{tres}}
    The argument follows the lines of \cite{AMPP}, with the modifications required by the local term. For completeness, we give a detailed proof. Let $\eta>0$. Then
    \begin{equation}\label{rep}
        \frac{1}{|B_r|_{d\mu}}\int_{B_r}v^{\eta}\,d\mu=\eta
        \int_0^\infty
        t^{\eta-1}\frac{|B_r\cap\{v>t\}|_{d\mu}}{|B_r|_{d\mu}}\,dt.
    \end{equation}
    Let $t > 0$ and $i \in \mathbb{N}$, and consider the set $A^i_t := \{ x \in B_r : v(x) > t \delta^i \}$, where $\delta$ is given by Lemma~\ref{lema2} with $\kappa=\bar\delta\,3^{-(N-2\gamma_1)}$. Notice that $A^{i-1}_t \subset A^i_t$.

    Fix $\rho > 0$ and $x \in B_r$ such that $B_{3\rho}(x) \cap B_r \subset [A^{i-1}_t]_{\bar{\delta}}$, so that
    \[
    |A^{i-1}_t \cap B_{3\rho}(x)|_{d\mu} > \bar{\delta}\, |B_\rho(x)|_{d\mu} = \frac{\bar{\delta}}{3^{N - 2\gamma_1}}\, |B_{3\rho}(x)|_{d\mu}.
    \]
    By Lemma~\ref{lema2}, applied with $k=t\delta^{i-1}$,
    \[
    v(z)>\delta \,(t\delta^{i-1})=t\delta^i \quad \text{for a.e. } z\in B_{3\rho}(x)\cap B_r,
    \]
    which gives the inclusion $[A^{i-1}_t]_{\bar{\delta}}\subset A^{i}_t$. By Lemma~\ref{cover}, we deduce that
    \begin{equation}\label{alt1}
        A^{i}_t=B_r\quad\hbox{or}\quad|A^{i}_t|_{d\mu}\geq
        \frac{\tilde{C}}{\bar{\delta}}|A^{i-1}_t|_{d\mu}.
    \end{equation}
    We choose $\bar\delta$ small enough so that $\tilde{C}/\bar{\delta}>1$. Assume that, for some $m\in \mathbb{N}$,
    \begin{equation}\label{nesr}
        |A^0_t|_{d\mu}>\left(\frac{\bar{\delta}}{\tilde{C}}\right)^m|B_r|_{d\mu}.
    \end{equation}
    We claim that $A^{m}_t=B_r$. Otherwise, $A^{i-1}_t\subset A^{m}_t\subsetneq B_r$ for all $i\leq m$, so that the second alternative in \eqref{alt1} holds at each step, and
    \[
    |A^{m-1}_t|_{d\mu}\geq
    \frac{\tilde{C}}{\bar{\delta}}|A^{m-2}_t|_{d\mu}\geq\cdots\geq
    \left(\frac{\tilde{C}}{\bar{\delta}}\right)^{m-1}|A^{0}_t|_{d\mu}>\frac{\bar{\delta}}{\tilde{C}}\,|B_r|_{d\mu}.
    \]
    Hence $|A^{m}_t|_{d\mu}\geq\frac{\tilde{C}}{\bar\delta}|A^{m-1}_t|_{d\mu}>|B_r|_{d\mu}$, which contradicts $A^{m}_t\subsetneq B_r$. Therefore $A^{m}_t=B_r$.

    Notice that \eqref{nesr} holds if
    \begin{equation}\label{mmm}
        m>\frac{1}{\log(\frac{\bar{\delta}}{\tilde{C}})}\log\left(\frac{|A^{0}_t|_{d\mu}}{|B_r|_{d\mu}}\right).
    \end{equation}
    Choosing $m$ as the smallest integer satisfying \eqref{mmm}, we have $m\geq 1$ and
    \[
    0\leq m-1\leq
    \frac{1}{\log(\frac{\bar{\delta}}{\tilde{C}})}\log\left(\frac{|A^{0}_t|_{d\mu}}{|B_r|_{d\mu}}\right).
    \]
    Therefore, since $\delta\in (0,\frac 12)$,
    \[
    \inf_{B_r}v\geq t\delta^m\geq t\delta\left(\frac{|A^{0}_t|_{d\mu}}{|B_r|_{d\mu}}\right)^{\frac{1}{\beta}},
    \qquad\text{with }
    \beta:=\frac{\log(\frac{\bar{\delta}}{\tilde{C}})}{\log\delta}>0.
    \]
    Set $\xi:=\inf_{B_r}v$. Then
    \[
    \frac{|B_r\cap\{v>t\}|_{d\mu}}{|B_r|_{d\mu}}=\frac{|A^{0}_t|_{d\mu}}{|B_r|_{d\mu}}\leq
    \delta^{-\beta}t^{-\beta}\xi^\beta.
    \]
    Going back to \eqref{rep}, for $0<\eta<\beta$ and $a>0$,
    \[
    \frac{1}{|B_r|_{d\mu}}\int_{B_r}v^{\eta}\,d\mu\leq
    \eta\int_0^a t^{\eta-1}\,dt +\eta \int_a^\infty
    t^{\eta-1}\delta^{-\beta}t^{-\beta}\xi^\beta\, dt
    = a^\eta+\frac{\eta}{\beta-\eta}\,\delta^{-\beta}\xi^{\beta}a^{\eta-\beta}.
    \]
    Choosing $a := \xi$ and $\eta := \min\{\frac{\beta}{2},\frac12\}$, we obtain \eqref{est31}.
\end{pfn}

With this result at hand, we can now prove the weighted weak Harnack inequality.

\begin{pfn}{Theorem~\ref{harnack}}
    By Lemma~\ref{tres}, applied in $B_{\frac32 r}$, there exists $\eta\in (0,1)$ such that
    \[
    \Big(\frac{1}{|B_{\frac32 r}|_{d\mu}}\int_{B_{\frac32 r}}v^{\eta}\,d\mu\Big)^{\frac{1}{\eta}}\leq
    C\inf_{B_{\frac32 r}} v.
    \]
    Let $q<\frac{N}{N-2}$. If $q\leq\eta$, the result follows from H\"older's inequality. Otherwise, Lemma~\ref{dos} with $\beta_1=\eta$ and $\beta_2=q$ gives
    \[
    \left(\frac{1}{|B_{r}|_{d\mu}}\int_{B_{r}}v^{q}\,d\mu\right)^{\frac{1}{q}} \leq C \left(\frac{1}{|B_{\frac32 r}|_{d\mu}}\int_{B_{\frac32 r}}v^{\eta}\,d\mu\right)^{\frac{1}{\eta}}.
    \]
    Combining the two inequalities yields \eqref{main}.
\end{pfn}

\medskip
Notice that one of the main difficulties of this work is the loss of regularity of supersolutions caused by the singular potential and by the mixed local--nonlocal nature of the operator. This makes the direct use of suitable test functions particularly delicate. To overcome this difficulty, we systematically rely on an approximation procedure, which allows us to recover the required energy framework. More precisely, given a nonnegative weak supersolution $u$ in the sense of Definition~\ref{veryweak}, there exists a nondecreasing sequence $\{w_n\}_n$ of energy supersolutions to suitable approximating problems, with $|x|^{\gamma_1}w_n\in Y^{1,\gamma_1}_0(\Omega)$, such that $w_n \uparrow w \leq u$ a.e.\ in $\Omega$, where $w$ is itself a weak supersolution of the limiting problem. Theorem~\ref{harnack} then applies to $|x|^{\gamma_1}w_n$, uniformly in $n$, and hence to $|x|^{\gamma_1}w$. This approximation procedure will be used throughout the paper; it provides a rigorous justification for the test functions required in the energy estimates.

The first application of the weak Harnack inequality is the following singular lower bound, which describes the behavior near the origin of any supersolution of \eqref{prob00}.

\begin{Proposition}\label{Lm}
    Assume that $0<\lambda<\Lambda_N$, and let $u\in L^1(\Omega)$ be a nonnegative function such that $\dfrac{u}{|x|^2}\in L^1(\Omega)$ and
    \begin{equation}\label{UUU0}
        \mathcal{L}_{s}u-\lambda\frac{u}{|x|^{2}}\geq 0
        \quad\text{in }\Omega
    \end{equation}
    in the sense of Definition~\ref{veryweak}, with $u\gneqq 0$. Then, for $r>0$ small enough such that $B_r(0)\subset \subset \Omega$, there exists a positive constant $C$ such that
    \[
    u(x)\geq C|x|^{-\gamma_1}\quad\text{for a.e. } x\in B_r(0),
    \]
    where $\gamma_1$ is given by \eqref{gamma1}.
\end{Proposition}

\begin{proof}
    By Proposition~\ref{aprggg}, there exists a supersolution $w$ to \eqref{UUU0}, obtained as the limit of solutions to approximating problems, such that $w\leq u$ and
    \[
    \left\{
    \begin{array}{rcll}
        \mathcal{L}_{s}w-\lambda \dfrac{w}{|x|^2} &=& \lambda \dfrac{u-w}{|x|^2}&\text{in }\Omega,\\[2mm]
        w &>& 0 &\text{in }\Omega,\\
        w &=& 0 &\text{in }\mathbb{R}^N\setminus\Omega.
    \end{array}
    \right.
    \]
    Set $v := |x|^{\gamma_1} w$. Since $u-w\geq 0$, we have $L_{\gamma_1} v \geq 0$, where $L_{\gamma_1}$ is defined in \eqref{LgO}. Since $v>0$ in $\Omega$, the left-hand side of \eqref{main} is positive, and Theorem~\ref{harnack}, applied through the approximation procedure described above, gives $\inf_{B_r(0)} v \geq C>0$. Hence $u(x)\geq w(x)\geq C |x|^{-\gamma_1}$ for a.e.\ $x\in B_r(0)$.
\end{proof}

\begin{remark}
{
When $\lambda=\Lambda_N$, we have
    $\gamma_1=\gamma_2=\frac{N-2}{2}$ and
    $2(\gamma_1+1)=N$. Although the operator
    $\mathcal{L}_{\gamma_1}(v)$ and the space
    $Y^{1,\gamma_1}_0(\Omega)$ remain well-defined, the standard analytical tools, such as weighted Hardy and Sobolev inequalities and, more generally, Caffarelli--Kohn--Nirenberg (CKN) inequalities, may no longer be applicable in their usual form. In particular, for a general test function $\phi$, one cannot expect
    $\frac{\phi^2}{|x|^{2(\gamma_1+1)}}=\frac{\phi^2}{|x|^N}
    \in L^1(B_r(0))$.
    This critical integrability issue prevents us from directly applying the arguments used in the subcritical case.}

    \medbreak

    To establish a lower bound for a supersolution, as in Proposition~\ref{Lm}, we proceed as follows. Let $u\in L^1(\Omega)$ be a nonnegative function such that
    $\dfrac{u}{|x|^2}\in L^1(\Omega)$ and
    \begin{equation}\label{LIM}
        \mathcal{L}_{s}u-\Lambda_N\frac{u}{|x|^{2}}\geq 0
        \quad\text{in }\Omega,
    \end{equation}
    in the sense of Definition~\ref{veryweak}, and assume that $u\not\equiv 0$.

    We claim that, for every $\epsilon>0$, there exist $r>0$ sufficiently small, with $B_r(0)\Subset\Omega$, and a constant $C(\epsilon)>0$ such that
    $$
    u(x)\geq C(\epsilon)|x|^{-\gamma_1+\epsilon}
    \quad\text{for a.e. }x\in B_r(0).
    $$
    Indeed, since $\lambda=\Lambda_N$, inequality \eqref{LIM} can be rewritten as
    $$
    \mathcal{L}_{s}u-(\Lambda_N-\epsilon^2)
    \frac{u}{|x|^{2}}\geq 0
    \quad\text{in }\Omega,
    $$
    where we have used $u\geq 0$. Setting
    $\lambda_\epsilon=\Lambda_N-\epsilon^2$, we obtain
    $$\gamma_1(\lambda_\epsilon)
    =\frac{N-2}{2}
    -\sqrt{\Lambda_N-\lambda_\epsilon}
    =\frac{N-2}{2}-\epsilon
    =\gamma_1-\epsilon.$$
    Choosing $\epsilon>0$ sufficiently small so that
    $\lambda_\epsilon\geq 0$, Proposition~\ref{Lm} applies to yield
    $$
    u(x)\geq C(\epsilon)|x|^{-\gamma_1(\lambda_\epsilon)}
    =C(\epsilon)|x|^{-\gamma_1+\epsilon}
    \quad\text{for a.e. }x\in B_r(0).$$
\end{remark}

\section{Existence and nonexistence results}\label{s4}

    As mentioned in the introduction, we present in this section several applications of the weak Harnack inequality established in Section~\ref{sec-harnack} to problems involving the Hardy potential. The section is organized according to the form of the term $f$ in Problem ~\eqref{tlem}, which we recall for the convenience of the reader:
    \begin{equation*}
        \left\{
        \begin{array}{rclll}
            -\Delta u+(-\Delta)^s u
            &=&
            \lambda \dfrac{u}{|x|^{2}}+f(x,u)
            &\hbox{ in } & \Omega,\\[4pt]
            u&>&0
            &\hbox{ in }& \Omega,
            \\[4pt]
            u&=&0
            &\hbox{ in } & \mathbb{R}^N\setminus\Omega.
        \end{array}
        \right.
    \end{equation*}
    In Subsection~\ref{subsec:linear}, we consider the linear case where $f$ does not depend on $u$. We give a necessary and sufficient condition for the existence of a weak solution, together with regularity properties of the solution. Subsection~\ref{subsec:application} is devoted to the case $f(x,u)=u^q+g(x)$, for which we identify a critical exponent $q^+(\lambda)$ governing the existence of solutions. Finally, in Subsection~\ref{sec:main}, we study the singular case $f(x,u)=h(x)u^{-\sigma}$, where $h$ is either a function or a nonnegative Radon measure.

\subsection{The linear case}\label{subsec:linear}

\

In this subsection, we assume that $f(x,u)\equiv f(x)$, so that Problem~\eqref{tlem} reduces to the linear Problem~\eqref{prob00}. We give a necessary and sufficient condition for the existence of a weak solution, which complements, in some sense, the result obtained in \cite{Biagi1}.
\smallbreak
\noindent Before stating the main result, we first establish the following regularity result.

\begin{Theorem} \label{singular}
    Let $f\in L^{m}(\Omega)$ with $m>\frac{N}{2}$, and let $u\in
    \mathcal{H}_0^1(\Omega)$ be the unique energy solution to Problem
    \eqref{prob00}, with $\lambda\le   \Lambda_{N}$. Then,
    $u(x)\leq C |x|^{-\gamma_1}\inn\mathbb{R}^{N} \setminus\{0\}$.
\end{Theorem}
\begin{proof}
    Without loss of generality, we may assume that $f\ge 0$.
    We set $v(x):= |x|^{\gamma_1} u(x)$. Then $v$ solves
    \begin{equation}\label{weightedProblem}
        \left\{
        \begin{array}{rclll}
            L_{\gamma_1} (v) &=& |x|^{-\gamma_1}f &{\rm in}&\Omega,\\
            v&=&0&{\rm in}& \mathbb{R}^{N}\setminus\Omega.
        \end{array}\right.
    \end{equation}
    Using $G_{k}(v)$ as a test function in  {\eqref{weightedProblem}}, we obtain
    \begin{equation}\label{stam1}
        \begin{array}{c}
            \dyle \int_{\Omega}|x|^{-2\gamma_1} \nabla v\nabla G_{k}(v)dx \\
            +\dyle \frac{a_{N,s}}{2} \iint_{D_\Omega} \dfrac{{\big(|x|^{-\gamma_1}v(x)-|y|^{-\gamma_1}v(y) \big)\big(|x|^{-\gamma_1}G_k (v (x))- |y|^{-\gamma_1}G_k(v (y)) \big) } }{|x-y|^{N+2s}} dxdy
            = \dyle\int_\Omega{|x|^{-\gamma_1}fG_k(v)\,dx} \,.
    \end{array}\end{equation}
    Since $\sigma=T_k (\sigma)+G_k (\sigma)$, using the properties of $T_k$ and $G_k$, we obtain
    $$
    \int_{\Omega}|x|^{-2\gamma_1} |\nabla G_{k}(v(x))|^2 \, dx + C
    \iint_{D_\Omega}\frac{|G_{k}(v(x)) - G_{k}(v(y))|^2}{|x - y|^{N + 2s}} \frac{dx}{|x|^{\gamma_1}} \, \frac{dy}{|y|^{\gamma_1}} \leq \int_{\Omega} f \, G_{k}(v(x)) \, \frac{dx}{|x|^{\gamma_1}},
    $$
    where $C$ is a positive constant depending only on the data.\\
    Setting $A_k :=\{ x \in \Omega\, ; \, v(x)  \geq k\}$ and applying the weighted Sobolev inequality of Lemma \ref{caff}, it follows that
    $$
    \mathcal{S} \| G_k (v) \|_{{L^{2^*} (\Omega,|x|^{-2^*\gamma_1}\,dx)}}^2 \leq \int_{A_{k}}\,{f}\,G_k (v  (x)) \frac{dx}{|x|^{\gamma_1}}.
    $$
    Applying H\"older's inequality on the right-hand side, we obtain
    $$
    \bigg|\int_{A_{k}}\,f\,G_k (v  (x)) \frac{dx}{|x|^{\gamma_1}}\bigg|
    \leq \|f \|_{L^{m}(\Omega)} \|G_k (v )\|_{{L^{2^*}(\Omega,{|x|^{-2^*\gamma_1}\,dx})}} |A_k|^{1-\frac1{2^*}-\frac1{m}}.
    $$
    Thus,
    $$
    \| G_k (v) \|_{L^{2^*} (\Omega,|x|^{-\gamma_1}\,dx)} \leq \mathcal{S}^{-1} \|f \|_{L^{m}(\Omega)}  |A_k|^{1-\frac1{2^*}-\frac1{m}} \,.
    $$
    On the other hand, since $\Omega$ is bounded, there exists a constant $c>0$ such that $$\| G_k (v   ) \|_{{L^{2^*} (\Omega,|x|^{-2^*\gamma_1}\,dx)}}\geq c\| G_k (v   ) \|_{L^{2^*} (\Omega)}.$$ Moreover, for any $\ell>k$ we have $A_\ell \subset A_k$ and $G_k (s) \chi_{A_\ell} \geq (\ell-k)\chi_{A_\ell}$. Thus,
    $$
    (\ell-k) |A_\ell|^{\frac1{2^*}} \leq \frac{1}{{c \mathcal{S}}} \|f \|_{L^{m}(\Omega)}  |A_k|^{1-\frac1{2^*}-\frac1{m}}.
    $$
    Hence,
    $$
    |A_\ell| \leq \frac{\|f \|^{2^*}_{L^{m}(\Omega)}  |A_k|^{2^* (1-\frac1{2^*}-\frac1{m}) } }{{(c  \mathcal{S} )}^{2^*} (\ell-k)^{2^*}}\,.
    $$
    Using a Stampacchia-type estimate, we obtain the existence of $k_0$ such that $\psi (k) \equiv 0 $ for every $k \geq k_0$, and thus
    $\displaystyle\text{ ess }\sup_{\Omega} v  \leq k_0.$ The proof is complete.
\end{proof}
In the case of a general datum $f\in L^1(\Omega)$, an additional condition is needed to ensure the existence of a solution, even in a weak sense. This is the main subject of the next result.

\begin{Theorem}\label{necesariosuficiente}
    Assume that $\l\le\Lambda_{N}$, and let $f\in L^1(\Omega)$ be a nonnegative function.
    Then, Problem \eqref{prob00} admits a nonnegative weak supersolution in the sense of Definition \ref{veryweak} if and only if
    \begin{equation}\label{oups}
        \dint_{B_{r}(0)}|x|^{-\gamma_1}f\,dx<+\iy,
    \end{equation}
    for some $B_{r}(0)\subset\subset\Omega$.
    Moreover, if \eqref{oups} holds, Problem \eqref{prob00} has a minimal solution $v$ such that  $v \in W^{1,q}_0(\Omega)$ for all $q<\frac{N}{N-1}$, and $v\in L^
    {\sigma}(\Omega)$ for all $\sigma<\frac{N}{N-2}$. In addition, $T_k (v)  \in \mathcal{H}^1_0 (\Omega)$ for all $k>0$.
\end{Theorem}

\begin{proof}

    \

    {\noindent$\bullet$ \textit{Necessary condition:} Assume $u$ is a positive weak supersolution to the Problem \eqref{prob00} in the sense of Definition \ref{veryweak}, then $\dfrac{u}{|x|^{2}}\in L^1(\Omega)$, $u=0\inn \mathbb{R}^N\setminus \Omega$. According to Proposition \ref{aprggg}, we can show the existence of $\bar{u}\le u$, a supersolution to Problem \eqref{prob00} that is obtained as a limit of approximation and that satisfies
        \begin{equation}\label{QQQ}
            \left\{\begin{array}{rclll}
                \mathcal{L}_{s}\bar{u}&=&\lambda\dfrac{\bar{u}}{|x|^{2}}+\lambda\dfrac{u-\bar{u}}{|x|^{2}}+f& \mbox{\rm  in } &\Omega,\\
                \bar{u}&=&0&\mbox{\rm in } &\ren\setminus\Omega,\\
            \end{array}\right.
        \end{equation}
        moreover using Proposition \ref{Lm}, we conclude that $\bar{u}\ge C|x|^{-\gamma_1}$ in $B_r(0)\subset\subset \Omega$.}\\

    \noindent Let $w_n \in \mathcal{H}_0^1(\Omega) \cap L^\infty(\Omega)$ be the unique positive solution to
    $$
    \left\{\begin{array}{rclll}
        \mathcal{L}_{s}(w_n)&=&\lambda\dfrac{w_{n-1}}{|x|^{2}+\frac{1}{n}}+1&\mbox{\rm in } &\Omega,\\
        w_n&=&0&\mbox{\rm in } &\ren\setminus\Omega,\\
    \end{array}\right.
    $$
    with
    $$
    \left\{\begin{array}{rclll}
        \mathcal{L}_{s}(w_0) &=&1&\mbox{in } &\Omega,\\
        w_0&=&0 &\mbox{in } &\ren\setminus\Omega.\\
    \end{array}\right.
    $$
    It is not difficult to show that the sequence $\{w_n\}_n$ is increasing in $n$, and that $w_n\to \bar{w}$, the unique positive  solution to the problem
    \begin{equation}\label{varphiProb}
        \left\{\begin{array}{rclll}
            \mathcal{L}_{s}(\bar{w}) &=&\lambda\dfrac{\bar{w}}{|x|^{2}}+1 & \mbox{\rm in } &\Omega,\\
            \bar{w}&=&0&\mbox{\rm in } & \ren\setminus\Omega.\\
        \end{array}\right.
    \end{equation}

    Using an approximation argument, we may choose $w_n$ as a test function in \eqref{QQQ} to obtain
    $$\dint_{\Omega} f w_{n}\,dx\leq\dint_{\Omega}
    \bar{u}\,dx=C<\infty.$$
    Since $\{fw_{n}\}_{n\in\mathbb{N}}$ is an increasing
    sequence bounded in $L^1(\Omega)$, applying the
    Monotone Convergence Theorem together with Proposition \ref{Lm}, we obtain
    $$
    \tilde{C} \dint_{B_{r}(0)} |x|^{-\gamma_1} f\,dx\leq \dint_{\Omega}
    f\bar{w}\,dx \leq C.
    $$
    \noindent$\bullet$ \textit{Sufficient condition:} Suppose that
    $$\int_{B_r(0)} |x|^{-\gamma_1} f \, dx < +\infty,$$
    for some $B_r(0) \subset\subset \Omega$. We now define the sequence of energy solutions $u_n \in \mathcal{H}_0^1(\Omega) \cap L^\infty(\Omega)$ to the corresponding approximated problems
    \begin{equation}\label{approx}
        \left\{
        \begin{array}{rclll}
            \mathcal{L}_{s}(u_n) &=&\lambda\dfrac{u_{n-1}}{|x|^{2}+\frac{1}{n}}+f_n&\text{\rm in } & \Omega,\\
            u_{n}(x)&=&0&\text{\rm in } & \ren\setminus\Omega,
        \end{array}
        \right.
    \end{equation}
    where
    $$
    \left\{
    \begin{array}{rclll}
        (-\Delta)^{s} u_{0}&=&f_1  &\text{\rm in } & \Omega,\\
        u_{0}(x)&=&0 &\text{\rm in } &\ren\setminus\Omega,\\
    \end{array}
    \right.
    $$
    with $f_n= T_{n}(f)$ and $u_{0}\leq u_{1}\leq u_{n-1}\leq
    u_{n}\inn\mathbb{R}^{N}$. Since $f_n\ge 0$, we have $u_n(x)\ge 0$ in $\Omega$. Using $\bar{w}\in
    \mathcal{H}^{1}_{0}(\Omega)$, the solution of Problem \eqref{varphiProb},
    as a test function in \eqref{approx}, and applying Theorem \ref{singular}, it follows that
    \begin{equation}\label{eq:millprob}
        \dint_{\Omega}
        u_n\,dx\leq \dint_{\Omega} f\bar{w}\,dx\leq \tilde{C}\int_\Omega f|x|^{-\gamma}\leq C.
    \end{equation}
    Since the sequence $ \{u_n\}_{n \in \mathbb{N}} $ is increasing, we may define $ u := \lim_{n \to \infty} u_n $, and it follows that $ u \in L^1(\Omega) $.\\
    We claim that $\dfrac{u}{|x|^{2}}\in L^1(\Omega)$. First, we observe that, by the Harnack inequality, there exists $r>0$ such that $w_0\geq C>0$ a.e.\ in $B_r(0)\subset\subset\Omega$. Using $w_0$ as a test function in \eqref{approx}, we obtain
    $$
    \int_{\Omega}\frac{u_{n-1}}{|x|^{2}+\frac{1}{n}}dx\le \frac 1C\int_{B_r(0)}\frac{w_0 u_{n-1}}{|x|^{2}+\frac{1}{n}}dx+ C(r)\int_{\Omega\backslash B_r(0)}u_ndx\le C,
    $$
    where we have used  \eqref{eq:millprob}.
    Therefore, $\{\dfrac{u_{n-1}}{|x|^{2}+\frac 1n}\}_n$ is bounded in $L^1(\Omega)$, and we conclude that  $\frac{u}{|x|^{2}}\in L^1(\Omega)$
    and
    $$\lambda\dfrac{u_{n-1}}{|x|^{2}+\frac{1}{n}}+f_n \nearrow \lambda\dfrac{u}{|x|^{2}}+f \mbox{ strongly in }L^1(\Omega).$$
    The rest of the proof follows by arguing as in \cite{Biagi1}.
\end{proof}

\subsection{Existence and nonexistence results in the case
    \texorpdfstring{$f(x,u):=u^q+g(x)$}{f(x,u)=u^q+g(x)}.}%
\label{subsec:application}
\

In this subsection, we study the following problem:
\begin{equation}\label{eq:problem1}
    \left\{\begin{array}{rclll}
        \mathcal{L}_s(u)&=&\lambda \dfrac{u}{|x|^{2}} + u^q+ g(x)
        &\text{\rm in } & \Omega,\\[4pt]
        u&>&0 &\text{\rm in }  &\Omega,\\[2pt]
        u&=&0 &\text{\rm in }&\mathbb{R}^N\setminus\Omega,
    \end{array}\right.
\end{equation}
where $\Omega$ is a bounded domain of class $\mathcal{C}^{1,1}$ in
$\mathbb{R}^N$ with $0\in \Omega$,
$0 < \lambda \leq \Lambda_{N}$, $q>0$, 
and $g$ is a nonnegative function satisfying additional hypotheses
specified below.
\smallbreak
\noindent As in the local case (see \cite{BDT}; see also \cite{AMPP} for the nonlocal
case), we establish the existence of a
critical exponent governing solvability of Problem~\eqref{eq:problem1}.
More precisely, we have the following result.

\begin{Theorem}\label{nonexist}
    Let $0<\lambda\le \Lambda_N$ be fixed and define
    $q^{+}(\lambda)=1+\dfrac{2}{\gamma_1}$.
    Assume $q\ge q^+(\lambda)$ and $g\ge 0$ in $\Omega$.
    Then, Problem~\eqref{eq:problem1} has no weak positive supersolution.
    In the case $g\equiv 0$, the unique nonnegative supersolution of
    Problem~\eqref{eq:problem1} is $u\equiv 0$.
\end{Theorem}

\begin{proof}
    \textbf{Step 1: Supercritical case $0<\lambda<\Lambda_N$, $q>q^+(\lambda)$.}
    Assume, for contradiction, that Problem~\eqref{eq:problem1} has a positive
    supersolution.  Then $u>0$ a.e.\ in $\Omega$.  By the weak Harnack
    inequality proved in Theorem~\ref{harnack}, there exists $r>0$ small
    enough such that $B_r(0)\Subset \Omega$ and
    $u(x)\ge C|x|^{-\gamma_1}$, where $C$ is a suitable positive constant.
    Let $\psi\in \mathcal{C}^\infty_0(B_r(0))$.  Applying Picone's
    inequality, we obtain
    \[
    \int_{B_r(0)}|\nabla \psi|^2\,dx
    +\frac{a_{N,s}}{2}\iint_{D_{B_r(0)}}
    \frac{|\psi(x)-\psi(y)|^2}{|x-y|^{N+2s}}\,dx\,dy
    \ge
    \int_{B_r(0)}\frac{\mathcal{L}_{s}(u)}{u}\,\psi^2(x)\,dx.
    \]
    Since $u$ is a supersolution of \eqref{eq:problem1} and
    $u\ge C|x|^{-\gamma_1}$, we deduce
    \[
    \int_{B_r(0)}|\nabla \psi|^2\,dx
    +\frac{a_{N,s}}{2}\iint_{D_{B_r(0)}}
    \frac{|\psi(x)-\psi(y)|^2}{|x-y|^{N+2s}}\,dx\,dy
    \ge
    \int_{B_r(0)} u^{q-1}\psi^2(x)\,dx
    \ge C\int_{B_r(0)}\frac{\psi^2(x)}{|x|^{(q-1)\gamma_1}}\,dx.
    \]
    Since $q>q^+(\lambda)$, we have
    $\gamma_1(q-1)>\gamma_1(q^+(\lambda)-1)= 2$,
    which contradicts the Hardy inequality.  This establishes Step~1.

    \medskip

    \medskip
    \textbf{Step 2: Critical case $0<\lambda<\Lambda_N$, $q=q^+(\lambda)$.}
    By the weak Harnack inequality, there exists $r>0$ such that
    $B_r(0)\Subset \Omega$ and $u(x)\ge C|x|^{-\gamma_1}$ for some constant
    $C>0$.  Note that $q^+(\lambda)\gamma_1=\gamma_1+2$, so
    \[
    \mathcal{L}_s(u)\ge \lambda\frac{u}{|x|^{2}}+\frac{C}{|x|^{\gamma_1+2}}
    \quad\text{in } B_r(0).
    \]
    Set $w(x)=|x|^{-\gamma_1}\ln\!\bigl(\tfrac{r}{|x|}\bigr)$ for $|x|\le r$
    and $w(x)=0$ for $|x|>r$.  By direct computation,
    \[
    -\Delta w-\lambda\frac{w}{|x|^2}
    =\frac{N-2(\gamma_1+1)}{|x|^{\gamma_1+2}}
    \quad\text{in } B_r(0)\setminus\{0\}.
    \]
    Concerning the nonlocal part, for $|x|<r $, we have (since the value of $a_{N,s}$ is irrelevant, we assume here that $a_{N,s}=1$)
    \begin{eqnarray*}
        x\neq0\qquad (-\Delta)^s w(x)=  \int_{\mathbb{R}^{N}}
        \frac{(w(x)-w(y))}{|x-y|^{N+2s}}dy&=&
        \int_{B_r(0)}
        \frac{(w(x)-w(y))}{|x-y|^{N+2s}}dy\\
        &+&w(x)\int_{\mathbb{R}^{N}\setminus B_r(0)}
        \frac{dy}{|x-y|^{N+2s}}\\
    \end{eqnarray*}
    We estimate the last integral.
    Since $|x|<r<|y|$, then
    $|y-x|\ge r-|x|$. Thus
    $$
    \begin{array}{rcl}
        \dyle \int_{\mathbb{R}^{N}\setminus B_r(0)}
        \frac{dy}{|x-y|^{N+2s}}\le
        \dint\limits_{|y-x|>r-|x|}\frac{dy}{|y-x|^{N+2s}}\le\frac{C}{(r-|x|)^{2s}}
    \end{array}
    $$
    Hence
    \begin{equation}\label{Nancy1}
        w(x)\int_{\mathbb{R}^{N}\setminus B_r(0)}
        \frac{dy}{|x-y|^{N+2s}}\le  \frac{Cw(x)}{(r-|x|)^{2s}}.
    \end{equation}
    On the other hand, to estimate the fist integral, we set
    $\rho=|y|$, then $x=|x| x', y=\rho y'$
    where $|x'|=|y'|=1$. Thus
    $$
    \begin{array}{rcl}
        \dyle \int_{B_r(0)}
        \frac{(w(x)-w(y))}{|x-y|^{N+2s}}dy&=&
        \dint\limits_0^{r}
        \frac{\rho^{N-1}}{|x|^{N+2s}}(w(|x|)-w(\rho))
        \left(
        \dint\limits_{|y'|=1}\dfrac{dH^{N-1}(y')}{|x'-\frac{\rho}{|x|}
            y'|^{N+2s}} \right) \,d\rho.
    \end{array}
    $$
    Let $\sigma=\dfrac{\rho}{|x|}$, then
    $$
    \begin{array}{rcl}
        \dyle \int_{B_r(0)}
        \frac{(w(x)-w(y))}{|x-y|^{N+2s}}dy&=&
        \dint\limits_0^{\frac{r}{|x|}}
        \frac{\sigma^{N-1}}{|x|^{2s}}(w(|x|)-w(\sigma |x|))
        K(\sigma) d\sigma,
    \end{array}
    $$
    where $
    K(\sigma)=\dint\limits_{|y'|=1}\dfrac{dH^{n-1}(y')}{|x'-\sigma
        y'|^{N+2s}}$. Notice that, from \cite{FV}, we know that
    \begin{equation}\label{kkk}
        K(\sigma)=2\frac{\pi^{\frac{N-1}{2}}}{\Gamma(\frac{N-1}{2})}\int_0^\pi
        \frac{\sin^{N-2}(\theta)}{(1-2\sigma \cos
            (\theta)+\sigma^2)^{\frac{N+2s}{2}}}d\theta.
    \end{equation}
    Thus
    $$
    \begin{array}{rcl}
        & &\dyle \int_{B_r(0)}
        \frac{(w(x)-w(y))}{|x-y|^{N+2s}}dy=
        \dfrac{1}{|x|^{2s+\gamma_1}}
        \dint\limits_0^{\frac{r}{|x|}}
        \sigma^{N-1}\bigg(\ln(\frac{r}{|x|})
        -\sigma^{-\gamma_1}\ln(\frac{r}{|x|\sigma})\bigg)
        K(\sigma) d\sigma\\
        &=& \dfrac{1}{|x|^{2s+\gamma_1}}
        \dint\limits_0^{1}
        \sigma^{N-1}\bigg(\ln(\frac{r}{|x|})
        -\sigma^{-\gamma_1}\ln(\frac{r}{|x|\sigma})\bigg)
        K(\sigma) d\sigma\\
        &+& \dfrac{1}{|x|^{2s+\gamma_1}}
        \dint\limits_1^{\frac{r}{|x|}}
        \sigma^{N-1}\bigg(\ln(\frac{r}{|x|})
        -\sigma^{-\gamma_1}\ln(\frac{r}{|x|\sigma})\bigg)
        K(\sigma) d\sigma\\
        &=& I_1+I_2.
    \end{array}
    $$
    Notice that $K(\frac{1}{\xi})=\xi^{N+2s}K(\xi)$ for any $\xi>0$, hence, setting $\xi=\frac{1}{\sigma}$, it follows that
    $$
    I_1=\dfrac{1}{|x|^{2s+\gamma_1}}
    \dint\limits_{1}^\infty
    \frac{k(\xi)}{\xi^{1-2s}}\bigg(\ln(\frac{r}{|x|})
    -\xi^{\gamma_1}\ln(\frac{r}{|x|}\xi)\bigg)
    d\xi.
    $$
    Hence
    $$
    \begin{array}{rcl}
        & &\dyle \int_{B_r(0)}
        \frac{(w(x)-w(y))}{|x-y|^{N+2s}}dy\\
        &=&
        \dyle\dfrac{1}{|x|^{2s+\gamma_1}}\int_1^{\frac{r}{|x|}} \frac{K(\sigma)}{\sigma}\bigg[\ln(\frac{r}{|x|})(\sigma^{N-\gamma_1}-\sigma^{2s})(\sigma^{\gamma_1}-1)+\ln (\sigma)  (\sigma^{N-\gamma_1}-\sigma^{\gamma_1+2s})\bigg]d\sigma\\
        &-& \dyle\dfrac{1}{|x|^{2s+\gamma_1}}
        \dint\limits_{\frac{r}{|x|}}^\infty
        \frac{K(\sigma)}{\sigma^{1-2s}}\bigg(\sigma^{\gamma_1}\ln(\frac{r}{|x|}\sigma)-\ln(\frac{r}{|x|})\bigg)
        d\sigma.
    \end{array}$$
    Then
    $$
    \begin{array}{rcl}
        & &\bigg|\dyle \int_{B_r(0)}
        \frac{(w(x)-w(y))}{|x-y|^{N+2s}}dy\bigg|\\
        &\le &
        \dyle\dfrac{1}{|x|^{2s+\gamma_1}}\int_1^{\frac{r}{|x|}} \frac{K(\sigma)}{\sigma}\bigg[\ln(\frac{r}{|x|})(\sigma^{N-\gamma_1}-\sigma^{2s})(\sigma^{\gamma_1}-1)+\ln (\sigma)  (\sigma^{N-\gamma_1}-\sigma^{\gamma_1+2s})\bigg]d\sigma\\
        &+& \dyle\dfrac{1}{|x|^{2s+\gamma_1}}
        \dint\limits_{\frac{r}{|x|}}^\infty
        \frac{K(\sigma)}{\sigma^{1-2s}}\bigg(\sigma^{\gamma_1}\ln(\frac{r}{|x|}\sigma)-\ln(\frac{r}{|x|})\bigg)
        d\sigma.
    \end{array}$$
    Now, fix \(x\) such that \(|x|\le \frac{r}{2}\). Then $\frac{r}{|x|}\ge 2$.
    Using the behavior of $K$ at infinity  $K(\sigma)\simeq \sigma^{-N-2s}$, see \cite{FV}, it is clear that the second integral in the previous estimate is finite.\\
    Concerning the first integral,  the main difficulty arises near \(\sigma=1\). In this case, see \cite{FV}, we know that
    \[
    K(\sigma)\simeq |\sigma-1|^{-1-2s}.
    \]
    Moreover, using the fact that
    \[
    \bigg[\ln\!\Big(\frac{r}{|x|}\Big)
    (\sigma^{N-\gamma_1}-\sigma^{2s})(\sigma^{\gamma_1}-1)
    +\ln (\sigma)\,
    (\sigma^{N-\gamma_1}-\sigma^{\gamma_1+2s})
    \bigg]
    \simeq (\sigma-1)^2,
    \]
    we deduce that
    \[
    k(\sigma)\bigg[\ln\!\Big(\frac{r}{|x|}\Big)
    (\sigma^{N-\gamma_1}-\sigma^{2s})(\sigma^{\gamma_1}-1)
    +\ln (\sigma)\,
    (\sigma^{N-\gamma_1}-\sigma^{\gamma_1+2s})
    \bigg]
    \simeq |\sigma-1|^{1-2s}\in L^1(1,2).
    \]
    Therefore, we conclude that
    \begin{equation}\label{Nancy2}
        \bigg|
        \int_{B_r(0)}
        \frac{w(x)-w(y)}{|x-y|^{N+2s}}\,dy
        \bigg|
        \le
        \frac{C}{|x|^{\gamma_1+2s}}.
    \end{equation}
    Combining estimates \eqref{Nancy1} and \eqref{Nancy2}, we obtain
    \begin{equation}\label{Nancy3}
        \big|(-\Delta)^s w(x)\big|
        \le
        \frac{C}{|x|^{\gamma_1+2s}}
        \qquad \text{in } B_{\frac{r}{2}}(0).
    \end{equation}
    Hence, setting $w_1(x)=c_2w$, with $c_2$ small enough, we obtain that
    $$
    \mathcal{L}_s(u)-\lambda \dfrac{u} {|x|^{2}}\ge \mathcal{L}_s(w_1)-\lambda \dfrac{w_1} {|x|^{2}}
    \hbox{  in } B_{\frac{r}{2}}(0).
    $$
    It is clear that $u\ge w_1$ in  $\mathbb{R}^N\setminus B_{\frac{r}{2}}(0)$. Hence $u\ge w_1$ in $B_r(0)$. Thus $u(x)\ge C|x|^{-\gamma_1}\ln(\frac{r}{|x|})$
    in $B_r(0)$.

    Fixed $\eta<<\frac{r}{2}$ and let $\psi\in \mathcal{C}^\infty_0(B_\eta(0))$. As in the supercritical case, using Picone's inequality, we obtain that
    \begin{eqnarray*}
        \dyle\int_{B_\eta(0)}|\nabla \psi|^2 dx +\frac{a_{N,s}}{2}\iint_{D_{B_\eta(0)}} \frac{|\psi(x)-\psi(y)|^2}{|x-y|^{N+2s}}\,dx\,dy&\ge&  C
        \int_{B_\eta(0)} u^{q_+(\lambda)-1}\psi^2(x)dx\\&\ge& \dyle C\int_{ B_\eta(0)} \dfrac{\psi^2(x)}{|x|^2}\bigg(\ln(\frac{r}{|x|})\bigg)^{\frac{2}{\gamma_1}}dx.
    \end{eqnarray*}
    It is clear that, for $|x|\to 0, \bigg(\ln(\frac{r}{|x|})\bigg)^{\frac{2}{\gamma_1}}\to \infty$. Hence for $\eta$ small enough we reach a contradiction with the optimality of the Hardy constant. Hence we conclude.

    \medskip
    \textbf{Step 3: The remaining case $\lambda=\Lambda_N$, $q=q^+(\Lambda_N)=2^*-1$.}

    In this case we have  $\gamma_1=\tfrac{N-2}{2}$ and
    $q^+(\Lambda_N)=1+\frac{2}{\gamma_1}=\frac{N+2}{N-2}=2^*-1$. Hence $q^+(\Lambda_N)\gamma_1=\gamma_1+2$.

    By the weak Harnack inequality, there exists $r>0$ with $B_r(0)\Subset\Omega$ and
    a constant $C>0$ such that $u(x)\ge C|x|^{-\gamma_1}$ in $B_r(0)$.
    Using the fact that $q^+(\Lambda_N)\gamma_1=\gamma_1+2$, we get
    $$
    \mathcal{L}_s(u)\ge \Lambda_N\frac{u}{|x|^2}+\frac{C}{|x|^{\gamma_1+2}}
    \quad\text{in }B_r(0).
    $$
    Hence as in the critical case, we can prove that
    \[
    u(x)\ge C|x|^{-\gamma_1}\ln\!\Bigl(\frac{r}{|x|}\Bigr)
    \quad\text{in }B_r(0).
    \]
    Now, as above, fix $\eta\ll r/2$ and let $\psi\in\mathcal{C}^\infty_0(B_\eta(0))$.
    using Picone's inequality, it holds that
    \begin{align*}
        &\int_{B_\eta(0)}|\nabla\psi|^2\,dx
        +\frac{a_{N,s}}{2}\iint_{D_{B_\eta(0)}}
        \frac{|\psi(x)-\psi(y)|^2}{|x-y|^{N+2s}}\,dx\,dy\ge C\int_{B_\eta(0)}\frac{\psi^2}{|x|^2}
        \Bigl(\ln\frac{r}{|x|}\Bigr)^{\!\frac{2}{\gamma_1}}\,dx,
    \end{align*}a contradiction with the Hardy inequality. This establishes Step~3 and completes the
    proof of Theorem~\ref{nonexist}.
\end{proof}

As a consequence, we obtain the following complete blow-up result.

\begin{Theorem}\label{blowup}
    Assume $1<\lambda\leq \Lambda_N$, $\alpha>0$, and
    $f\in L^{\infty}(\Omega)$ with $f\gneqq 0$ in $\Omega$.
    Assume $q\geq q^{+}(\lambda)$.
    If $u_n$ denotes the minimal solution of the truncated problem
    \begin{equation}\label{trunc}
        \left\{
        \begin{array}{rcll}
            \mathcal{L}_s(u_n) &= &\lambda b_n(x)u_n
            +\dfrac{u^q_n}{1+\frac{1}{n}u_n^q}+\alpha f
            &\text{in}\; \Omega,\\[4pt]
            u_n &> & 0 &\text{in}\; \Omega, \\[2pt]
            u_n &= &0 &\text{in }\,\, \ren\setminus\Omega,
        \end{array}
        \right.
    \end{equation}
    with $b_n(x)=\dfrac{1}{|x|^2 +\frac{1}{n}}$,
    then $u_n(x_0)\to +\infty$ for all $x_0\in \Omega$.
\end{Theorem}

\begin{proof}
    We proceed by contradiction.  Let $u_n$ be the minimal positive
    solution of \eqref{trunc}; the sequence $\{u_n\}$ is nondecreasing
    in $n$.  Suppose, for the sake of contradiction, that
    $u_n(x_0)\le C$ for all $n$, for some $x_0\in \Omega$.
    Setting
    \[
    h_n(x):=\lambda b_n(x)u_n
    +\frac{u^q_n}{1+\frac{1}{n}u_n^q}+\alpha f,
    \]
    Proposition~\ref{BC} gives, for all $n\ge 1$,
    \[
    \int_{\Omega}h_n(y)\,\delta(y)\,dy\le C\,u_n(x_0)\le C_1.
    \]
    Let $\xi$ be the unique positive solution of
    \begin{equation}\label{xi}
        \left\{
        \begin{array}{rcll}
            \mathcal{L}_s(\xi) &= &1&\text{in}\; \Omega,\\
            \xi &> & 0 &\text{in}\; \Omega, \\
            \xi &= &0 &\text{in }\,\, \ren\setminus\Omega,
        \end{array}
        \right.
    \end{equation}
    so that $\xi\in L^{\infty}(\Omega)\cap\mathcal{C}^{1,\sigma}(\overline{\Omega})$
    for some $\sigma\in(0,1)$.

    {From \cite{abatan} and according to
        Remark~\ref{BBB}}, we have $\xi\asymp\delta$.
    Using $\xi$ as a test function in \eqref{xi}, we get
    \[
    \int_{\Omega}u_n\,dx
    \le\int_{\Omega}h_n(x)\,\delta(x)\,dx\le C.
    \]
    Since $\{u_n\}$ and $\{b_n(x)u_n\}$ are nondecreasing sequences,
    the Monotone Convergence Theorem yields a positive function
    $u\in L^1(\Omega)$ such that
    $u_n\to u$ strongly in $L^1(\Omega)$,
    $b_n(x)u_n\to\dfrac{u}{|x|^2}$, and
    $\dfrac{u^q_n}{1+\frac{1}{n}u_n^q}\to u^q$
    strongly in $L^1(\Omega,\delta(x)\,dx)$.
    It follows readily that $u$ is a supersolution of
    \eqref{eq:problem1}, which contradicts
    Theorem~\ref{nonexist}.  This completes the proof.
\end{proof}

To complete the picture, we now establish an existence result showing, in some sense, the optimality of the critical exponent $q^+(\lambda)$. To this end, we assume that $g=\alpha g_0$, where $\alpha\geq 0$ is a parameter and $g_0$ is a nonnegative bounded function.

\begin{Theorem}\label{exist}
    Let $1<\lambda\leq \Lambda_N$.  Assume $1<q<q^{+}(\lambda)$.  Then:
    \begin{enumerate}
        \item If $1<q<2^*-1$, Problem~\eqref{eq:problem1} with
        $\alpha=0$ has a solution
        $u\in\mathcal{H}^1_0(\Omega)$ when $\lambda<\Lambda_N$,
        and $u\in W^{1,p}_0(\Omega)$ for all $p<2$
        when $\lambda=\Lambda_N$.
        \item If $0<\lambda<\Lambda_N$ and $2^*-1\le q<q^+(\lambda)$,
        then Problem~\eqref{eq:problem1} has a positive weak
        solution for some regular datum $g_0$, provided $\alpha$
        is sufficiently small.
    \end{enumerate}
\end{Theorem}

\begin{proof}
    \
    \

    \textbf{Part~(1).}
    If $1<q<2^*-1$ and $\lambda<\Lambda_N$, the mountain-pass theorem
    (see \cite{Biagi}) yields a nontrivial positive solution
    $u\in\mathcal{H}^1_0(\Omega)$.  The same approach applies when
    $\lambda=\Lambda_N$ and $q<q^+(\Lambda_N)=2^*-1$, working in the
    space $\mathcal{H}(\Omega)$ defined in Remark~\ref{HHH}.

    \medskip
    \textbf{Part~(2).}
    {Assume $0<\lambda<\Lambda_N$ and
        $2^*-1\le q<q^+(\lambda)$.}
    Let $R>0$ be such that $\Omega\subset B_R(0)$.
    Since $q<q^+(\lambda)$, we have $q\gamma_1<\gamma_1+2$.
    {Taking into account the definition of $\gamma_1$ and
        $\gamma_2$ in \eqref{gamma1}}, we may choose $\beta\in(\gamma_1,\gamma_2)$ such
    that $q\beta<\beta+2$ and $2\beta<N-2$. \\
    With $\beta$ as above, we show that there exist $A>0$ and $\alpha_0>0$
    for which $w(x):=A|x|^{-\beta}$
    is a supersolution of Problem~\eqref{eq:problem1}:
    \begin{equation}\label{supersolution1}
        \left\{
        \begin{array}{rcll}
            \mathcal{L}_s(w)&\geq& \lambda\dfrac{w}{|x|^2}+w^q+\alpha_0 g_0
            & \text{in } \Omega, \\[4pt]
            w&>& 0 & \text{in } \Omega.
        \end{array}
        \right.
    \end{equation}
    Indeed, a direct computation gives
    \[
    \mathcal{L}_s(w)-\lambda\frac{w}{|x|^2}
    =A\Bigl[-\beta^2+(N-2)\beta-\lambda\Bigr]|x|^{-\beta-2}
    +A\,\Lambda(\beta)|x|^{-\beta-2s},
    \]
    where
    \[
    \Lambda(\beta)
    =\int_{1}^{+\infty}K(\sigma)(\sigma^\beta-1)
    (\sigma^{N-1-\beta}-\sigma^{2s-1})\,d\sigma,
    \]
    with
    \[
    K(\sigma)=2\frac{\pi^{\frac{N-1}{2}}}{\Gamma\!\left(\frac{N-1}{2}\right)}
    \int_{0}^{\pi}
    \frac{\sin^{N-2}(\theta)}
    {(1-2\sigma\cos\theta+\sigma^2)^{\frac{N+2s}{2}}}\,d\theta.
    \]
    Since $\beta<N-2s$ and
    $K(\sigma)\simeq\sigma^{-N-2s}$ as $\sigma\to+\infty$ and
    $K(\sigma)\simeq|\sigma-1|^{-1-2s}$ as $\sigma\to 1$, the integral $\Lambda(\beta)$ is well defined.
    Since $q\beta<\beta+2$, for any bounded nonnegative function $g_0$
    we may choose $A,\alpha_0>0$ such that
    \[
    A\le\frac{1}{2}\Bigl[\bigl(-\beta^2+(N-2)\beta-\lambda\bigr)
    R^{q\beta-\beta-2}
    +\Lambda(\beta)R^{q\beta-\beta-2s}\Bigr]^{\frac{1}{q-1}}
    \]
    and $\alpha_0\|g_0\|_{L^\infty(\Omega)}\le A^q R^{-\beta q}$.
    With these choices, $w=A|x|^{-\beta}$ satisfies
    \eqref{supersolution1}.
    \\
    Now define $u_0$ as the unique positive solution of
    \begin{equation}\label{eq:subsol}
        \left\{\begin{array}{rclll}
            \mathcal{L}_s(u_0)&=&\alpha_0 g_0
            &\text{in } & \Omega,\\[4pt]
            u_0&>&0 &\text{in } &\Omega,\\[2pt]
            u_0&=&0 &\text{in } &\mathbb{R}^N\setminus\Omega.
        \end{array}\right.
    \end{equation}
    Then $u_0$ is a subsolution of \eqref{eq:problem1}.
    By the sub-supersolution method, Problem~\eqref{eq:problem1} has a
    solution $u$ with $u_0\le u\le w$.

    \medskip
    It remains to show that the set
    \[
    I:=\bigl\{\alpha>0:\text{ Problem \eqref{eq:problem1} admits a
        solution}\bigr\}
    \]
    is a bounded interval.
    \smallbreak
    \noindent {Let $\alpha_0\in I$ and let $u_{\alpha_0}$ be a
        corresponding solution of \eqref{eq:problem1}.  For $0<\alpha_1<\alpha_0$,
        the monotone dependence on $\alpha$ implies that $u_{\alpha_0}$ is a
        supersolution of \eqref{eq:problem1} with parameter $\alpha_1$, while $0$
        is a strict subsolution.  Hence the sub-supersolution method yields a
        solution for $\alpha=\alpha_1$, so $\alpha_1\in I$ and $I$ is an interval.}
    \smallbreak
    \noindent
    To show $I$ is bounded, let $u$ be a solution of \eqref{eq:problem1} and
    use the first eigenfunction $\phi_1$, normalized in $L^1(\Omega)$, as a
    test function:
    \[
    \lambda_1\int_\Omega u\,\phi_1\,dx
    =\int_\Omega u^q\phi_1\,dx+\alpha\int_\Omega g_0\phi_1\,dx.
    \]
    Since $q>1$, Young's inequality gives
    \[
    \int_\Omega u^q\phi_1\,dx+\alpha\int_\Omega g_0\phi_1\,dx
    \le\frac{1}{q}\int_\Omega u^q\phi_1\,dx
    +\frac{q-1}{q}\lambda_1^{\frac{q}{q-1}}.
    \]
    Hence $\alpha\le\frac{\lambda_1^{\frac{q}{q-1}}}%
    {\displaystyle\int_\Omega g_0\phi_1\,dx}$, then $I$ is bounded.
\end{proof}


\begin{remark}
    In the spirit of the works \cite{BP}, \cite{BC}, and \cite{KV}, it
    would be of interest to find a necessary and sufficient condition,
    relating $q$ to the regularity of $g$, for the existence of a
    solution to Problem~\eqref{eq:problem1}, at least for small
    $\lambda$.
    \\
    Let $\phi\in\mathcal{T}(\Omega)$, where
    \[
    \mathcal{T}(\Omega):=\Bigl\{\phi\in \mathcal{H}^1_0(\Omega):
    \phi>0\text{ in }\Omega,\;
    \mathcal{L}_s(\phi)=\psi\in L^\infty(\Omega),\;\psi\ge 0\Bigr\},
    \]
    and let $u$ be a weak solution of \eqref{eq:problem1}.  Using $\phi$
    as a test function and the self-adjointness of $\mathcal{L}_s$, we
    obtain
    \[
    \int_{\Omega}\mathcal{L}_{s,\lambda}(\phi)\,u\,dx
    =\int_{\Omega}u^q\phi\,dx+\int_{\Omega}g\phi\,dx,
    \qquad
    \mathcal{L}_{s,\lambda}(\phi)
    :=\mathcal{L}_{s}(\phi)-\lambda\frac{\phi}{|x|^2}.
    \]
    Applying Young's inequality with exponents $q$ and
    $q'=q/(q-1)$, for every $\varepsilon>0$,
    \[
    \mathcal{L}_{s,\lambda}(\phi)\,u
    \le\varepsilon\,u^q\phi
    +C(\varepsilon)\,\frac{|\mathcal{L}_{s,\lambda}(\phi)|^{q'}}{\phi^{q'-1}}
    \quad\text{on }\{\phi>0\}.
    \]
    Choosing $\varepsilon$ small enough to absorb $\int_\Omega u^q\phi\,dx$
    and using $g\ge 0$, we conclude that a \emph{necessary} condition for
    existence is
    \[
    \int_\Omega g\phi\,dx
    \le C\int_{\{\phi>0\}}
    \frac{|\mathcal{L}_{s,\lambda}(\phi)|^{q'}}{\phi^{q'-1}}\,dx
    \qquad\text{for every }\phi\in\mathcal{T}(\Omega).
    \]
    In a forthcoming work we shall show that this condition is, in a
    suitable sense, also sufficient for existence, and we will study its
    connection with a Brezis--Cabr\'e type condition, see \cite[Theorem~5.1]{BC}), and with the weighted norm inequalities
    of Kalton--Verbitsky \cite{KV}. $\square$
\end{remark}

    We conclude this subsection with a closely related problem without source term, in which the singular weight $|x|^{-2}$ acts on the nonlinearity itself rather than on a linear term. This is the analogue, for the mixed operator $\mathcal{L}_s$, of the problem studied by Brezis and Cabr\'e
    in \cite{BC}. More precisely, we consider the following  problem:
\begin{equation}\label{PBC}
    \left\{\begin{array}{rclll}
        \mathcal{L}_s(u)&=&\dfrac{u^q}{|x|^{2}}
        &\text{in } &\Omega,\\[4pt]
        u&>&0 &\text{in } &\Omega,\\[2pt]
        u&=&0 &\text{in } &\mathbb{R}^N\setminus\Omega,
    \end{array}\right.
\end{equation}
where $q>1$.  The following nonexistence result holds.

\begin{Theorem}\label{Cabre}
    Assume $q>1$. Then, Problem~\eqref{PBC} has no positive weak solution in the sense
    that $u^q/|x|^{2}\in L^1(\Omega)$ and the equation holds in the
    weak sense.
\end{Theorem}
\begin{proof}
    Assume, by contradiction, that Problem~\eqref{PBC} has a positive weak
    solution $u$ such that $u^q/|x|^{2}\in L^1(\Omega)$.
    By the weak Harnack inequality, there exist $r>0$ and $c_0>0$ such that
    $B_r(0)\Subset\Omega$ and $u\ge c_0$ in $B_r(0)$. Thus
    \[
    \mathcal{L}_s(u)\ge\frac{c_0^q}{|x|^{2}}
    \quad\text{in } B_r(0).
    \]
    Define $w(x):=\ln(r/|x|)$ for $|x|\le r$ and $w(x):=0$ for $|x|>r$.
    By a direct computation analogous to Step~2 of the proof of
    Theorem~\ref{nonexist}, one obtains
    \[
    \mathcal{L}_s(w)\le\frac{c_1}{|x|^{2}}
    \quad\text{in }B_r(0).
    \]
    Set $w_1:=c_2 w$, with $c_2>0$ small enough so that
    $c_2 c_1\le c_0^q$. Since $w_1=0\le u$ in $\mathbb{R}^N\setminus B_r(0)$,
    the comparison principle gives $u\ge w_1$ in $B_r(0)$, whence
    $\lim_{|x|\to 0}u(x)=+\infty$.

    Consequently, since $\mathcal{L}_s(u)=u^{q-1}\dfrac{u}{|x|^{2}}$ and $q>1$,
    for any constant $c_3>\Lambda_N$ we may fix $r_1\in(0,r)$ small enough so that
    \[
    \mathcal{L}_s(u)\ge c_3\,\frac{u}{|x|^{2}}
    \quad\text{in }B_{r_1}(0).
    \]
    Let $\psi\in\mathcal{C}^\infty_0(B_{r_1}(0))$.
    Picone's inequality (Theorem~\ref{Picone}) then gives
    \[
    \int_{B_{r_1}(0)}|\nabla\psi|^2\,dx
    +\frac{a_{N,s}}{2}\iint_{D_{B_{r_1}(0)}}
    \frac{|\psi(x)-\psi(y)|^2}{|x-y|^{N+2s}}\,dx\,dy
    \ge\int_{B_{r_1}(0)}\frac{\mathcal{L}_s(u)}{u}\,\psi^2\,dx
    \ge c_3\int_{B_{r_1}(0)}\frac{\psi^2}{|x|^{2}}\,dx.
    \]
    Since $c_3>\Lambda_N$, this contradicts the optimality of the Hardy
    constant $\Lambda_N$, which does not depend on the domain. This completes
    the proof.
\end{proof}

\subsection{Existence and nonexistence  in the \texorpdfstring{$f(x,u)=h(x)u^{-\sigma}$}{f(x,u)=h(x)/u\^{}sigma}}\label{sec:main}
\

In this subsection,  we consider the following problem:
\begin{equation}\label{mainp}
    \left\{
    \begin{array}{rclll}
        \mathcal{L}_s(u) &= &\lambda \dfrac{u}{|x|^{2}}+\dfrac{h(x)}{u^{\sigma}}&{\rm in} & \Omega,\\
        u &> & 0 &{\rm in} & \Omega, \\
        u &= &0 &{\rm in } & \ren\setminus\Omega,
    \end{array}
    \right.
\end{equation}
  where $\sigma>0$, $0<\lambda\leq\Lambda_N$, and $h$ is a nonnegative datum.\\
    We first consider the case where $h$ is a function, and then the case where $h$ is a nonnegative Radon measure. In both cases, solutions are obtained as limits of solutions $u_n$ to approximating problems with bounded data $h_n$.

\subsubsection{\bf The case where $h$ is a function}

    \begin{Definition}
        A positive function $u \in \mathcal{H}_{loc}^1(\Omega)\cap L^{1}(\Omega)$ is a weak solution to problem \eqref{mainp} if
        $$\dfrac{u}{|x|^{2}}\in L^{1}_{{\rm loc}}(\Omega),\qquad\frac{h(x)}{u^\sigma} \in L^{1}_{{\rm loc}}(\Omega),
        $$
        and
        \begin{equation}
            \label{v-f}
            \dyle \int_{\Omega} \nabla u\nabla \varphi\,dx
            +\frac{a_{N,s}}{2}\iint_{\mathbb{R}^N\times \mathbb{R}^N}
            \frac{(u(x)-u(y))(\varphi(x)-\varphi(y))}{|x-y|^{N+2s}}\,dx\,dy
            =\lambda\,\dint_{\Omega}\frac{u\varphi}{|x|^2}\,dx+
            \int_{\Omega} \frac{h(x)}{u^\sigma}\,\varphi\, dx,
        \end{equation}
        for every $\varphi \in C^\infty_c(\Omega)$.
    \end{Definition}

    Since the interpretation of the boundary condition is not unambiguous, we introduce the following:

    \begin{Definition}\label{rediri}
        Let $u$ be such that $u  = 0$ in    $\mathbb{R}^N \setminus \Omega$.
        We say
        that $u\leq 0$ on $\partial\Omega$ if, for every $\e>0$, it follows that
        \[
        (u-\e)^+\in \mathcal{H}^1_0(\Omega)\,.
        \]
        We say that $u= 0$ on $\partial\Omega$ if $u$ is nonnegative and $u\leq 0$ on $\partial\Omega$.
    \end{Definition}

    Our first existence result is the following.

    \begin{Theorem}\label{m1}
        Assume that $\sigma >0 $ and $ \lambda< \Lambda_N $. Then, for every $ 0\lneq h \in L^1(\Omega)$, Problem \eqref{mainp} has a positive weak solution. Moreover:
        \begin{enumerate}
            \item let $0<\sigma\leq 1$ and assume that
            \begin{equation*}
                h\in L^{m}(\Omega) \text{ where }
                m:=\frac{2N}{N+2+\sigma(N-2)}.
            \end{equation*}
            Then,  Problem \eqref{mainp} has a weak solution $u\in \mathcal{H}^1_0(\Omega)$ with  ${\rm{ess\,inf}}_K\, u>0$
            for any compact $K\Subset\Omega$.
            \item if $ \sigma > 1 $ and $\lambda<\Lambda_N$, then $ u \in \mathcal{H}^1_{\text{loc}}(\Omega) $, with $ G_k(u) \in \mathcal{H}^1_0(\Omega) $ and $ T_k^{\frac{\sigma+1}{2}}(u) \in \mathcal{H}^1_0(\Omega) $. Additionally, if
            $$
            \left[ \frac{4\sigma}{(\sigma+1)^2} - \frac{\lambda}{\Lambda_N} \right] > 0,
            $$
            $ u^{\frac{\sigma+1}{2}} \in \mathcal{H}^1_0(\Omega) $.
        \end{enumerate}
    \end{Theorem}

    \begin{proof}
        The result follows by an approximation argument. 
        Using a minimizing argument, we obtain the existence of $u_n$, the unique positive solution to the problem
        \begin{equation}\label{apro}
            \left\{
            \begin{array}{rclll}
                \mathcal{L}_s(u_n) &= &\lambda \dfrac{u_{n-1}}{|x|^{2}}+\dfrac{h_n(x)}{(u_n+\frac{1}{n})^{\sigma}} &{\rm in} &\Omega,\\
                u_n &> & 0 &{\rm in} &\Omega, \\
                u_n &= & 0&{\rm in} &\ren\setminus\Omega,
            \end{array}
            \right.
        \end{equation}
        with $u_0\equiv 0$. Uniqueness follows from Lemma \ref{compar}. Since $h_n$ is increasing in $n$, the sequence $\{u_n\}_{n \in \mathbb{N}}$ is likewise increasing. We split the proof into two cases according to the value of $\sigma$.
        \medskip

        {\bf Case 1: $0<\sigma\leq 1$.}
        We begin with the case $\sigma=1$. Using $u_n$ as a test function in \eqref{apro} and applying the Hardy inequality, we obtain
        \begin{equation*}\label{ess}
            \left(1-\frac{\l}{\Lambda_{N}}\right)||u_n||^2_{\mathcal{H}^1_0(\Omega)}\leq\int_{\Omega}\dfrac{h_n
                u_n}{u_n+\frac{1}{n}}\,dx\leq \int_{\Omega}h \,dx=C.
        \end{equation*}
        Hence the sequence $\{u_n\}_{n \in \mathbb{N}}$ is bounded in $\mathcal{H}^1_0(\Omega)$. Thus there exists $u \in \mathcal{H}^1_0(\Omega)$ such that, up to a subsequence, $u_n$ converges weakly to $u$ in $\mathcal{H}^1_0(\Omega)$ and strongly to $u$ in $L^r(\Omega)$ for all $r < 2^*$.

        Moreover, since $\{u_n\}$ is increasing, by Proposition \ref{BC}, for all $K\Subset\Omega$ there exists $\sigma_K>0$ such that $u(x)\geq u_1(x)\geq \sigma_K>0$ for a.e.\ $x\in K$.

        Finally, since $\mathcal{L}_s(u_n)\geq 0$, applying the monotonicity of the sequence $\{u_n\}_{n \in \mathbb{N}}$ together with the compactness result from Lemma \ref{compa}, it follows that $u_n \to u$ strongly in $\mathcal{H}^1_0(\Omega)$. Consequently, $u$ is a solution to Problem \eqref{mainp}.

        \medskip

        Notice that the case $0<\sigma<1$ follows by the same argument, taking into account that $h\in L^m(\Omega)$ with $(\sigma-1)m'=2^*$.

        \medskip

        {\bf Case 2: $\sigma>1$.}
        We use $G_k(u_n)$ as a test function in \eqref{apro}, which gives
        \begin{equation*}\label{estim1}
            \begin{array}{lll}
                &\dyle \int_{\Omega}|\nabla G_k(u_n)|^2 dx+\frac{a_{N,s}}{2}\iint_{D_\Omega}{\frac{|G_k(u_n(x))-G_k(u_n(y))|^2}{|x-y|^{N+2s}}}\,
                dx\, dy
                -\lambda\int_{\Omega}\dfrac{u_nG_k(u_n)}{|x|^{2}}\,dx\\
                &\le \dyle \int_{\Omega}\dfrac{h_n G_k(u_n)}{(u_n+\frac{1}{n})^\sigma}\,dx\le
                \frac{1}{k^{\sigma-1}}\int_{\Omega} h \,dx.
        \end{array}\end{equation*}
        Since $u_n G_k(u_n) = G_k^2(u_n) + k G_k(u_n)$, we obtain
        \begin{equation*}\label{estim20}
            \begin{split}
                \int_{\Omega}|\nabla G_k(u_n)|^2 dx+\frac{a_{N,s}}{2}&\iint_{D_\Omega}\frac{|G_k(u_n(x))-G_k(u_n(y))|^2}{|x-y|^{N+2s}}\,
                dx\, dy-\lambda\int_{\Omega}\dfrac{G_k^2(u_n)}{|x|^{2}}\,dx\\
                & \qquad \;\;\;\leq\lambda k\int_{\Omega} \dfrac{G_k(u_n)}{|x|^{2}}\,dx
                +\frac{1}{k^{\sigma-1}}\int_{\Omega} h \,dx.
            \end{split}
        \end{equation*}
        Since $\lambda < \Lambda_N$, applying the Hardy--Sobolev inequality to the left-hand side yields
        \begin{equation*}\label{estim4}
            C\|G_k(u_n)\|^2_{\mathcal{H}^1_0(\Omega)}\leq \lambda k\int_{\Omega} \dfrac{G_k(u_n)}{|x|^{2}}\,dx
            +C(k,h,\sigma).
        \end{equation*}
        Applying Young's inequality and the Hardy--Sobolev inequality to the right-hand side, we get
        \begin{equation*}\label{estim5}
            \|G_k(u_n)\|^2_{\mathcal{H}^1_0(\Omega)}\leq C(h, k,\lambda,\Lambda_{N},\sigma,\Omega).
        \end{equation*}
        Thus $\{G_k(u_n)\}_{n\in\mathbb{N}}$ is uniformly bounded in $\mathcal{H}^1_0(\Omega)$, and the Hardy--Sobolev inequality further gives
        $$
        \int_{\Omega}\dfrac{G^2_k(u_n(x))}{|x|^{2}}\,dx\le
        C(h, k,\lambda,\Lambda_{N},\sigma,\Omega).
        $$
        Since $t=T_k(t)+G_k(t)$, we obtain
        \begin{equation}\label{eq:Hardyest}
            \int_{\Omega}\dfrac{u^2_n(x)}{|x|^{2}}\,dx\le
            C(h, k,\lambda,\Lambda_{N},\sigma,\Omega).
        \end{equation}
        \noindent Next, testing with $ T_k^\sigma(u_n) $ in \eqref{apro} gives
        \begin{equation*}\label{estim6}
            \begin{split}
                \dyle \dfrac{4\sigma}{(\sigma+1)^2}\int_{ \Omega}\left|\nabla T_k^{\frac{\sigma+1}{2}}(u_n)\right|^2dx +
                \frac{a_{N,s}}{2}&\iint_{D_\Omega}{\frac{(T_k^\sigma(u_n(x))-T_k^\sigma(u_n(y)))(u_n(x)-u_n(y))}{|x-y|^{N+2s}}}
                dx\, dy\\
                &\leq \displaystyle
                \lambda\int_{\Omega}\dfrac{u_nT^{\sigma}_k(u_n)}{|x|^{2}}\,dx+\int_{\Omega}\dfrac{h_n
                    T^\sigma_k(u_n)}{(u_n+\frac{1}{n})^\sigma}\,dx\\
                &\le \displaystyle k^{\sigma-1}\lambda\int_{\Omega}\dfrac{u^2_n}{|x|^{2}}\,dx+\int_{\Omega}h_n
                \,dx\le C(h, k,\lambda,\Lambda_{N},\sigma,\Omega),
            \end{split}
        \end{equation*}
        and hence, using \eqref{eq:Hardyest},
        $$
        \left\|T^{\frac{\sigma+1}{2}}_k(u_n)\right\|^2_{\mathcal{H}^1_0(\Omega)}\leq C(h, k,\lambda,\Lambda_{N},\sigma,\Omega).
        $$
        Thus $\{T^{\frac{\sigma+1}{2}}_k(u_n)\}_{n\in\mathbb{N}}$ is bounded in $\mathcal{H}^{1}_0(\Omega)$. By Proposition \ref{BC}, for all $K\Subset\Omega$ there exists $\sigma_K>0$ such that $u(x)\geq u_1(x)\geq \sigma_K>0$ for a.e.\ $x\in K$. Consequently, the sequence $ \{T_k(u_n)\}_{n \in \mathbb{N}} $ is uniformly bounded in $H^1_{\text{loc}}(\Omega)$.
        Therefore, up to a subsequence, there exists $u\in \mathcal{H}^{1}_{loc}(\Omega)$ such that
        $T^{\frac{{\sigma}+1}{2}}_k(u), G_k(u)\in \mathcal{H}^{1}_{0}(\Omega)$,
        with $G_k(u_n)\rightharpoonup G_k(u)$ weakly in $\mathcal{H}^{1}_{0}(\Omega)$,
        $T^{\frac{{\sigma}+1}{2}}_k(u_n) \rightharpoonup T^{\frac{{\sigma}+1}{2}}_k(u)$ weakly in $\mathcal{H}^{1}_{0}(\Omega)$,
        and $u_n\uparrow u$ strongly in $L^r(\Omega)$ for all $1\le r< 2^*$.
        \\
        Testing with $\varphi \in\mathcal{C}^1_c(\Omega)$ in \eqref{apro} yields
        \begin{equation}\label{lastt}
            \dyle\int_{\Omega}\nabla u_n \nabla \varphi\, dx +\frac{a_{N,s}}{2}\iint_{\mathbb R^N}{\frac{(u_n(x)-u_n(y))(\varphi(x)-\varphi(y))}{|x-y|^{N+2s}}}\,dx\,dy=\lambda \int_\Omega
            \dfrac{u_n\varphi }{|x|^{2}}\,dx+\int_\Omega\dfrac{\varphi
                h_n(x)}{(u_n+\frac{1}{n})^\sigma}\,dx.
        \end{equation}
        By the above estimates, letting $n\to \infty$, we deduce that
        $$
        \lambda \int_\Omega \dfrac{u_n\varphi
        }{|x|^{2}}\,dx+\int_\Omega\dfrac{\varphi
            h_n(x)}{(u_n+\frac{1}{n})^{\sigma}}\,dx\to \lambda \int_\Omega \dfrac{u\varphi
        }{|x|^{2}}\,dx+\int_\Omega\dfrac{\varphi h(x)}{u^{\sigma}}\,dx<+\infty.
        $$
        Hence, passing to the limit in \eqref{lastt}, we obtain
        $$
        \dyle\int_{\Omega}\nabla u \nabla \varphi\, dx +\frac{a_{N,s}}{2}\iint_{\mathbb R^N}{\frac{(u(x)-u(y))(\varphi(x)-\varphi(y))}{|x-y|^{N+2s}}}\,dx\,dy=\lambda \int_\Omega
        \dfrac{u\varphi }{|x|^{2}}\,dx+\int_\Omega\dfrac{\varphi
            h(x)}{u^\sigma}\,dx,
        $$
        so $u$ is a weak solution to Problem~\eqref{mainp}.
        \\
        In this case, since $G_k(u)\in \mathcal{H}^{1}_{0}(\Omega)$, choosing $k=\varepsilon>0$ gives $G_\varepsilon(u)=(u-\varepsilon)^+\in \mathcal{H}^{1}_{0}(\Omega)$, that is, $u=0$ on $\partial \Omega$ in the sense of Definition \ref{rediri}.
        \smallbreak
        \noindent Finally, if we take $\sigma$ such that
        \[\left[\dfrac{4\sigma}{(\sigma+1)^2}-\dfrac{\lambda}{\Lambda_{N}}\right]>0\]
        and using $u_n^{\sigma}$ as a test function in \eqref{apro}, we obtain
        \begin{eqnarray*}
            \int_{\Omega}\nabla u_n \nabla u_n^{\sigma}\,dx+\dfrac{a_{N,s}}{2}\iint_{\mathbb R^N}\dfrac{(u_n(x)-u_n(y))(u^{\sigma}_n(x)-u^{\sigma}_n(y))}{|x-y|^{N+2s}}\,dx\,dy\le
            \lambda \int_{\Omega}\dfrac{u^{\sigma+1}_n}{|x|^{2}}\,dx+
            \int_{\Omega} h_n\, dx.
        \end{eqnarray*}
        By Lemma \ref{alg} and Hardy's inequality, we obtain
        $$
        \left(\dfrac{4\sigma}{(\sigma+1)^2}-\frac{\lambda}{\Lambda_{N}}\right)\left[\int_\Omega\left|\nabla u_n^{\frac{\sigma+1}{2}}\right|^2 dx+\dfrac{a_{N,s}}{2}\iint_{D_\Omega}\dfrac{\left(u^{\frac{\sigma+1}{2}}_n(x)-u^{\frac{\sigma+1}{2}}_n(y)\right)^2}{|x-y|^{N+2s}}
        \,dx\,dy\right]
        \le \|h\|_{L^1(\Omega)}.$$
        Therefore $u^{\frac{\sigma+1}{2}}\in \mathcal{H}^{1}_{0}(\Omega)$, which completes the proof.
    \end{proof}

    \begin{remark}\label{RRU}

 We point out here three further properties of the solution: additional regularity in the critical case $\lambda=\Lambda_N$, uniqueness under suitable assumptions, and minimality when the solution does not belong to $\mathcal{H}^1_0(\Omega)$.
        \begin{itemize}
            \item[(i)] If $0<\sigma\leq 1$ and $\lambda=\Lambda_N$, then arguing as in Case~1 of the proof of Theorem \ref{m1}, using the improved Hardy inequality \eqref{imp00}, we deduce that the sequence $\{u_n\}_{n \in \mathbb{N}}$ is bounded in $W^{1,q}_0(\Omega)$ for all $q<2$. Since $\{u_n\}_{n \in \mathbb{N}}$ is monotone in $n$, the compactness result of Lemma \ref{compa} gives $u_n \to u$ strongly in $W^{1,q}_0(\Omega)$ for all $q<2$. It is not difficult to verify that $u\in W^{1,q}_0(\Omega)$ satisfies
            \begin{equation*}
                \dyle \int_{\Omega} \nabla u\nabla \varphi\,dx
                +\frac{a_{N,s}}{2}\int_{\mathbb{R}^N}\int_{\mathbb{R}^N}
                \frac{(u(x)-u(y))(\varphi(x)-\varphi(y))}{|x-y|^{N+2s}}\,dx\,dy
                =\Lambda_N\,\dint_{\Omega}\frac{u\varphi}{|x|^2}\,dx+\int_{\Omega}
                \frac{f(x)}{u^\sigma}\,\varphi\, dx,
            \end{equation*}
            for every $\varphi \in C^\infty_c(\Omega)$.

            \item[(ii)] In the case $0<\sigma\leq 1$ and $\lambda<\Lambda_N$, the solution satisfies $u\in \mathcal{H}^1_0(\Omega)$. Applying Lemma~\ref{compar}, one can then prove that the solution is unique. More generally, whenever a solution belongs to $\mathcal{H}^1_0(\Omega)$, Lemma~\ref{compar} ensures uniqueness.

            \item[(iii)] When the solution does not belong to $\mathcal{H}^1_0(\Omega)$, one can show that the solution obtained through the above approximation procedure is minimal. $\square$
        \end{itemize}
    \end{remark}

    \subsubsection{\bf The case where $h$ is a measure}
    \

    We now consider the case where $h$ is a nonnegative Radon measure. In this situation, the monotonicity of the approximating sequence $\{u_n\}_n$ is no longer guaranteed, and the limit problem may fail to admit a solution. This phenomenon occurs in particular when the measure is singular with respect to the Sobolev capacity introduced below.

    For the reader's convenience, let us  recall the following definitions.

    \begin{Definition}\label{2cap}
        Let $K$ be a compact subset of $\Omega$. The \emph{capacity} of $K$ with respect to $\Omega$ is defined by
        \[
        \operatorname{cap}_{2,s}(K) = \inf \left\{ \|\varphi\|_{\mathcal{H}^1_0(\Omega)}^2 : \varphi \in C_0^\infty(\Omega), \, \varphi \geq \chi_K \right\},
        \]
        where $\chi_K$ denotes the characteristic function of $K$ (by convention, $\inf \varnothing = +\infty$).
        \\
        If $U \subset \Omega$ is open, we set
        \[
        \operatorname{cap}_{2,s}(U) = \sup \left\{ \operatorname{cap}_{2,s}(K) : K \text{ compact}, \; K \subset U \right\},
        \]
        and for an arbitrary Borel set $B \subset \Omega$,
        \[
        \operatorname{cap}_{2,s}(B) = \inf \left\{ \operatorname{cap}_{2,s}(U) : U \text{ open}, \; B \subset U \right\}.
        \]
    \end{Definition}
    \smallbreak
    In view of the definition of $\mathcal{H}^1_0(\Omega)$, $\operatorname{cap}_{2,s}$ is equivalent to the classical Sobolev capacity obtained by replacing $\mathcal{H}_0^1(\Omega)$ with $H^1_0(\Omega)$ in Definition \ref{2cap}.
    \\
    Let $\mu$ be a Radon measure on $\Omega$. We say that $\mu$ is \emph{concentrated} on a Borel set $E\subseteq \Omega$ if $\mu(B) = \mu(B\cap E)$ for every Borel set $B\subseteq\Omega$.

    \begin{Definition}\label{singularMM}
        A positive Radon measure $\mu$ on $\Omega$ is said to be \emph{singular} if it is concentrated on a subset $E \subset \Omega$ such that $\operatorname{cap}_{2,s}(E)=0$.
        We denote the set of positive singular measures by $\mathcal{M}_s$.
    \end{Definition}

    We are now in a position to state the following nonexistence result. For simplicity, we give the proof for $\sigma=1$. The other cases are analogous, see Remark \ref{EXT00}".

       \begin{Theorem}\label{sing}
        Assume that $\lambda<\Lambda_N$. Let $u_n$ be the unique solution to the problem
        \begin{equation}\label{apr01pp}
            \left\{\begin{array}{rclll}
                \mathcal{L}_s(u_n) & = & \lambda\dfrac{u_n}{|x|^2}+\dfrac{h_n}{u_n} & \text{\rm in } & \Omega, \\[2mm]
                u_n & > & 0 & \text{\rm in } & \Omega, \\[2mm]
                u_n & = & 0 & \text{\rm in } & \mathbb{R}^N\setminus\Omega,
            \end{array}\right.
        \end{equation}
        where $\{h_n\}_n\subset L^\infty(\Omega)$ is a sequence of nonnegative functions such that $\|h_n\|_{L^1(\Omega)}\le C$. Assume that $h_n\rightharpoonup \mu$ in the narrow topology of measures, where $\mu$ is a nonnegative singular Radon measure, {concentrated on a closed set $E\subset\Omega$}. Then
        \[
        u_n \rightharpoonup 0 \quad \text{weakly in } \mathcal{H}^1_0(\Omega).
        \]
    \end{Theorem}

    \begin{proof}
        Using $u_n$ as a test function in \eqref{apr01pp} (by approximation), we deduce that $\{u_n\}_n$ is bounded in $\mathcal{H}_0^{1}(\Omega)$. Hence there exists a nonnegative function $u\in \mathcal{H}^1_0(\Omega)$ such that, up to a subsequence, $u_n\rightharpoonup u$ weakly in $\mathcal{H}^1_0(\Omega)$. To prove the theorem it suffices to show that $u=0$ a.e.\ in $\Omega$.
        \\
        Since $\mu$ is singular with respect to $\operatorname{cap}_{2,s}$, for any $\varepsilon>0$ there exist a compact set $K_\varepsilon\subset E$ and a function $\phi_\varepsilon\in C_0^{\infty}(\Omega)$ with $0\leq \phi_\varepsilon\leq 1$ and $\phi_\varepsilon\equiv 1$ on $K_\varepsilon$, such that
        \[
        \mu(E\setminus K_\varepsilon)<\varepsilon \qquad \text{and} \qquad \|\phi_\varepsilon\|^{2}_{\mathcal{H}^1_0(\Omega)}<\varepsilon.
        \]
        By Picone's inequality (Theorem \ref{Picone}), we obtain
        \[
        \int_{\Omega}\frac{\mathcal{L}_s(u_n)}{u_n}|\phi_\varepsilon|^{2} \,dx \le \|\phi_\varepsilon\|^{2}_{\mathcal{H}^1_0(\Omega)}\le \varepsilon.
        \]
        Thus $\dyle \int_{K_\varepsilon}\frac{h_n}{u_n^{2}}\,dx \le \varepsilon$. Applying H\"{o}lder's inequality yields
        \begin{align*}
            \int_{K_\varepsilon}\frac{h_n}{u_n}\,dx
            &\le \left(\int_{K_\varepsilon}\frac{h_n}{u_n^{2}}\,dx\right)^{\frac{1}{2}} \left(\int_{K_\varepsilon} h_n\,dx\right)^{\frac{1}{2}} \\
            &\le C\varepsilon^{\frac{1}{2}}.
        \end{align*}

        We claim that for every $\phi\in C^\infty_0(\Omega)$,
        \[
        \lim_{n\to \infty}\int_{\Omega}\frac{h_n}{u_n}\phi\,dx = 0.
        \]
        Indeed,
        \begin{align}\label{eq:1}
            \int_{\Omega}\frac{h_n}{u_n}\phi \,dx
            &\le \|\phi\|_\infty \int_{K_\varepsilon}\frac{h_n}{u_n}\,dx + \int_{\Omega\setminus K_\varepsilon}\frac{h_n}{u_n}\phi\,dx \\
            &\nonumber\le \|\phi\|_\infty \varepsilon^{\frac{1}{2}} + \int_{\Omega\setminus K_\varepsilon}\frac{h_n}{u_n}|\phi|\,dx.
        \end{align}
        By outer regularity of the measure, let $U_\varepsilon$ be an open set such that $K_\varepsilon\subset U_\varepsilon\Subset\Omega$ and $\mu(U_\varepsilon\setminus K_\varepsilon)<\varepsilon$.
        Since the measure is concentrated on $E$, we may choose the sequence $\{h_n\}_n$ such that $\operatorname{supp}(h_n)\subset U_\varepsilon$ for $n\ge n_0$. Hence
        \begin{equation}\label{eq:2}
            \int_{\Omega\setminus U_\varepsilon}
            \frac{h_n}{u_n}|\phi|\,dx=0.
        \end{equation}
        Since $\overline{U_\varepsilon}\Subset\Omega$, there exists $c_\varepsilon>0$, independent of $n$, such that
        $u_n\geq c_\varepsilon$ in $U_\varepsilon$.
        We note that, since the approximating sequence $\{h_n\}$ is not monotone, this lower bound cannot be derived from the maximum principle. Instead, it follows from Harnack's inequality (Theorem \ref{harnack}); see in particular Proposition \ref{Lm}. Indeed, $u_n$ is a supersolution satisfying
        \[
        \mathcal{L}_s(u_n) \geq \lambda\frac{u_n}{|x|^2} \quad \text{in } \Omega.
        \]
        Therefore,
        \begin{equation}\label{eq:3}
            \int_{U_\varepsilon\setminus K_\varepsilon}
            \frac{h_n}{u_n}|\phi|\,dx
            \leq
            \frac{1}{c_\varepsilon}
            \int_{U_\varepsilon\setminus K_\varepsilon}|\phi|h_n\,dx.
        \end{equation}
        Since $U_\varepsilon$ is open and $K_\varepsilon$ is compact, the set $U_\varepsilon \setminus K_\varepsilon$ is open. By outer regularity of $\mu$, there exists an open set $V_\varepsilon \supset U_\varepsilon \setminus K_\varepsilon$ such that $\mu(V_\varepsilon \setminus (U_\varepsilon \setminus K_\varepsilon)) < \varepsilon$.
        Define a continuous function $\psi_\varepsilon \in C_0(\Omega)$ satisfying $0 \le \psi_\varepsilon \le 1$, with $\psi_\varepsilon \equiv 1$ on $U_\varepsilon \setminus K_\varepsilon$ and $\operatorname{supp}(\psi_\varepsilon) \subset V_\varepsilon$. Since $\psi_\varepsilon |\phi| \in C_0(\Omega)$, applying the narrow convergence $h_n \rightharpoonup \mu$ yields
        \begin{align}\label{eq:4}
            \limsup_{n\to\infty} \frac{1}{c_\varepsilon} \int_{U_\varepsilon \setminus K_\varepsilon} |\phi| h_n \, dx
            &\le \limsup_{n\to\infty} \frac{1}{c_\varepsilon} \int_{\Omega} \psi_\varepsilon |\phi| h_n \, dx = \frac{1}{c_\varepsilon} \int_{\Omega} \psi_\varepsilon |\phi| \, d\mu \\
            &\nonumber\le \frac{\|\phi\|_\infty}{c_\varepsilon} \mu(V_\varepsilon) < \frac{\|\phi\|_\infty}{c_\varepsilon} \big(\mu(V_\varepsilon \setminus (U_\varepsilon \setminus K_\varepsilon))+ \mu(U_\varepsilon \setminus K_\varepsilon)  \big)\leq  \frac{2\varepsilon\|\phi\|_\infty}{c_\varepsilon}.
        \end{align}
        Combining \eqref{eq:1}--\eqref{eq:4}, we conclude that
        \[
        \lim_{n\to \infty}\int_{\Omega}\frac{h_n}{u_n}\phi \,dx = 0,
        \]
        which proves the claim. Recalling that $u_n\rightharpoonup u$ weakly in $\mathcal{H}^1_0(\Omega)$, we get
        \[
        \lim_{n\to \infty}\int_{\Omega}\mathcal{L}_s(u_n)\phi \,dx = \int_{\Omega}\mathcal{L}_s(u)\phi \,dx.
        \]
        Hence $u\in \mathcal{H}^1_0(\Omega)$ solves
        \begin{equation}\label{DDD}
            \mathcal{L}_s(u) = \lambda\frac{u}{|x|^2} \quad \text{in } \mathcal{D}'(\Omega).
        \end{equation}
        Since $\lambda<\Lambda_{N}$, it follows that $u\equiv 0$, which concludes the proof.
    \end{proof}

\begin{remark}
        Let $\mu=\mu_s+\theta$, where $\theta\in L^1(\Omega)$ is a nonnegative function and $\mu_s\in \mathcal{M}_s$ is a nonnegative singular measure. Let $\{\rho_n\}_n$ and $\{\theta_n\}_n$ be two sequences of bounded functions such that $\rho_n\rightharpoonup \mu_s$ in the narrow topology of measures and $\theta_n\to \theta$ strongly in $L^1(\Omega)$.
        If $u_n$ is the unique solution to Problem~\eqref{apr01pp} with $h_n=\rho_n+\theta_n$, then $u_n\rightharpoonup v$ weakly in $\mathcal{H}^1_0(\Omega)$, where $v$ is the unique solution to the problem
        \begin{equation*}
            \left\{\begin{array}{rclll}
                \mathcal{L}_s(v)&=&\lambda\dfrac{v}{|x|^2}+\dfrac{\theta}{v} &\hbox{\rm in } &\Omega,\\[2mm]
                v&>&0 &\hbox{\rm in } &\Omega, \\[2mm]
                v&=&0 &\hbox{\rm in } &\mathbb{R}^N\setminus\Omega.
            \end{array}\right.
        \end{equation*}
        In other words, the singular part $\mu_s$ of the measure has no effect on the limit problem.
    \end{remark}

\begin{remark}\label{EXT00}
Let $\{h_n\}_n\subset L^\infty(\Omega)$ be such that $(\lambda,\sigma,h_n)$ satisfy the
hypotheses of Theorem~\ref{m1}, and let $u_n$ be the solution of
\begin{equation}\label{apr01ppp1}
\left\{
\begin{array}{rclll}
\mathcal{L}_s(u_n) & = & \lambda\dfrac{u_n}{|x|^2}+\dfrac{h_n}{u_n^\sigma} & \text{in} & \Omega,\\
u_n & > & 0 & \text{in} & \Omega,\\
u_n & = & 0 & \text{in} & \mathbb{R}^N\setminus\Omega,
\end{array}
\right.
\end{equation}
given by Theorem~\ref{m1}. Assume further that $h_n\rightharpoonup\mu$ narrowly, where $\mu$
is a nonnegative Radon measure, singular with respect to the Lebesgue measure and
concentrated on a closed set $E\subset\Omega$. Then:
\begin{enumerate}
\item[(i)] if $0<\sigma<1$ and $\lambda<\frac{4\Lambda_N}{(\sigma+1)^2}$, then
$u_n\rightharpoonup 0$ weakly in $\mathcal{H}^1_0(\Omega)$;
\item[(ii)] if $\sigma>1$ and $\lambda<\Lambda_N$, then, for every $k>0$, both $G_k(u_n)$ and
$\big(T_k(u_n)\big)^{\frac{\sigma+1}{2}}$ converge weakly to $0$ in $\mathcal{H}^1_0(\Omega)$.
If, in addition, $\lambda<\frac{4\sigma\Lambda_N}{(\sigma+1)^2}$, then
$u_n^{\frac{\sigma+1}{2}}\rightharpoonup 0$ weakly in $\mathcal{H}^1_0(\Omega)$.
\end{enumerate}
The proof follows the same arguments as those used in the preceding theorem. More precisely, we choose suitable test functions adapted to the hypotheses on $\lambda$ and $\sigma$ in order to establish the required uniform estimates for the sequence $\{u_n\}_n$. The asserted weak convergence results then follow by applying the same compactness and limiting arguments as in the proof of the preceding theorem.
\end{remark}

\section{Perspectives and open problems}\label{sec:last}

    The results obtained in this work naturally lead to several open problems
    and directions for future research.

    \subsection{Existence theory for general power-type nonlinearities}
    \

    The question of existence for Problem~\eqref{eq:problem1} with more general data will be investigated in a forthcoming work. Theorem~\ref{exist} provides existence for $1<q<q^+(\lambda)$ when $g=\alpha g_0$ with $g_0$ bounded and $\alpha$ small. A natural question is to determine, for unbounded data, the optimal relation between the summability of $g_0$ and the exponent $q$ ensuring existence. We expect such results to rely on suitable weighted Sobolev inequalities adapted to the mixed local--nonlocal framework, with particular attention to the interplay between the Hardy potential and the regularizing effects of the operator.

    \subsection{Gradient-dependent and nonlocal-gradient nonlinearities}
    \

        A natural direction, still largely unexplored in the mixed setting with a Hardy potential, concerns problems in which the nonlinear term also depends on the gradient of the solution or on a nonlocal counterpart of it. Specifically, one can consider
        \[
        \left\{\begin{array}{rcll}
            \mathcal{L}_s(u)&=&\lambda \dfrac{u}{|x|^{2}}
            + H\!\left(x,u,\nabla u,\mathcal{D}^s u\right)
            &\hbox{in } \Omega,\\[4pt]
            u&=&0 &\hbox{in }\mathbb{R}^N\setminus\Omega,
        \end{array}\right.
        \]
        where $H$ is a Hamiltonian with superlinear growth in the gradient variables, and $\mathcal{D}^s u$ denotes a suitable nonlocal differential quantity, such as the half-power $(-\Delta)^{s/2}u$, the Riesz fractional gradient, or the difference quotient $(u(x)-u(y))/|x-y|^{(N+2s)/2}$. The presence of the Hardy singularity introduces a genuine difficulty: the coercivity of the quadratic form associated with $\mathcal{L}_s-\lambda|x|^{-2}$ deteriorates as $\lambda$ approaches $\Lambda_N$, and solutions lose regularity near the origin. The interplay between the Hardy weight and the gradient growth therefore requires new weighted Calder\'on--Zygmund-type estimates for $\mathcal{L}_s$. For the purely local operator $-\Delta$, problems of this type have been extensively studied; see, e.g., \cite{BBM} and the references therein. Extending these results to the mixed setting with a Hardy potential, and identifying the growth conditions on $H$ that ensure existence, uniqueness or blow-up, constitutes a challenging open problem.

\subsection{The quasilinear case}
    \

    Another significant extension concerns the quasilinear analogue of the operator $\mathcal{L}_s$, namely
    \[
    \mathcal{L}_{s,p}(u) := -\Delta_p u + (-\Delta)^s_p u,
    \]
    where $1<p<N$, $-\Delta_p u = -\operatorname{div}(|\nabla u|^{p-2}\nabla u)$ is the $p$-Laplacian and $(-\Delta)^s_p$ is the fractional $p$-Laplacian, $s\in(0,1)$  (see, e.g., \cite{IMS, CKP1}).\\
    When $1<p<N$, the associated Hardy inequality reads
    \[
    \int_\Omega |\nabla u|^p\,dx
    + c_{N,s,p}\iint_{\mathbb{R}^N\times\mathbb{R}^N}
    \frac{|u(x)-u(y)|^p}{|x-y|^{N+sp}}\,dx\,dy
    \geq \Lambda_{N,p}\int_\Omega\frac{|u|^p}{|x|^p}\,dx,
    \]
    where $\Lambda_{N,p}=\bigl(\frac{N-p}{p}\bigr)^p$ is the sharp Hardy
    constant in the mixed $p$-Laplacian framework.

    The analysis of the problem
    \begin{equation*}
        \left\{\begin{array}{rcll}
            \mathcal{L}_{s,p}(u) &=& \lambda\dfrac{u^{p-1}}{|x|^p}
            + f(x,u) &\hbox{in }\Omega,\\[3pt]
            u &>& 0 &\hbox{in }\Omega,\\[3pt]
            u &=& 0 &\hbox{in }\mathbb{R}^N\setminus\Omega,
        \end{array}\right.
    \end{equation*}
    raises several new difficulties. In the quasilinear setting, ground-state representations of Frank--Seiringer type only provide inequalities, not identities, see e.g. \cite{FS}. So the change of variables $w=|x|^{\gamma_1}u$, which reduces Problem~\eqref{tlem} to a weighted operator without Hardy singularity, is no longer exact. Determining the behavior of supersolutions near the origin and the critical exponent governing existence therefore requires new arguments and refined estimates to compensate for the lack of linearity of the operator.

    \subsection{The parabolic problem and its extensions}
    \

    The parabolic counterpart of Problem~\eqref{tlem}, namely
    \begin{equation*}
        \left\{\begin{array}{rcll}
            \partial_t u + \mathcal{L}_s(u) &=& \lambda\dfrac{u}{|x|^2}
            + f(x,t,u) &\hbox{in }\Omega\times(0,T),\\[3pt]
            u &\geq& 0 &\hbox{in }\Omega\times(0,T),\\[2pt]
            u &=& 0 &\hbox{in }(\mathbb{R}^N\setminus\Omega)\times(0,T),\\[2pt]
            u(\cdot,0) &=& u_0\geq 0 &\hbox{in }\Omega,
        \end{array}\right.
    \end{equation*}
    as well as its variant with a gradient-dependent term $H(x,u,\nabla u,\mathcal{D}^s u)$, raises questions that remain largely open
    for the mixed operator $\mathcal{L}_s$. A first natural question is to obtain a counterpart of the result of Baras and Goldstein \cite{BG} for the heat equation, namely to show that $\lambda=\Lambda_N$ is the sharp threshold separating the existence of a positive local solution from instantaneous and complete blow-up. A second direction concerns the existence of a critical exponent separating global existence from finite-time blow-up, both for a power source term $f(x,t,u)=u^q$
    and for a gradient term $H=|\nabla u|^r$. In the spirit of the classical Fujita phenomenon \cite{Fujita}, in particular when $\Omega$ is replaced by $\mathbb{R}^N$, the interaction between the Hardy potential and the mixed diffusion is expected to shift the classical thresholds. Finally, the asymptotic behavior of global solutions, such as decay rates and convergence to stationary solutions, and its dependence on $\lambda$, deserves further investigation.

\vspace{1cm}

\

\begin{center}
{\sc Data availability statement}

\

All data generated or analyzed during this study are included in this published article.
\end{center}

\

\begin{center}
	{\sc Conflict of interest statement}
	
	\
	
	The authors declare that they have no competing interest.
\end{center}

\end{document}